\documentclass[preprint,12pt]{elsarticle}

\usepackage{amsmath}
\usepackage{amssymb}
\usepackage{amsthm}
\usepackage{mathtools}
\usepackage{bm}
\usepackage{algorithm}
\usepackage[linkcolor=blue,colorlinks=true]{hyperref}
\usepackage{graphicx}
\usepackage{comment}
\usepackage{arydshln}
\usepackage{epsfig}
\usepackage{verbatim} 
\usepackage{tikz}
\usepackage[export]{adjustbox}
\usepackage{mathrsfs}
\usepackage{mathtools}
\usepackage{todonotes}

\theoremstyle{definition}
\newtheorem{definition}{Definition}[section]
\newtheorem{theorem}{Theorem}[section]
\newtheorem{lemma}{Lemma}[section]

\newtheorem{problem}{Problem}[section]
\newtheorem{remark}{Remark}[section]

\newtheorem{assumption}{Assumption}[section]

\newtheorem*{Problem*}{Problem}
\DeclareMathOperator*{\wrap}{wrap}
\newcommand{\bb}{\boldsymbol}
\newcommand{\reals}{\mathbb{R}}
\newcommand{\s}{\mathcal{S}}
\newcommand{\classK}{\mbox{class-$\K$}}
\newcommand{\K}{\mathcal{K}}
\DeclareMathOperator{\^T}{^\top}
\DeclareMathOperator*{\argmax}{arg\,max}
\DeclareMathOperator*{\argmin}{arg\,min}
\newcommand{\highdegree}{\mathcal{G}}

\newcommand{\pth}[1]{\left(#1\right)}
\newcommand{\bkt}[1]{\left[#1\right]}
\newcommand{\rarr}{\rightarrow}
\newcommand{\R}{\mathbb{R}}
\newcommand{\N}{\mathcal{N}}
\newcommand{\A}{\mathcal{A}}
\newcommand{\eps}{\epsilon}
\newcommand{\brc}[1]{\left \{#1\right \}}
\newcommand{\bmx}[1]{\begin{bmatrix}#1\end{bmatrix}}

\journal{Annual Reviews on Control}

\begin{document}

\begin{frontmatter}



\title{Advances in the Theory of Control Barrier Functions: Addressing Practical Challenges in Safe Control Synthesis for Autonomous and Robotic Systems}

\author[KG]{Kunal Garg}
\author[JU]{James Usevitch}
\author[UM2]{Joseph Breeden}
\author[MB]{Mitchell Black}
\author[UM2]{Devansh Agrawal}
\author[UM1]{Hardik Parwana}
\author[UM0]{Dimitra Panagou\corref{cor1}}
\cortext[cor1]{Corresponding author}
\ead{dpanagou@umich.edu}

\affiliation[KG]{organization={Department of Aeronautics and Astronautics, Massachusetts Institute of Technology},
            addressline={77 Massachusetts Avenue}, 
            city={Cambridge},
            postcode={02139}, 
            state={MA},
            country={USA}}

\affiliation[JU]{organization={Department of Electrical and Computer Engineering, Brigham Young University},
            addressline={450 Engineering Building}, 
            city={Provo},
            postcode={84602}, 
            state={UT},
            country={USA}}

\affiliation[UM2]{organization={Department of Aerospace Engineering, University of Michigan, Ann Arbor},
            addressline={1320 Beal Avenue}, 
            city={Ann Arbor},
            postcode={48109}, 
            state={MI},
            country={USA}}  
            
\affiliation[MB]{organization={Toyota North America Research and Development},
            addressline={1555 Woodridge Ave}, 
            city={Ann Arbor},
            postcode={48105}, 
            state={MI},
            country={USA}}

\affiliation[UM1]{organization={Department of Robotics, University of Michigan, Ann Arbor},
            addressline={2505 Hayward St.}, 
            city={Ann Arbor},
            postcode={48109}, 
            state={MI},
            country={USA}}

\affiliation[UM0]{organization={Department of Robotics and Department of Aerospace Engineering, University of Michigan, Ann Arbor},
            addressline={2505 Hayward St.}, 
            city={Ann Arbor},
            postcode={48109}, 
            state={MI},
            country={USA}}
            
\begin{abstract}
This tutorial paper presents recent work of the authors that extends the theory of Control Barrier Functions (CBFs) to address practical challenges in the synthesis of safe controllers for autonomous systems and robots. We present novel CBFs and methods that handle safety constraints (i) with time and input constraints under disturbances, (ii) with high-relative degree under disturbances and input constraints, and (iii) that are affected by adversarial inputs and sampled-data effects. We then present novel CBFs and adaptation methods that prevent loss of validity of the CBF, as well as methods to tune the parameters of the CBF online to reduce conservatism in the system response. We also address the pointwise-only optimal character of CBF-induced control inputs by introducing a CBF formulation that accounts for future trajectories, as well as implementation challenges such as how to preserve safety when using output feedback control and zero-order-hold control. Finally we consider how to synthesize non-smooth CBFs when discontinuous inputs and multiple constraints are present. 
\end{abstract}


\end{frontmatter}


\section{Introduction}
\label{Introduction}

Control Barrier Functions (CBFs) have been developed in recent years as a tool to verify and synthesize trajectories for nonlinear constrained control systems. Their principle is as follows: Given a constraint function, termed barrier function thereafter, whose zero super-level (or sub-level) sets define a constrained set, termed also the safe set, the key idea is that one restricts the rate of change of the barrier function along the system trajectories using a class $\mathcal K$ function of the barrier function~\citep{ames2017control, ames2019control}. If such a condition can be satisfied everywhere in the constrained set under the given dynamics and control input constraints, then the barrier function is called a Control Barrier Function (CBF), and the constrained set is forward invariant. This method, in conjunction with Control Lyapunov functions (CLFs) for stability, has been employed to design safe controllers for several applications.

\subsection{Challenges}
However, similarly to Lyapunov methods, some of the major challenges of verifying safety and synthesizing safe controllers are that 1) finding a valid CBF for arbitrary system dynamics is not trivial, 2) safety constraints of high relative degree compared to the system dynamics, as well as input constraints make the problem even more challenging, 3) modeling/parametric uncertainty dictates the formulation of tools and techniques from robust control, adaptation and learning in order to define valid CBFs that account for uncertainty. Among other challenges, the fact that the control inputs derived due to the CBF condition are only \emph{pointwise} optimal (also often called myopic control inputs), has given rise to considerations on under which conditions one can guarantee optimality and feasibility of the resulting control policies.  

In the recent 3-4 years, the literature has seen an abundance of papers with a variety of methodologies that aim to address some of the aforementioned challenges. Finding valid CBFs for example has been addressed with offline~\cite{dawson2023safe,jagtap2020control} and online methods, searching for either a valid function $h$ over the constrained set, searching for some of the parameters of the CBF condition~\cite{robey2020learning, lindemann2021learning}, or adapting for those parameters online in order to render the candidate function a valid CBF~\cite{xiao2021adaptive}. High-relative degree constraint functions have been first addressed in ~\cite{nguyen2016exponential}, which considers the class of Exponential CBFs when the $\classK$ functions used in CBF derivative condition are linear in their argument; then~\cite{hocbf} generalizes Exponential CBFs to generic nonlinear $\classK$ functions in the form of Higher-Order CBFs (HOCBF). Time constraints and specifications (beyond state constraints) and cooperative multi-agent systems have also been considered~\cite{lindemann2020barrier, garg2022fixed}. In a relatively less explored area, CBFs for noncooperative multi-agent systems have also started being studied recently~\cite{usevitch2021adversarial}. Adaptive, robust, and learning-based formulations have also appeared in order to deal with various sources of uncertainty (stochastic uncertainty in the system dynamics, parametric uncertainty, deterministic additive external disturbances), see for example~\cite{taylor2020acc, jankovic2018robust, dhiman2023tac}.  

\subsection{Overview and Organization}
The scope of this tutorial paper is not to provide a thorough literature survey and review of recent CBF techniques, but rather to focus on some of the authors' own work on safety verification and control, presented in a roughly chronological and thematic order. More specifically, Section \ref{section::FxTCBF} shows how time constraints can be encoded as \emph{novel forms of timed CBFs}, called Fixed-Time Barriers, how to concurrently handle time, safety and input constraints using novel forms of FxT-CLF-CBF-QPs, and how such concepts can be used to solve problems ranging from spatiotemporal control, to integrated planning and control with safety and recursive feasibility guarantees. Section \ref{section::HOCBF} presents our \emph{constructive} methods for constraints with high-relative degree under disturbances and input constraints. We also introduce Input-Constrained CBFs, which are generalizations of High-Order CBFs. Then, Section \ref{section::adaptation} presents novel adaptation methods so that either the control-input coefficient is tuned online to prevent loss of controllability, or the parameters of the CBF condition are tuned online in order to reduce conservatism in the system response. Then, Section \ref{section::predictive} addresses the pointwise optimal character of CBF-induced control inputs by accounting for future trajectories, in a computationally-efficient way that checks for possible future safety violations, and adjusts the control action as needed. In Section \ref{section::implementation} we address implementation challenges such as how to preserve safety when using output feedback control and zero-order hold control, while Section \ref{section::adversarial} covers the definition of Adversarially-Robust CBFs for multi-robot control. Finally, we present our approach on how to synthesize non-smooth CBFs when multiple constraints are present in \ref{section::stronginvariance}. Concluding, we note some of our more recent and ongoing work in Section \ref{section::conclusions}. Again, while we have cited relevant work of many of our fellow colleagues in the field, the references list is vastly incomplete. It is out of the scope of this paper to provide a thorough literature review. Interested readers are referred to~\cite{ames2019control, review2021} for recent comprehensive reviews on various topics related to CBFs, as well as to the survey papers in this special issue.

\section{Preliminaries: Definition of Control Barrier Functions, Set Invariance, and Basic Quadratic Program for Safe Control}
\subsection{Notations}
\label{section::notation}
The set of real numbers is denoted as $\mathbb R$ and the non-negative real numbers as $\mathbb R^+$. Given $x\in \mathbb R$, $y \in \mathbb R^{n_i}$, and $z\in \mathbb R^{n_i\times m_i}$, $|x|$ denotes the absolute value of $x$ and $||y||$ denotes $L_2$ norm of $y$. The interior and boundary of a set $\mathcal C$ are denoted by $\textrm{Int}(\mathcal C)$ and $\partial \mathcal C$. The distance of a point $x$ from a set $\mathcal C$ is denoted $|x|_{\mathcal C} = \min_{y \in \mathcal C} || x - y ||$. For $a\in \mathbb R^+$, a continuous function $\alpha:[0,a)\rightarrow[0,\infty)$ is a class $\mathcal K$ function if it is strictly increasing and $\alpha(0)=0$. A continuous function $\alpha:(-b,a)\rightarrow (-\infty,\infty)$ for $a,b\in \mathbb R^+$ is an extended class $\mathcal K$ function if it is strictly increasing and $\alpha(0)=0$. Furthermore, if $a=\infty$ and $\lim_{r\rightarrow \infty} \alpha(r)=\infty$, then it is called extended class-$\mathcal{K}_\infty$. The $k^{th}$ time derivative of a function $h(t,x):\mathbb R^{+}\times \mathbb R^{n}\rightarrow \mathbb R$ is denoted as $h^{(k)}$. 
For brevity, we will refrain from mentioning explicit arguments whenever the context is clear. For example, $h(x)$ may simply be denoted as $h$. The Lie derivative of a function $h$ w.r.t a function $f$ is denoted as $L_f h = \frac{\partial h}{\partial x}f$. 

\subsection{Control Barrier Functions}
Consider the nonlinear control-affine dynamical system:
\begin{align}
\dot{x} = f(x) + g(x)u,
\label{eq::dynamics_general}
\end{align}
where $x\in \mathcal X \subset \mathbb R^{n}$ and $u \in \mathcal U \subset \mathbb R^{m}$ represent the state and control input, and $f:\mathcal X\rightarrow \mathbb R^{n}$ and $g:\mathcal X\rightarrow \mathbb R^{m}$ are locally Lipschitz continuous functions. The set $\mathcal S(t)$ of allowable states at time $t$ is specified as an intersection of $N$ sets $\mathcal S_i(t),i\in\{1,2,..,N\}$, each of which is defined as the zero-superlevel set\footnote{Note that in certain sections of the current paper, as well as in many references in the related literature, the constrained set $\mathcal S_i(t)$ is defined as the zero-sublevel set of a constraint function $h_i$.} of a (sufficiently smooth) function $h_i:\mathbb R^+ \times \mathcal{X} \rightarrow \mathbb R$ as:
\begin{subequations}
\label{eq::safeset}
\begin{align}
        \mathcal S_i(t) &= \{ x \in \mathcal X \;|\; h_i(t,x) \geq 0 \}, \label{eq::safeset1} \\
        \partial \mathcal S_i(t) &= \{ x\in \mathcal X \;|\; h_i(t,x)=0 \}, \label{eq::safeset2} \\
        \textrm{Int} (\mathcal S_i)(t) &= \{ x \in \mathcal X \;|\; h_i(t,x)>0 \}. \label{eq::safeset3}
\end{align}
\end{subequations}

\begin{definition}
\label{definition::cbf_definition}
(Control Barrier Function)\citep{ames2017control} 
For the dynamical system (\ref{eq::dynamics_general}), $h_i: \mathbb R^+ \times \mathcal X\rightarrow \mathbb R$ is a control barrier function (CBF) on the set $\mathcal S_i(t)$ defined by (\ref{eq::safeset1})-(\ref{eq::safeset3}) for $t\geq 0$ if there exists a class-$\mathcal{K}$ function $\alpha_i: \mathbb R \rightarrow \mathbb R^+$ 
such that 
\begin{align}
 \sup_{u\in \mathcal{U}} \left[ \frac{\partial h_i(t,x)}{\partial t} + L_f h_i(t,x)+ L_g h_i(t,x)u \right] \geq -\alpha_i(h_i(t,x)) \nonumber \\
 \forall x\in \mathcal{S}_i, \forall t>0.
  \label{eq::cbf_derivative}
\end{align}
\end{definition}

Henceforth, we refer to~\eqref{eq::cbf_derivative} as the \textit{CBF derivative condition}.

\begin{theorem} (Set Invariance)~\citep{lindemann2018control}\label{thm: CBF safety}
Given the dynamical system (\ref{eq::dynamics_general}) and a set $\mathcal S_i(t)$ defined by (\ref{eq::safeset1})-(\ref{eq::safeset3}) for some continuously differentiable function $h_i:\mathbb R^+ \times \mathbb R^n\rightarrow \mathbb R$, if $h_i$ is a control barrier function on the set $\mathcal S_i(t)$, and there exists a $u:\mathbb R_+\times\mathbb R^n\to \mathbb R^m$, piecewise continuous in $t$ Lipschitz continuous in $x$, that satisfies
\begin{equation}
    \frac{\partial h_i(t,x)}{\partial t} + L_f h_i(t,x)+ L_g h_i(t,x)u \geq -\alpha_i(h_i(t,x)) ,\;
 \forall x\in \mathcal{S}_i(t), \forall t>0, \label{eq:cbf_online_condition}
\end{equation}\color{black}
then $\mathcal S_i(t)$ is forward invariant.
\end{theorem}

If $L_g h_i(t,x) \equiv 0$ $\forall (t,x)\in\mathbb R^+\times \mathcal X$, then the control input $u$ does not appear in the left-hand side of the CBF condition~\eqref{eq:cbf_online_condition}. Suppose the relative degree of the function $h_i$ w.r.t. the control input $u$ under the dynamics~\eqref{eq::dynamics_general} is equal to $r_i\geq 2$. We can then define $r_i$ functions as follows:
\begin{subequations}
    \begin{align}
   \psi_i^0(t,x) &= h_i(t,x), \\
   \psi_i^k(t,x) &= \dot \psi_i^{k-1}(t,x) + \alpha_{i}^k (\psi^{k-1}_i(t,x)), \; k\in \{1,2,\dots,r_i-1\} \label{eq::hocbf_derived_barriers},
\end{align}
\label{eq::hobcf_barriers}
\end{subequations}
and denote their zero-superlevel sets respectively, as: 
\begin{align}
  \label{eq::hobcf_intersection}
  \mathscr{C}_i(t) = \{x \; | \; \psi_i^k(t,x)\geq 0, \quad \forall k\in \{0,..,r_i-1\} \}.
\end{align}

\begin{definition}\label{definition::hocbf_definition}(Higher-Order CBF)\footnote{Definitions \ref{definition::cbf_definition} and \ref{definition::hocbf_definition} were presented in their original papers for the time-invariant safe sets $\mathcal S_i$. We note that an extension to the time-varying case can be proven with Nagumo's theorem applied to non-autonomous systems~\citep[Theorem 3.5.2]{carja2007viability} and hence we directly present that. This follows also the notation in~\citep{lindemann2018control}. }\citep{tan2021high}
The function $h_i(t,x):\mathbb R^+\times \mathbb R^{n}\rightarrow \mathbb R$ is a Higher-Order CBF (HOCBF) of $r_i$-th order on the set $\mathscr{C}_i(t)$ if there exist $r_i$ extended class-$\mathcal K$ functions $\alpha_i^k:\mathbb R \rightarrow \mathbb R, \;  k\in\{1,2,..,r_i\}$, and an open set $\mathscr{D}_i(t)\subset \mathbb R^+ \times \mathbb R^n$ with $\mathscr{C}_i(t) \subset \mathscr{D}_i(t) \subset \mathcal{X}$ such that
\begin{align}
      \dot \psi_i^{r_i-1}(t,x,u) \geq -\alpha_i^{r_i} (\psi_i^{r_i-1}(t,x)), ~ \forall x \in \mathscr{D}_i(t), \forall t\geq 0.
      \label{eq::hobcf_cbf_derivatrive}
    \end{align}
\end{definition}

Enforcing multiple constraints encoded via Control Lyapunov Functions~\cite{ames2017control} and (in general, high-order) Control Barrier Functions has commonly been addressed via the following class of controllers:

\noindent \textbf{(CLF-(HO)CBF-QP)}
\begin{subequations}
        \begin{align}
     (u, \delta) = \arg \min_{u\in \mathcal{U}, \delta \geq 0} \quad &  ||u-u_r(t,x)|| + M \delta^2 \\
         \textrm{s.t.} \quad & \dot{V}(t,x,u) \leq - k V(t,x) + \delta, \label{eq::qp_hoclf_constraint} \\
          & \dot \psi_i^{r_i-1}(t,x,u) \geq -\alpha_i^{r_i} (\psi_i^{r_i-1}(t,x)), \label{eq::qp_hocbf_constraint} \; \\ & \quad \quad \quad \quad  \quad \quad \quad \quad \quad 
 i\in \{1,2,..,N\} \nonumber
        \end{align}
        \label{eq::HOCBF-CLF-QP}
\end{subequations}
where $u_r: \mathbb R^+\times \mathbb R^n\rightarrow \mathbb R^m$ is the reference control input, often designed without any regard to constraints, $M\in \mathbb R^{+}$ is positive definite weighting matrix, $V(t,x)$ a control Lyapunov function (CLF) encoding convergence objectives for the system trajectories, $k \in \mathbb{R}^+$ is the exponential rate of convergence, and $\delta \in \mathbb R^+$ is a slack variable used to relax the CLF constraint~\eqref{eq::qp_hoclf_constraint}. The optimization~\eqref{eq::HOCBF-CLF-QP} is a QP when the dynamics is control-affine as in~\eqref{eq::dynamics_general} and $\mathcal{U}$ can be expressed in the form of a polytope $Au\leq b, A\in \mathbb R^{q\times m}, b\in \mathbb R^{q\times 1}, q>0$.


\section{Fixed-Time Control Barrier Functions: Synthesis under Time, Input and Safety Constraints}
\label{section::FxTCBF}

In this section, we present a method to address temporal constraints (e.g., convergence to a goal region within a given time horizon) in addition to safety constraints (realized via CBFs) for nonlinear systems with bounded inputs. The main references for this section are~\cite{garg2019control,garg2021characterization,garg2021advances,garg2022fixed,garg2021robust}.


\subsection{Fixed-time Stability (FxTS) under Input Constraints}
To encode time constraints, we utilize a relatively newer notion of stability, termed fixed-time stability (FxTS)~\citep{polyakov2012nonlinear}, which requires that the system trajectories converge to the equilibrium within a given \textit{fixed} time $T < \infty$. The following definition of FxTS and the corresponding Lyapunov conditions are adapted from~\cite{polyakov2012nonlinear}. Consider the autonomous dynamical system: 
\begin{align}\label{eq: sys auto}
\dot x(t) = f(x(t)),
\end{align}
where $x\in \mathbb R^n$, $f: \mathcal D \rightarrow \mathbb R^n$ is continuous on an open neighborhood $\mathcal D\subseteq \mathbb R^n$ of the origin and $f(0)=0$.

\begin{definition}[\textbf{FxTS}]\label{Def:FxTS auto}
The origin is a {FxTS} equilibrium of~\eqref{eq: sys auto} if it is Lyapunov stable and there exists a fixed time $T$ such that $\lim_{t\to T}x(t) = 0$ for all $x(0)\in \mathbb R^n$, i.e., the trajectories converge to the origin within a fixed time $T$. 
\end{definition}

The authors of~\cite{polyakov2012nonlinear} also presented Lyapunov conditions for the equilibrium of the \textit{uncontrained} system~\eqref{eq: sys auto} to be FxTS. 

\begin{theorem}[\textbf{FxTS conditions for unconstrained systems}]\label{Thm: FxTS Polyakov}
Suppose there exists a continuously differentiable, positive definite, radially unbounded function $V:\mathbb R^n\rightarrow\mathbb R$ such that 
\begin{align}\label{eq: dot V FxTS old}
    \dot V(x) \leq -\alpha_1V(x)^{\gamma_1}-\alpha_2V(x)^{\gamma_2},
\end{align}
holds for all $x\in \mathbb R^n \setminus\{0\}$, with $\alpha_1, \alpha_2>0$, $\gamma_1>1$ and $0<\gamma_2<1$. Then, the origin of~\eqref{eq: sys auto} is {FxTS} with continuous settling-time function $T$ that satisfies:
\begin{align}\label{eq: T bound old}
    T \leq \frac{1}{\alpha_1(\gamma_1-1)} + \frac{1}{\alpha_2(1-\gamma_2)}.
\end{align}
\end{theorem}

As illustrated in~\cite{garg2022fixed}, this Lyapunov result cannot be used for systems with input constraints. The modified Lyapunov conditions were given in~\cite{garg2022fixed,garg2021characterization} for FxTS under input constraints. 

\begin{theorem}[\textbf{New Lyapunov conditions for FxTS}]\label{Th: FxTS new}
Let $V:\mathbb R^n\rightarrow \mathbb R$ be a continuously differentiable, positive definite, radially unbounded function, satisfying
\begin{align}\label{eq: dot V new ineq}
     \dot V(x) \leq     -\alpha_1V(x)^{\gamma_1}-\alpha_2V(x)^{\gamma_2}+\delta_1V(x),
\end{align}
for all $x\in \mathbb R^n\setminus\{0\}$ along the trajectories of~\eqref{eq: sys auto} with $\alpha_1, \alpha_2>0$, $\delta_1\in \mathbb R$, $\gamma_1 = 1+\frac{1}{\mu}$, $\gamma_2 = 1-\frac{1}{\mu}$ and $\mu>1$. Then, there exists a neighborhood $D\subseteq \mathbb R^n$ of the origin such that for all $x(0)\in D$, the closed-trajectories of~\eqref{eq: sys auto} reach the origin within a fixed time $T$, where
\begin{align}
    D & = \begin{cases} \; \mathbb R^n; &  \frac{\delta_1}{2\sqrt{\alpha_1\alpha_2}}< 1,\\
    \left\{x\; |\; V(x)\leq k^\mu\left(\frac{\delta_1-\sqrt{\delta_1^2-4\alpha_1\alpha_2}}{2\alpha_1}\right)^\mu\right\}; & \frac{\delta_1}{2\sqrt{\alpha_1\alpha_2}}\geq 1, 
    \end{cases},\label{eq: domain of attraction}\\
    T & \leq \begin{cases}\frac{\mu\pi}{2\sqrt{\alpha_1\alpha_2}};& \quad \; \frac{\delta_1}{2\sqrt{\alpha_1\alpha_2}}\leq 0,\\
    \frac{\mu}{\alpha_1k_1}\left(\frac{\pi}{2}-\tan^{-1}k_2\right); & \quad \; 0 \leq \frac{\delta_1}{2\sqrt{\alpha_1\alpha_2}}<1,\\
    \frac{\mu}{\alpha_1(b-a)}\left(\log\left(\frac{b-ka}{a(1-k)}\right)-\log\left(\frac{b}{a}\right)\right); & \quad \;\frac{\delta_1}{2\sqrt{\alpha_1\alpha_2}}\geq 1,
    \end{cases},\label{new FxTS T est}
\end{align}
where $0<k<1$, $a< b$ are the solutions of $\gamma(z) \coloneqq \alpha_1z^2-\delta_1z+\alpha_2 = 0$,  $k_1 = \sqrt{\frac{4\alpha_1\alpha_2-\delta_1^2}{4\alpha_1^2}}$ and $k_2 = -\frac{\delta_1}{\sqrt{4\alpha_1\alpha_2-\delta_1^2}}$.
\end{theorem}

For a constrained control system, a relation between the domain of attraction, the time of convergence, and the input bounds using the new Lyapunov conditions~\eqref{eq: dot V new ineq} was developed in~\cite{garg2022fixed}. In brief, it was shown that the domain of attraction grows as the bounds on the input increases, or the required time of convergence increases, which also matches the basic intuition. Interested readers on the proof of this theorem and a more detailed discussion on this topic are referred to~\cite{garg2021advances}. Next, we illustrate how this modified Lyapunov condition naturally fits in a QP formulation for the concurrent problem of FxTS and safety, in the presence of input constraints.

\subsection{Concurrent FxTS and Safety}

Consider the nonlinear, control-affine system
\begin{align}\label{cont aff sys}
    \dot x = f(x) + g(x)u, \quad x(0) = x_0, 
\end{align}
where $x\in \mathbb R^n$ is the state vector, $f:\mathbb R^n\rightarrow \mathbb R^n$ and $g:\mathbb R^n\rightarrow \mathbb R^{n\times m}$ are system vector fields, continuous in their arguments, and $u\in \mathcal U\subset \mathbb R^m$ is the control input vector where $\mathcal U$ is the input constraint set. Let $h_S:\mathbb R^n\rightarrow \mathbb R$ and $h_G:\mathbb R^n\rightarrow \mathbb R$ be continuously differentiable functions. Define the safe set $S_S \coloneqq \{x\; |\; h_S(x)\leq 0\}$ such that its boundary and its interior $S_S$ are given as $\partial S_S \coloneqq \{x\; |\; h_S(x) = 0\}$ and $\textnormal{int}(S_S) \coloneqq \{x\; |\; h_S(x)<0\}$, respectively, to be rendered forward invariant under the closed-loop dynamics of~\eqref{cont aff sys}. Similarly, define the goal set $S_G \coloneqq  \{x\; |\; h_G(x)\leq 0\}$ such that its boundary and its interior $S_G$ are given as $\partial S_G \coloneqq \{x\; |\; h_G(x) = 0\}$ and $\textnormal{int}(S_G) \coloneqq \{x\; |\; h_G(x)<0\}$, respectively, to be reached by the closed-loop trajectories of~\eqref{cont aff sys} in a user-defined fixed time $T_{ud}>0$. 

\begin{assumption}\label{assum: sets Ss Sg}
 $S_G\bigcap S_S\neq \emptyset$, {the set $S_G$ is compact}, and the sets $S_S$ and $S_G$ have non-empty interiors. There exists a class-$\mathcal K_\infty$ function $\alpha_G$ such that $h_G(x)\geq \alpha_G(|x|_{S_G})$, for all $x\notin S_G$.
\end{assumption} 

The QP formulation in~\cite{garg2019control} uses the old FxTS Lyapunov conditions from Theorem \ref{Thm: FxTS Polyakov} along with the CBF condition from Definition \ref{definition::cbf_definition} for concurrent safety and FxTS. However, that formulation is incapable of handling input constraints. The formulation in~\cite{garg2022fixed} uses the new FxTS Lyapunov conditions from Theorem \ref{Th: FxTS new}, allowing incorporation of input constraints in the QP. The function $h_G$ is termed as \textbf{FxT-CLF} if it satisfies the new FxTS Lyapunov conditions in Theorem \ref{Th: FxTS new}, while the function $h_S$ is termed as a CBF if it satisfies the conditions in Definition \ref{definition::cbf_definition}. Next, we define the notion of the fixed-time domain of attraction for a compact set $S\subset \mathbb R^n$: 

\begin{definition}[\textbf{FxT-DoA}]\label{def:FxT-DoA} 
For a compact set $S_G\subset\mathbb R^n$, the set $D_S\subset \mathbb R^{n}$, satisfying $S_G\subset D_S$, is a Fixed-Time Domain of Attraction (FxT-DoA) with time $T>0$ for the closed-loop system~\eqref{cont aff sys} under $u$, if 
\begin{itemize}
    \item[i)] for all $x(0) \in D_S$, $x(t) \in D_S$ for all $t\in [0, T]$, and
    \item[ii)] there exists $T_0\in [0, T]$ such that $\lim_{t\to T_{0}}x(t) \in S_G$.
\end{itemize}
\end{definition}

\begin{figure}[!ht]
    \centering
    \includegraphics[ width=0.5\columnwidth]{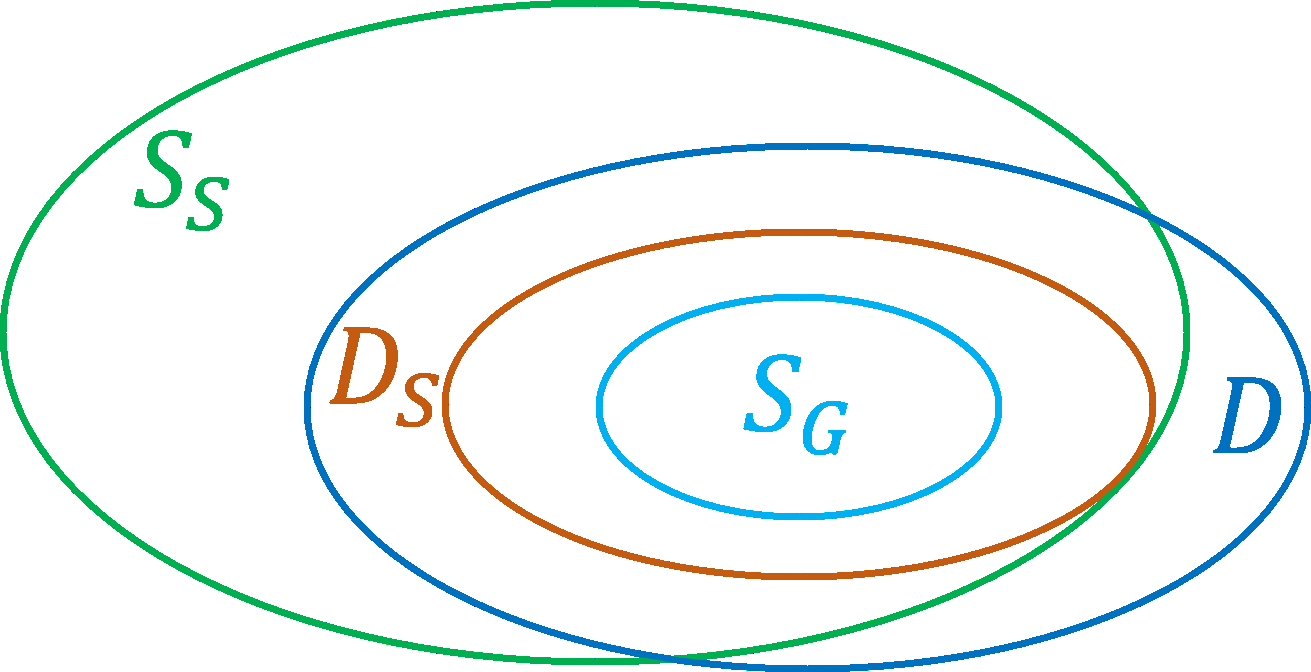}
    \caption{Illustration of the safe set $S_S$ (shown in green), the goal set $S_G$ (shown in light blue), FxT DoA $D$ (shown in dark blue) and the domain $D_S$ (shown in brown).}\label{fig:set DS}
\end{figure}

\begin{problem}\label{P reach S}
Design a control input $u \in \mathcal U \coloneqq \{v\in \mathbb R^m\; |\; A_uv\leq b_u\}$ and compute $D\subset\mathbb R^n$, so that for all $x_0\in D\subseteq S_S$, the closed-loop trajectories $x(t)$ of~\eqref{cont aff sys} satisfy $x(t) \in S_S$ for all $t\geq 0$, and $x(T_{ud})\in S_G$, where $T_{ud}>0$ is a user-defined fixed time and $D$ is a FxT-DoA for the set $S_G$.\footnote{Note that this problem can also be formulated using Signal Temporal Logic~\cite{lindemann2018control}, as stated in~\cite{garg2019control}.}
\end{problem}


In~\cite{garg2022fixed}, a QP-based feedback synthesis approach is presented to address Problem \ref{P reach S}. 
Define $z = \begin{bmatrix}v^T & \delta_1 & \delta_2\end{bmatrix}^T\in \mathbb R^{m+2}$, and consider the QP:
\\
\noindent\textbf{(FxT-CLF-CBF-QP)}
\begin{subequations}\label{QP gen}
\begin{align}
\min_{z\in \mathbb R^{m+2}} \; \frac{1}{2}z^THz & + F^Tz\\
    \textrm{s.t.} \quad \; A_uv  \leq & \; b_u, \label{C1 cont const}\\
    L_fh_G(x) + L_gh_G(x)v  \leq & \; \delta_1h_G(x)-\alpha_1\max\{0,h_G(x)\}^{\gamma_1} \nonumber\\
    & -\alpha_2\max\{0,h_G(x)\}^{\gamma_2} \label{C2 stab const}\\
    L_fh_S(x) + L_gh_S(x)v \leq &-\delta_2h_S(x),\label{C3 safe const}
\end{align}
\end{subequations}
where $H = \textrm{diag}\{w_{u_1},\ldots, w_{u_m}, w_1, w_2\}$ is a diagonal matrix consisting of positive weights $w_{u_i}, w_i>0$, $F = \begin{bmatrix}\mathbf 0_m^T & q_1 & 0\end{bmatrix}^T$ with $q_1>0$ and $\mathbf 0_m\in \mathbb R^m$ a column vector consisting of zeros. The parameters $\alpha_1, \alpha_2, \gamma_1, \gamma_2$ aare chosen as $\alpha_1 = \alpha_2 = \frac{\mu\pi}{2T_{ud}}$, $\gamma_1 = 1+\frac{1}{\mu}$ and $\gamma_2 = 1-\frac{1}{\mu}$ with $\mu>1$. The linear term $F^Tz = q_1\delta_1$ in the objective function of~\eqref{QP gen} penalizes the positive values of $\delta_1$. Constraint~\eqref{C1 cont const} imposes control input constraints. Constraint~\eqref{C2 stab const} is imposed for convergence of the closed-loop trajectories of~\eqref{cont aff sys} to the set $S_G$, and the constraint~\eqref{C3 safe const} is imposed for forward invariance of the set $S_S$.
The slack terms corresponding to $\delta_1, \delta_2$ allow the upper bounds of the time derivatives of $h_S(x)$ and $h_G(x)$, respectively, to have a positive term for $x$ such that $h_S(x)<0$ and $h_G(x)>0$. With this setup and under certain conditions, it was shown in~\cite{garg2022fixed} that the QP~\eqref{QP gen} is feasible (ensuring a control input exists), has a continuous solution (ensuring applicability of Nagumo's theorem for forward invariance) and guarantees both safety and FxTS from a domain that depends on the maximum value of the slack variable $\delta_1$. For simultaneous safety and FxT convergence, a subset $D_S\subset S_S$ of the FxT-DoA $D$ of the set $S_G$ can be defined so that its forward invariance per Lyapunov theorem results in safety and it being subset of the FxT-DoA results in FxT convergence (see Figure \ref{fig:set DS}).  

\begin{figure}[t]
    \centering
    \includegraphics[ width=0.6\columnwidth,clip]{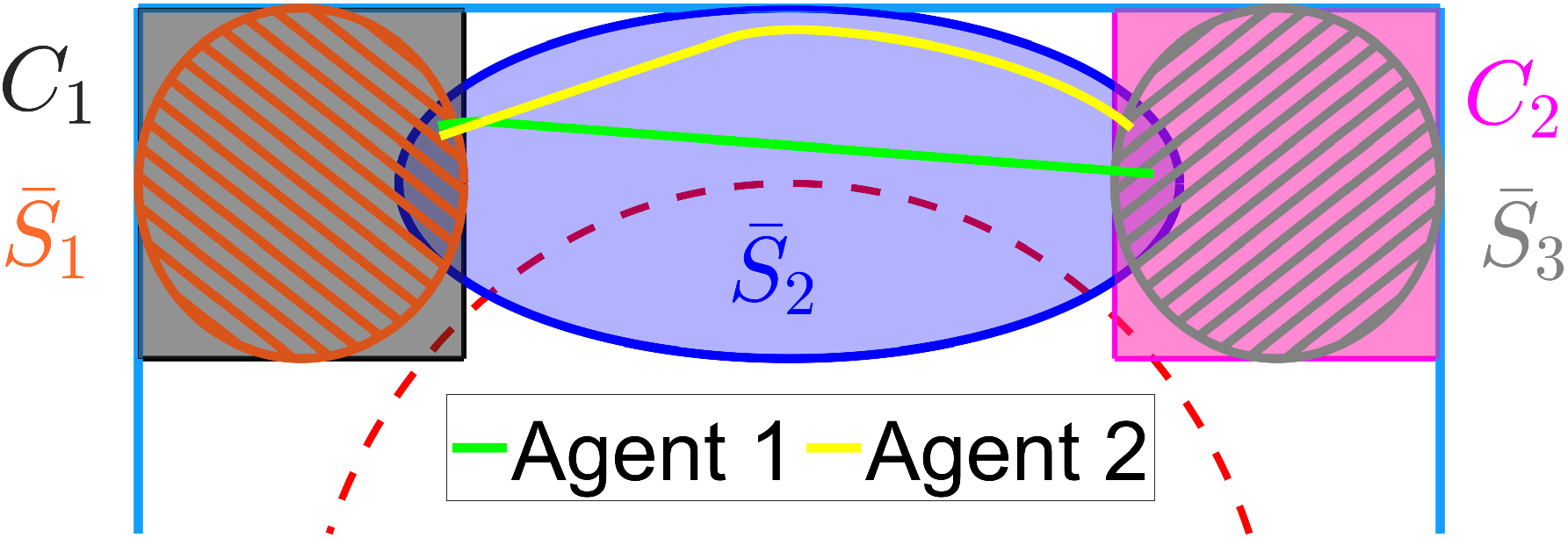}
    \caption{The resulting closed-loop paths of the agents.}\label{fig: unic traj}
\end{figure}

We present a two-agent motion planning example under spatiotemporal specifications, where the robot dynamics are modeled under constrained unicycle dynamics as $\dot x_i  = u_i\cos(\theta_i), \; \dot y_i  = u_i\sin(\theta_i), \; \dot \theta_i  = \omega_i,$ where $[x_i\; \ y_i]^T\in \mathbb R^2$ is the position vector of the agent $i$ for $i\in \{1, 2\}$, $\theta_i\in \mathbb R$ its orientation and $[u_i\; \ \omega_i]^T\in \mathbb R^2$ the control input vector comprising of the linear speed $u_i\in [0, u_M]$ and angular velocity $|\omega_i|\leq \omega_M$. The closed-loop trajectories for the agents, starting from $[x_1(0)\; \ y_1(0)]^T\in C_1 = \{z\in \mathbb R^2\; |\; \|z-[-1.5\; 1.5]^T\|_\infty\leq 0.5\}$ and $[x_2(0)\; \ y_2(0)]^T\in C_2 = \{z\in \mathbb R^2\; |\; \|z-[1.5\; 1.5]^T\|_\infty\leq 0.5\}$, respectively, are required to reach to sets $C_2$ and $C_1$, while staying inside the blue rectangle $\{z\in \mathbb R^2\; |\; \|z\|_\infty\leq 2\}$, and outside the red-dotted circle $\{z\in \mathbb R^2\; |\; \|z\|_2 \leq 1.5\}$, as shown in Figure \ref{fig: unic traj}. The agents are also maintaining an inter-agent distance $d_m>0$ at all times. 

\subsection{Robust FxTS and Robust Safety}
Next, we discuss how robustness to unmodeled phenomenal and measurement noise can be taken into consideration during control design. For this, consider a perturbed dynamical system:
\begin{align}\label{eq: NL pert cont affine}
    \dot x(t) & = f(x(t)) + g(x(t))u + d(t,x),
\end{align}
where $x\in \mathbb R^n, u \in \mathcal U\subset \mathbb R^m$ are the state and the control input vectors, respectively, with $\mathcal U$ the control input constraint set, $f:\mathbb R^n\rightarrow\mathbb R^n$ and $g: \mathbb R^n\rightarrow\mathbb R^{n\times m}$ are continuous functions and $d:\mathbb R_+\times \mathbb R^n\rightarrow \mathbb R^n$ is an unknown additive disturbance. The following assumption is made. 

\begin{assumption}[\textbf{Disturbance bound}]\label{assum d}
There exists $\gamma>0$ such that for all $t\geq 0$ and $x\in \mathcal D\subset \mathbb R^n$, $ \|d(t,x)\|\leq \gamma$, where $\mathcal D$ is a compact domain.
\end{assumption}

Encoding safety in the presence of disturbances can be done using robust CBFs~\cite{jankovic2018robust,cortez2019control, black2020quadratic}. In these works, however, only added process noise, or uncertainty in the state dynamics as in~\eqref{eq: NL pert cont affine}, is considered, and robust variants of FxT-CLF and CBF are introduced to guarantee convergence to a neighborhood of the goal set and safety. Here we take into account the effect of sensor noise and measurement uncertainties. More specifically, consider that only an estimate of the system state denoted as $\hat x$, is available, that satisfies:
\begin{align}
    \dot{\hat x} = f(\hat x) + g(\hat x)u.
\end{align}
The following assumption is made on the state-estimation error $\|x-\hat x\|$.  

\begin{assumption}[\textbf{Estimation error bound}]\label{assum hat x eps}
There exists an $\epsilon>0$ such that $\|\hat x(t)-x(t)\|\leq \epsilon$, for all $t\geq 0$. 
\end{assumption}

Then a robust variant of FxT-CLF and a robust variant of CBF is proposed in~\cite{garg2021robust} as follows: Corresponding to the set $S(t) = \{x\; |\; h(t,x)\leq 0\}$ where $h:\mathbb R_+\times\mathbb R^n\rightarrow\mathbb R$ is continuously differentiable, define $\hat S_\epsilon(t) = \{\hat x\; |\; h(t,\hat x)\leq -l\epsilon\}$, where $l = \sup \|\frac{\partial h(t,x)}{\partial x}\|$ is the Lipschitz constant of the function $h$. Inspired from~\cite{xu2015robustness}, the notion of a robust CBF is defined as follows.

\begin{definition}[\textbf{Robust CBF}]
A continuously differentiable function $h:\mathbb R_+\times\mathbb R^n\rightarrow\mathbb R$ is called a robust CBF for~\eqref{eq: NL pert cont affine} with respect to a disturbance $d$ satisfying Assumption \ref{assum d} if there exists a locally Lipschitz class-$\mathcal K$ function $\alpha$ such that the following condition holds
\begin{align}\label{eq: robust safety cond}
\begin{split}
    &\inf_{u\in \mathcal U}\left\{L_{f}h(t,x(t))+L_{g}h(t,x(t))u+\frac{\partial h}{\partial t}(t,x(t))\right\} \leq \alpha(-h(t, x(t)))-l\gamma,
\end{split}
\end{align}
for all $x(t)\in S(t)$ and $t\geq 0$. 
\end{definition}

\noindent Note that the worst-case bound $l\gamma$ of the term $\|L_dh(t,x)\|$ can be relaxed if more information than just the upper bound of the disturbance is known, or can be adapted online to reduce the conservatism. Some relevant work has been presented in~\cite{black2020quadratic, black2021ecc, black2022adaptation}. The existence of a robust CBF implies forward invariance of the set $S(t)$ for all $t\geq 0$, assuming that the system trajectories start with an initial $x(0)$ such that the measured or estimated state satisfies $\hat x(0)\in \hat S_\epsilon(0)$. 

Similarly we can define the notion of a robust FxT-CLF to guarantee {FxTS} of the closed-loop trajectories to the goal set. Consider a continuously differentiable function $V:\mathbb R^n\rightarrow\mathbb R$ with Lipschitz constant $l_V$.

\begin{definition}[\textbf{Robust FxT-CLF-$S$}]
A continuously differentiable function $V:\mathbb R^n\rightarrow \mathbb R$ is called a Robust {FxT-CLF}-$S$ for a set $S$ with respect to a disturbance $d$ satisfying Assumption \ref{assum d} if $V$ is positive definite and radially unbounded {wrt} the set $S$, $V(x)<0$ for $x\in \textnormal{int}(S)$, and satisfies
\begin{align}\label{eq: dot V new ineq robust}
\begin{split}
    \hspace{-5pt}\inf_{u\in \mathcal U}\{L_fV(x)+L_gV(x)u\} \leq -\alpha_1V(x)^{\gamma_1} -\alpha_2V(x)^{\gamma_2} + \delta_1V(x)-l_V\gamma,
\end{split}
\end{align}
with $\alpha_1, \alpha_2>0$, $\delta_1\in \mathbb R$, $\gamma_1 = 1+\frac{1}{\mu}$, $\gamma_2 = 1-\frac{1}{\mu}$ for $\mu>1$, along the trajectories of~\eqref{eq: NL pert cont affine}. 
\end{definition}
Using the mean value theorem, the following inequality can be obtained:
\begin{align}\label{eq: V hat V rel}
    V(x) & \leq V(\hat x) + l_V\epsilon,
\end{align}
which implies that if $V(\hat x) \leq -l_V\epsilon$, then $V(x)\leq 0$. Based on this, it is shown in~\cite{garg2021robust} that existence of a robust FxT-CLF for the set $S_G$ implies existence of neighborhood $D$ of the set $S_G$ such that fixed time convergence of the closed-loop trajectories of $x$ is guaranteed for all initial conditions $x(0)$ such that the estimated state satisfies $\hat x(0)\in D$. 

Note that to encode safety with respect to a general time-varying safe set, let $h_{T} :\mathbb R_+\times\mathbb R^n\rightarrow \mathbb R$ be a continuously differentiable function defining the \textit{time-varying} safe set $S_{T}(t) = \{x\; |\; h_{T}(t,x)\leq 0\}$. Now we are ready to present the QP formulation for a robust control synthesis under input constraints. For the sake of brevity, we omit the arguments $\hat x$ and $(t,\hat x)$. Define $z = \begin{bmatrix}v^T & \delta_{1} & \delta_{2}& \delta_{3}\end{bmatrix}^T\in \mathbb R^{m+3}$, and consider the following optimization problem:
\\
\noindent\textbf{(Robust FxT-CLF-CBF QP)}
{\small{
\begin{subequations}\label{QP gen ideal}
\begin{align}
\min_{z\in \mathbb R^{m+3}} \;&  \frac{1}{2}z^THz  + F^Tz\\
    \textrm{s.t.} \quad \quad  A_{u}v  \leq & \; b_{u}, \label{C1 cont const ideal}\\
    L_{f}\hat h_{G} + L_{g}\hat h_{G}v  \leq & \; \delta_{1}\hat h_{G}-\alpha_{1}\max\{0,\hat h_{G}\}^{\gamma_{1}} \nonumber\\
    &\quad -\alpha_{2}\max\{0,\hat h_{G}\}^{\gamma_{2}}-l_{G}\gamma \label{C2 stab const ideal}\\
    L_{f}\hat h_S + L_{g}\hat h_Sv \leq &-\delta_{2}\hat h_S-l_S\gamma,\label{C3 safe const ideal}\\
    L_{f}\hat h_{T} + L_{g}\hat h_{T}v \leq &-\delta_{3}\hat h_{T}-\frac{\partial \hat h_{T}}{\partial t}-l_T\gamma,\label{C4 tv safe const ideal}
\end{align}
\end{subequations}}}\normalsize
where $H = \textrm{diag}\{\{w_{u_l}\}, w_1, w_2, w_3\}$ is a diagonal matrix consisting of positive weights $w_{u_l}, w_1, w_2,
w_3>0$ for $l = 1, 2, \dots, m$, $F = \begin{bmatrix}\mathbf 0_m^T & q & 0 & 0\end{bmatrix}^T$ with $q>0$ and functions $\hat h_G,\hat h_S$ (respectively, $\hat h_T$) are functions of $\hat x$ (respectively, $(t,\hat x)$) defined as follows. For any function $\phi:\mathbb R_+\times \mathbb R^n\rightarrow\mathbb R$ with Lipschitz constant $l_\phi$, define 
\begin{align}\label{eq: hat func def}
\hat \phi(t,\cdot) = \phi(t,\cdot) + l_\phi\epsilon.
\end{align}
The parameters $\alpha_{1}, \alpha_{2}, \gamma_{1}, \gamma_{2}$ are chosen as $\alpha_{1} = \alpha_{2} = \mu\pi/(2\bar T)$, $\gamma_{1} = 1+\frac{1}{\mu}$ and $\gamma_{2} = 1-\frac{1}{\mu}$ with $\mu>1$ and $\bar T$ the user-defined time in Problem \ref{P reach S}. With this robust control design framework, under some technical assumptions and conditions, it is shown in~\cite{garg2021robust} that the QP~\eqref{QP gen ideal} is feasible, its solution is continuous and results in both safety of $S_S$ and FxTS of $S_G$, from domain of attraction that depends on the maximum value of the slack variable $\delta_1$. 

\begin{figure}[t]
    \centering
    \includegraphics[width=0.7\columnwidth,clip]{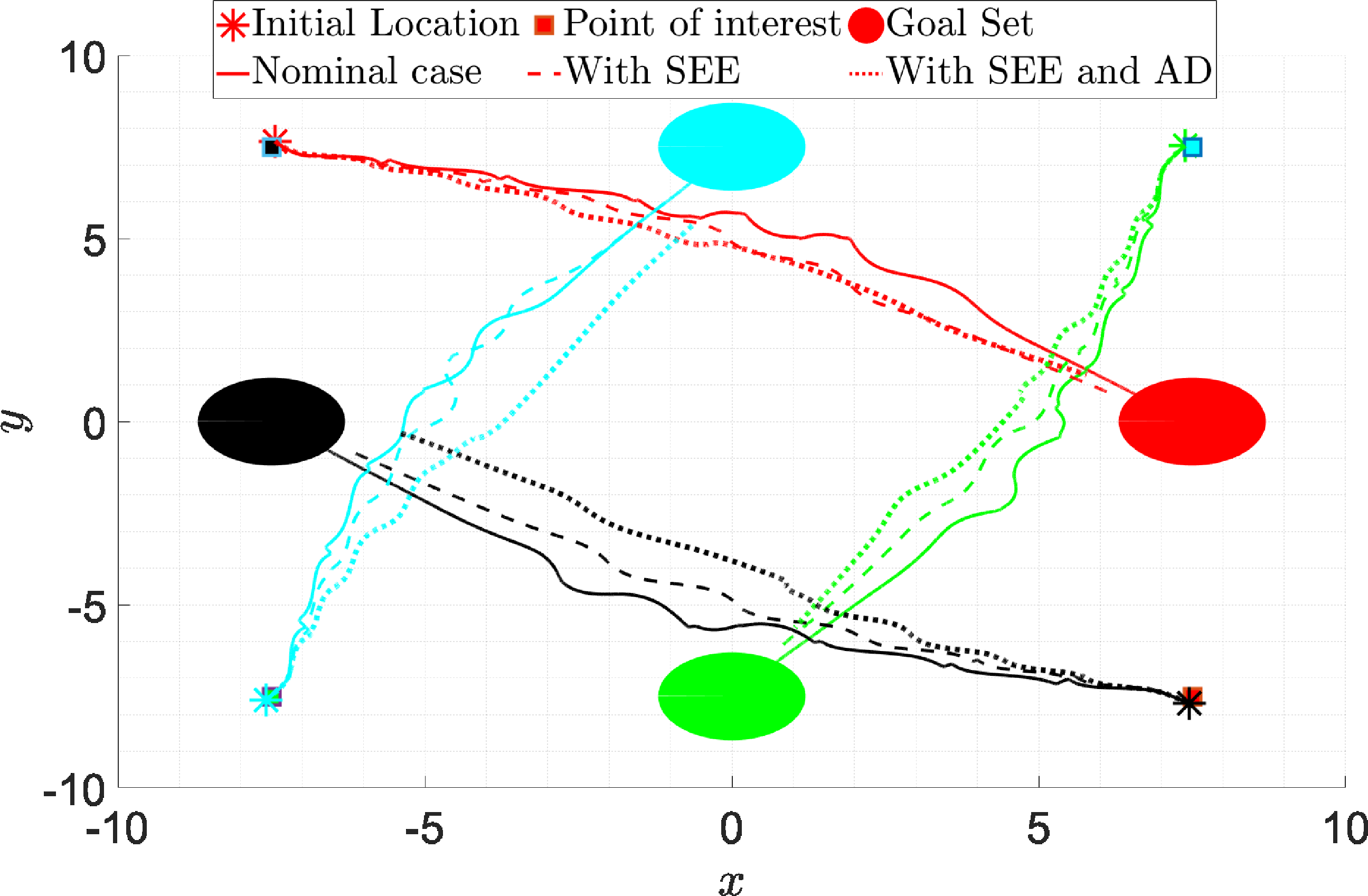}
    \caption{Closed-loop paths traced by agents in a 4 agents scenario for the nominal case, i.e., without any disturbance (solid lines), with only state estimation error (SEE) (dashed lines) and with both SEE and additive disturbance (AD) (dotted lines).}\label{fig:4 agent traj}
\end{figure}

In the interest of space, we are not including a detailed case study here, however interested readers are referred to~\cite{garg2021robust}, which includes the problem of navigating multiple nonlinear underactuated marine vehicles while respecting visual sensing constraints, avoiding collisions (encoding safety constraints), and moving towards desired destinations (encoding convergence constraints) under additive disturbances (currents) and navigation (state estimation) error. The closed-loop paths of four vehicles are shown in Fig. \ref{fig:4 agent traj}.

\subsection{Extensions and Relevant Work}
One of the limitations of CLF-CBF QPs for convergence and safety, as analyzed in~\cite{reis2020control}, is the existence of stable undesirable equilibrium. Furthermore, while CBF-based approaches can guarantee step-wise safety (i.e., safety at each step) and hence are termed \textit{myopic} in nature~\cite{cohen2020approximate}, they cannot guarantee that system trajectories will not enter a region in the future from where safety cannot be guaranteed. Without proper knowledge of the \textit{control invariant} set, a wrongly chosen barrier function might lead to infeasibility of the CBF-QP and as a result, violation of safety. To circumvent these issues, combining a high-level planner with a low-level controller has become a popular approach~\cite{herbert2017fastrack,yin2019optimization, smith2019continuous,singh2017robust,rosolia2020multi}. The underlying idea in these strategies is to design low-level controllers to track a reference trajectory, which is computed by a high-level planner using a simplified model. However, it is important that the low-level controller is able to track trajectory generated by high-level in a given time dictated by the update frequency of the high-level planner. To this end, the notion of FxTS is utilized in~\cite{garg2021multi}, where a FxT-CLF-CBF-QP-based low-level controller guarantees that the trajectories remain in the domain of attraction of the next waypoint, and reach there before the next high-level-planning update occurs. In turn, a model predictive control (MPC)-based high-level planner utilizes the FxT-DoA to generate trajectories so that the low-level QP is guaranteed to remain feasible. This way, the low-level controller helps guarantee the recursive feasibility of the MPC, and the high-level planner helps guarantee the feasibility of the QP, thereby guaranteeing that the underlying problem can be solved. In~\cite{garg2021multi}, we also introduced a new notion of safety, termed Periodic Safety, where the system trajectories are required to enter or visit a set (say, $\mathcal X_T$) periodically (say, with period $T>0$), while remaining in a safe set $\mathcal X$ at all times. 

\begin{figure}[t]
\centering
\includegraphics[width=0.98\linewidth,height = 110pt]{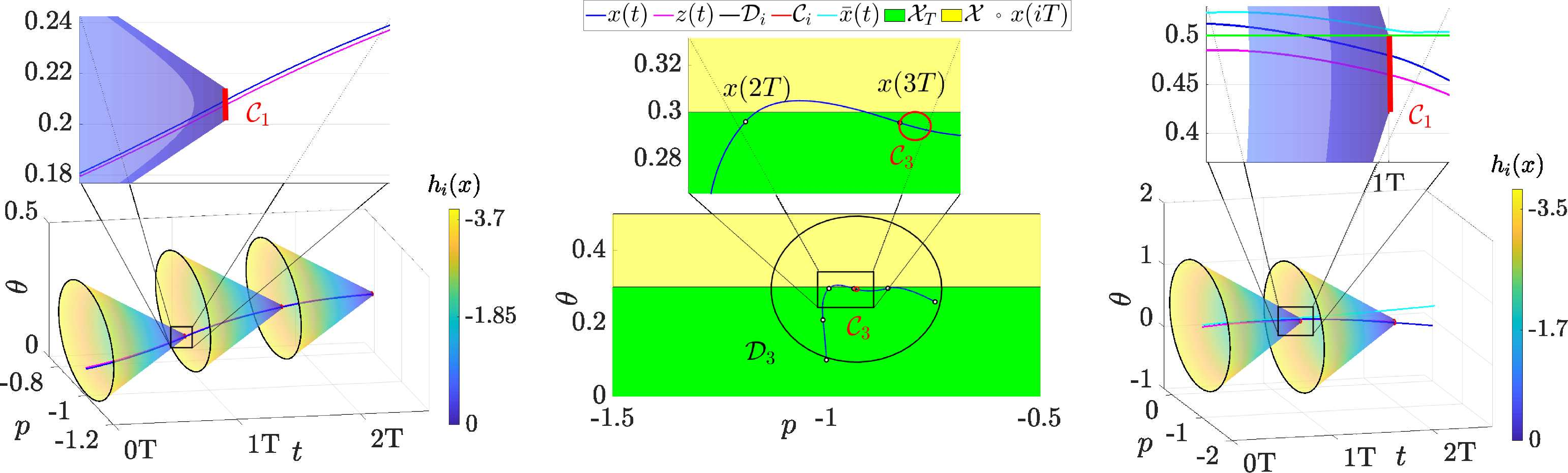}
\caption{Case study demonstrating the proposed method. \textbf{(Left)} The trajectory $x(t)$ projected onto $p-\theta-t$ axes. In each interval $[iT, (i+1)T)$, FxT-DoA $\mathcal{D}_i$ is depicted as the colored region which decays to the set of feasibility of MPC ($\mathcal C_i$) by the end of the interval. \textbf{(Middle)} Projection of the closed-loop trajectory on the $p-\theta$ plane. The trajectory leaves and enters the set $\mathcal X_T$ after $t = 2T$ and before $t = 3T$, respectively. \textbf{(Right)} Illustration of a scenario where the trajectory $\bar x(t)$ generated using an exponentially stabilizing CLF fails to enter $\mathcal C_1$, leading to infeasibility of the MPC at $t = 1T$, whereas the trajectory $x(t)$ generated by the proposed method enters the set before $t = 1T$.}
\label{fig:multi rate}
\end{figure}

In the interest of space, we skip the technical details of the hierarchical framework and briefly discuss the case study that illustrates utility of such as approach. We use the proposed strategy to steer a Segway to the origin.\footnote{\scriptsize{Code available at
\href{https://github.com/kunalgarg42/fxts_multi_rate}{\texttt{github.com/kunalgarg42/fxts\_multi\_rate}}
}}
The state of the system are the position $p$, the velocity $v$, the rod angle $\theta$ and the angular velocity $\omega$. The control action is the voltage commanded to the motor and the equations of motion used to simulate the system can be found in~\cite[Section~IV.B]{gurriet2018towards}. In this simulation, we run the high-level MPC planner at $5$Hz and the low-level controller at $10$kHz. We choose the set $\mathcal X_T = \{x = [p, v, \theta, \omega]^T\; |\; |p|\leq 10, |v|\leq 5, |\theta|\leq 0.3, |\omega|\leq 10\pi\}, \mathcal X_F = \{0\}$, input bounds $\|u\|\leq 25$ with $\|u_m\|\leq 15$. From Figure \ref{fig:multi rate}, the main takeaway is that periodically, using the proposed FxT-CLF-QP, the closed-loop trajectories reach the set from where feasibility of MPC is guaranteed (denoted $\mathcal C_i$ for $i-$th MPC step). However, an exponentially stabiling controller fails to do so, resulting in infeasibility of the MPC. This demonstrates the efficacy of the proposed framework over the existing methods that use exponentially stabilizing controllers.

\section{Input-Constrained Control Barrier Functions: Synthesis under High Relative Degree and Disturbances}
\label{section::HOCBF}

\subsection{Constructive Methods for Higher-Order CBFs under Disturbances and Input Constraints}
High order CBFs (HOCBFs) were first introduced in~\cite{xiao2019control} and extended to be robust in~\cite{tan2021high,xiao2021high}. However, a limitation of all these works is that it is unclear how to choose the class-$\mathcal{K}$ functions in the HOCBF construction. When there are no input constraints, choosing these functions is equivalent to tuning the control law. When there are input constraints, these functions determine the size and shape of the CBF set, and thus must be chosen carefully to ensure satisfaction of~\eqref{eq::cbf_derivative}, here modified in Definition~\ref{def:rcbf} to include robustness. Thus, the objective of~\cite{Breeden2021InputConstraints,breeden2021robust} is to develop constructive methods to choose these functions. This section presents one such method from our work in~\cite{breeden2021robust}, and the interested reader is referred to~\cite[Sec.~3]{breeden2021robust} for two additional methods. Related works also include, non-exhaustively,~\cite{ShawCortez2022euler,ShawCortez2022multiple,squires2018backup,gurriet2018towards}, and the following method is further extended to high-order robust sampled-data CBFs in~\cite{breeden2023aiaa}.

Consider the time-varying control-affine model
\begin{equation}
    \dot{x} = \underbrace{f(t,x) + g(t,x) (u+w_u) + w_x}_{= F(t,x,u,w_u,w_x)} \,, \label{eq:modeldist}
\end{equation}
with time $t\in \mathcal{T}=[t_0,{\color{black}\infty)}$, state $x\in\reals^n$, control input $u\in \mathcal U \subset \reals^m$ where $\mathcal{U}$ is compact, unknown disturbances $w_u \in \reals^m$ and $w_x\in\reals^n$ that are continuous in time, and functions $f:\mathcal{T}\times\reals^n\rightarrow\reals^n$ and $g:\mathcal{T}\times\reals^n\rightarrow\reals^{n\times m}$ that are piecewise continuous in $t$ and locally Lipschitz continuous in $x$. Let $w_u$ and $w_x$ be bounded as $\|w_u\| \leq w_{u,\textrm{max}}$ and $\|w_x\| \leq w_{x,\textrm{max}}$ for some $w_{u,\textrm{max}},w_{x,\textrm{max}} \in \reals_{\geq 0}$, and define the set of allowable disturbances $\mathcal{W} \triangleq \{ w_u \in \reals^m \mid \|w_u\|\leq w_{u,\textrm{max}} \} \times \{ w_x \in \reals^n \mid \|{\color{black}w_x}\|\leq w_{x,\textrm{max}}\}$. 
{\color{black}Assume a unique solution to~\eqref{eq:modeldist} exists for all $t\in\mathcal{T}$.}
Given dynamics~\eqref{eq:modeldist}, a function $\eta : \mathcal{T}\times\reals^n \rightarrow \reals$ is said to be of relative-degree $r$ if it is $r$-times total differentiable in time and $\eta^{(r)}$ is the lowest order derivative in which $u$ {\color{black}and $w_u$} appear explicitly. Denote the set of all relative-degree $r$ functions as $\mathcal{G}^r$.

\begin{figure}
    \centering
    \includegraphics[width=0.8\columnwidth]{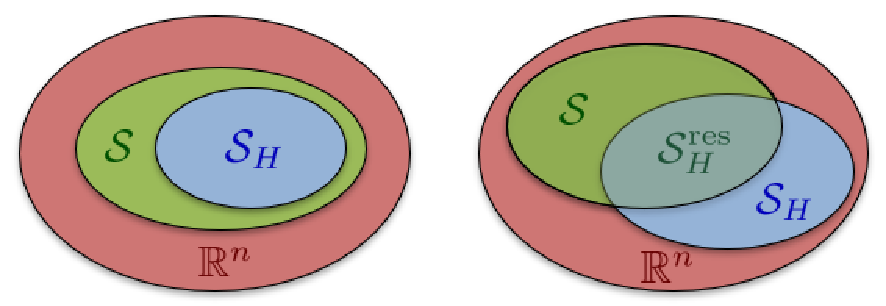}
    \caption{\color{black}Given a safe set $\mathcal{S} \subset \reals^n$ as in~\eqref{eq:safe_set_dist} and input constraints $\mathcal{U}$, this section presents a method of finding a viability domain $\mathcal{S}_H$ as in~\eqref{eq:inner_safe_set} (left) or $\mathcal{S}_H^\textrm{res}$ as in~\eqref{eq:restricted_safe_set} (right). Theorem~\ref{thm:variable_cbf} then guarantees that $\mathcal{S}_H^\textrm{res}$ can be rendered forward invariant under any disturbances $(w_u,w_x)\in\mathcal{W}$ while always satisfying the input constraints $\mathcal{U}$.}
    \label{fig:sets_fig}
\end{figure}

Let $h : \mathcal{T}\times\reals^n \rightarrow \reals$, $h \in \highdegree^r$, denote the \textit{constraint function}, and define a safe set $\mathcal{S}$ as
\begin{equation}
    \mathcal S(t) \triangleq \{x\in\reals^n\mid h(t,x) \leq 0\} \,, \label{eq:safe_set_dist}
\end{equation}
where we will henceforth drop the argument $t$ for compactness. Also, denote the safe set across time as $\mathcal{S}^\mathcal{T} \triangleq \{ (t,x) \in \reals\times\reals^n \mid t\in \mathcal{T}, x \in \mathcal S(t) \}$. Our aim is to develop methods for rendering the state trajectory always inside the safe set $\mathcal S$ {\color{black}in the presence of any allowable disturbances $(w_u,w_x)\in\mathcal{W}$.} 
We will do this by constructing functions $H:\mathcal{T}\times\reals^n\rightarrow\reals$ 
{\color{black}that}
generate sets of the form
\begin{gather}
    \mathcal S_H(t) \triangleq \{ x \in \reals^n \mid H(t,x) \leq 0 \} \,,  \label{eq:inner_safe_set} \\
    \color{black}\mathcal S_H^{\textrm{\textnormal{res}}}(t) \triangleq \{ x \in \reals^n \mid H(t,x) \leq 0 \textrm{\textnormal{ and }} h(t,x) \leq 0 \} {\color{black}\,,} \label{eq:restricted_safe_set}
\end{gather}
{\color{black}visualized in Fig.~\ref{fig:sets_fig}.} We refer to the set $\mathcal S_H$ as an \textit{inner safe set} (or also a \emph{CBF set}){\color{black}, and to the set $\mathcal{S}_H^\textrm{res}$ as a \textit{restricted safe set}. Note that if $H(t,x) \geq h(t,x)$ for all $(t,x)\in\mathcal{T}\times\reals^n$, then $\mathcal{S}_H\equiv \mathcal{S}_H^\textrm{res}$.} A controller is said to render ${\color{black}\mathcal S_H^\textrm{res}}$ forward invariant, if given any $x(t_0) \in {\color{black}\mathcal S_H^\textrm{res}(t_0)}$, the closed-loop trajectory satisfies $x(t)\in{\color{black}\mathcal{S}_H^\textrm{res}(t)}, \forall t\in\mathcal{T}$. 
In general, {\color{black}there may exist points $x(t_0)\in \mathcal{S}(t_0)$, from which} we will not be able to render $\mathcal S$ forward invariant under~\eqref{eq:modeldist}. Nevertheless, {\color{black}if we can render the subset}
${\color{black}\mathcal S_H^\textrm{res} \subseteq \mathcal{S}}$ forward invariant, {\color{black}then} we {\color{black}can} ensure that the closed loop trajectories of~\eqref{eq:modeldist} are safe {\color{black}(i.e. always stay in $\mathcal{S}$)} for initial conditions lying in the set ${\color{black}\mathcal S_H^\textrm{res}}$. Thus, a crucial requirement is that $x(t_0)\in{\color{black}\mathcal{S}_H^\textrm{res}}(t_0)$. We also define the {\color{black}domains $\mathcal{S}_H^\mathcal{T} \triangleq\{(t,x)\in\reals\times\reals^n\mid t\in\mathcal{T},x\in\mathcal{S}_H(t)\}$ and $\mathcal{S}_H^{\textrm{res},\mathcal{T}} \triangleq\{(t,x)\in\reals\times\reals^n\mid t\in\mathcal{T},x\in\mathcal{S}_H^\textrm{res}(t)\}$ similar} to $\mathcal{S}^\mathcal{T}$.

\subsubsection{Robust HOCBF Definition}

Here, we break from the HOCBF convention, and instead work with first-order CBFs. We also only consider relative-degree 2 constraint functions presently, though~\cite{breeden2021robust} also presents one method for greater relative-degrees. Let $\partial_t$ denote the partial derivative in time $t$ and $\nabla$ the gradient in state $x$.

\begin{definition}[\textbf{Robust CBF}] \label{def:rcbf}
    For the system~\eqref{eq:modeldist}, a continuously differentiable function $H:\mathcal{D}\times\reals^n\rightarrow\reals$ is a \emph{robust control barrier function (RCBF)} on a time-varying set $\mathcal{X}$ if there exists a function $\alpha\in\mathcal{K}$ such that $\forall x \in \mathcal X(t), t\in \mathcal{T}$,
    \begin{multline}
        \max_{(w_u,w_x)\in\mathcal{W}} \Big[ \inf_{u\in \mathcal U} {\color{black}\Big(}\partial_t H(t,x) +  \nabla H(t,x) {\color{black}\big[}f(t,x) \\ + g(t,x) (u+w_u) + w_x{\color{black}\big] \Big)} \Big] \leq \alpha(-H(t,x)) \,. \label{eq:rcbf_definition}
    \end{multline}
\end{definition}
Based on Definition~\ref{def:rcbf}, we also define for compactness
\begin{align}
    W(t, x) &\triangleq \max_{(w_u,w_x)\in\mathcal{W}} \nabla H(t,x) (g(t,x) w_u + w_x) \label{eq:max_error} \\ 
    &\equiv \|\nabla H(t,x) g(t,x)\| w_{u,\textrm{max}} + \|\nabla H(t,x)\| w_{x,\textrm{max}} \,. \nonumber
\end{align}
The set of control inputs such that~\eqref{eq:rcbf_definition} is satisfied is then
\begin{multline} \label{eq:valid_rcbf_control}
    \boldsymbol\mu_\textrm{rcbf} (t,x) \triangleq \{ u\in\mathcal{U} \mid \partial_t H(t,x) + \nabla H(t,x)(f(t,x) \\ +g(t,x)u) \leq \alpha(-H(t,x)) - W(t,x) \} \,.
\end{multline}
Note that since Definition~\ref{def:rcbf} considers the allowable control set $\mathcal{U}$, if $H$ is a RCBF on $\mathcal{X}$, then $\boldsymbol\mu_\textrm{rcbf}(t,x)$ is nonempty for all $x\in\mathcal{X}(t),t\in\mathcal{T}$.

\subsubsection{One Method for Constructing an HOCBF}

For the system~\eqref{eq:modeldist}, if $h$ is of relative-degree 2, note that $\dot{h}$ is a function of $w_x$, and $\ddot{h}$ is a function of $w_u$ and $w_x$, and thus are not precisely known. Thus, define the following upper bound on $\dot{h}$:
\begin{align}
    \dot{h}_w(t,x) \triangleq& \max_{||w_x||\leq w_{x,\textrm{max}}} \dot{h}(t,x,w_x)  \label{eq:def_hw} \\ =& \partial_t h(t,x) + \nabla h(t,x) f(t,x) + \|\nabla h(t,x)\|w_{x,\textrm{max}} \nonumber
\end{align}
and its derivative
\begin{align}
    \ddot{h}_w(t,x,u,w_u,w_x) &= \frac{d}{dt}\dot{h}_w(t,x) \\ &= \partial_t \dot{h}_w(t,x) + \nabla \dot{h}_w(t,x) F(t,x,u,w_u,w_x) \,. \nonumber 
\end{align}
Note that $\dot{h}_w$ is a known quantity, while $\ddot{h}_w$ is still a function of the unknown quantities $w_u,w_x$ in $F$.

For a relative-degree 2 constraint, we can intuitively describe $h$ as the position of an agent with respect to an obstacle, $\dot{h}$ its velocity, and $\ddot{h}$ its acceleration, where acceleration is the controlled variable.
Given some maximal amount of control authority encoded in $\mathcal{U}$, suppose that there exists some function $\phi:\reals\rightarrow\reals$ such that
\begin{equation}
        \max_{(w_u,w_x)\in \mathcal{W}} \inf_{u\in \mathcal U} \ddot{h}_w(t,x,u,w_u,w_x) \leq \phi(h(t,x)) \color{black}<\color{black} 0 \,. \label{eq:phi_requirement}
\end{equation}
This is a reasonable assumption for many systems, since intuitively $\phi$ represents the effects of other forces/accelerations in the environment. Given models of these forces, one can often read the function $\phi$ directly from the dynamics. If no such function $\phi$ exists, it may instead be possible to find such a function $\phi$ for a tighter constraint function, e.g. $h^\dagger = h + \gamma$. We then have the following theorem.

\begin{theorem}[\textbf{Method to Construct a RCBF}] \label{thm:variable_cbf}
Let $h \in \highdegree^2$ define a safe set as in~\eqref{eq:safe_set_dist}. Suppose there exists an invertible, continuously differentiable, and strictly monotone decreasing function $\Phi:\reals\rightarrow\reals$, whose derivative is $\Phi'=\phi$ for $\phi:\reals\rightarrow\reals$, such that~\eqref{eq:phi_requirement} holds $\forall (t,x) \in \mathcal{S}^\mathcal{T}$. Let $\Phi^{-1}$ be the function for which $\Phi^{-1}(\Phi(\lambda)) = \lambda, \forall \lambda\in\reals$. Then the function 
    \begin{equation}
        H(t,x) = \Phi^{-1}\left( \Phi(h(t,x)) - \frac{1}{2}\dot{h}_w(t,x)|\dot{h}_w(t,x)| \right) \label{eq:variable_cbf}
    \end{equation}
    is a RCBF on $\mathcal{S}_H^{\textrm{\textnormal{res}}}$ in~\eqref{eq:restricted_safe_set} for the system~\eqref{eq:modeldist} for any $\alpha\in\mathcal{K}$. Moreover, any control law $u(t,x)$ such that $u(t,x) \in \boldsymbol\mu_\textrm{\textnormal{rcbf}}(t,x), \forall (t,x)\in\mathcal{S}_H^{\textrm{\textnormal{res}},\mathcal{T}}$ also renders $\mathcal{S}_H^\textrm{\textnormal{res}}$ forward invariant.
\end{theorem}

That is, if~\eqref{eq:phi_requirement} holds, then we have a constructive way to find a CBF $H$ as a function of $h$ and $\dot{h}_w$. Note that the intersection of sets $\mathcal{S}_H^\textrm{res} = \mathcal{S}_H \cap \mathcal{S}$ is analogous to the intersection $C_1 \cap C_2$ in~\cite[Def.~7]{xiao2019control}. See~\cite{breeden2021robust} for more information.

\subsubsection{Case Study and Remarks} \label{sec:sims_rhocbf}

\newcommand{\transpose}[0]{^\textrm{\textnormal{T}}}

To see Theorem~\ref{thm:variable_cbf} in practice, consider the system with state $x = [r\transpose, \; \dot{r}\transpose]\transpose \in \reals^6$ with dynamics
\begin{equation}
    \ddot{r} = -\frac{\mu}{\|r\|^3}r + u
\end{equation}
for $\mu = 6.26(10)^{10}$, and constraint function
\begin{equation}
    h = \rho - \|r\|
\end{equation}
for $\rho = 4.76(10)^5$. Let $\mathcal{U} = \{ u\in\reals^3 \mid \| u\|_\infty \leq u_\textrm{max}\}$. For this constraint function, it holds that
\begin{equation}
    \ddot{h}_w \leq \frac{\mu}{\|r\|^2} - \frac{r\transpose (u + w_x + w_u)}{\|r\|}
\end{equation}
It follows that
\begin{equation}
    \max_{(w_u,w_x)\in\mathcal{W}} \inf_{u\in\mathcal{U}} \ddot{h}_w \leq \underbrace{\frac{\mu}{(\rho - h)^2} + w_{x,\textrm{max}} + w_{u,\textrm{max}} - u_\textrm{max}}_{=\phi(h)} \label{eq:phi_real}
\end{equation}
Assuming that~\eqref{eq:phi_real} is always negative for $h \leq 0$, then let $\Phi$ be any anti-derivative of $\phi$ in~\eqref{eq:phi_real}, such as
\begin{equation}
    \Phi(\lambda) = \frac{\mu}{\rho - \lambda} + (w_{u,\textrm{max}} + w_{x,\textrm{max}} - u_\textrm{max}) \lambda \label{eq:implemented_phi} \,.
\end{equation}
and Theorem~\ref{thm:variable_cbf} guarantees that $H$ as in~\eqref{eq:variable_cbf} with~\eqref{eq:implemented_phi} is a RCBF. Thus, we have a constructive method of constructing a RCBF for this system. This RCBF was then used in simulation in~\cite{breeden2021robust} and~\cite{breeden2022docking}.

We now consider how the above approach relates to the more widely used HOCBF formulation in~\cite{xiao2019control} and to the Exponential CBF (ECBF) formulation in~\cite{nguyen2016exponential} (which is a special case of~\cite{xiao2019control}). Given a relative-degree 2 constraint function $h$ meeting the assumptions of Theorem~\ref{thm:variable_cbf}, $h$ is also an HOCBF 
\begin{subequations}\begin{align}
    \psi_0 &= h \\[-2pt]
    \psi_1 &= \dot{h} - \alpha_1(-h) \\[-2pt]
    \psi_2 &= \ddot{h} + \alpha_1'(-h) \dot{h} - \alpha_2(-\psi_1)
\end{align}\end{subequations}
with choice
\begin{equation}
    \alpha_1(\lambda) = \sqrt{2 (\Phi(-\lambda) - \Phi(0))} \label{eq:equiv_alpha}
\end{equation}
and with $\alpha_2\in\mathcal{K}$ as a free variable. Thus, one can map between the approaches in~\cite{breeden2021robust} and~\cite{xiao2019control}. However, the choice of $\alpha_1$ in~\eqref{eq:equiv_alpha} is 1) non-obvious without the above analysis, and 2) violates the Lipschitzness assumptions present in~\cite{xiao2019control}. Also, while our method is constructive, it is conservative in the sense that it results in a CBF that is valid for any class-$\mathcal{K}$ function $\alpha$ in~\eqref{eq:valid_rcbf_control}, and as a result is valid for any $\alpha_2$ using the conventions of~\cite{xiao2019control}; Theorem~\ref{thm:variable_cbf} could potentially be applicable to a wider class of systems if this was relaxed.

Lastly, the ECBF is an HOCBF that uses only linear class-$\mathcal{K}$ functions. We note that if $\ddot{h}$ is bounded (as is usually the case when $\mathcal{U}$ is compact), then the ECBF can only be used with compact safe sets. This is because the ECBF, similar to a linear control law, requires stronger accelerations $\ddot{h}$, and hence larger $u$, as the state moves further from the safe set boundary. By contrast, the CBF in Theorem~\ref{thm:variable_cbf} works with an unbounded safe set and admits an unbounded inner safe set while still only commanding signals $u\in\mathcal{U}$ everywhere in the CBF set.

\subsubsection{Robustly Reachable Sets}

Finally, since the system~\eqref{eq:modeldist} is uncertain and the control~\eqref{eq:valid_rcbf_control} always considers the worst-case disturbance, it is worth considering the set of states that the system might reach. Frequently, the system evolution can be divided into arcs where either 1) the CBF condition is inactive, or 2) the CBF condition is satisfied with equality. Consider the behavior under the latter case.

\begin{theorem}[\textbf{RCBF Asymptotic Set}] \label{lemma:approach}
    Suppose $H:\mathcal{T}\times\reals^n\rightarrow\reals$ is a RCBF on the set $\mathcal{S}_H(t)$ (or $\mathcal{S}_H^\textrm{res}$ if using $H$ in~\eqref{eq:variable_cbf}) for the system~\eqref{eq:modeldist}. Suppose there exists constants $\eta_1, \eta_2 > 0$ such that $W$ in~\eqref{eq:max_error} satisfies $W(t,x) \in [\eta_1,\eta_2], \forall x\in\mathcal \mathcal{S}_H(t), t\in\mathcal{T}$. Let $\alpha_w\in\mathcal{K}$. Suppose $H(t_0,x(t_0)) \leq 0$. Then any control law $u(t,x)$ that satisfies 
    \begin{multline}
       \partial_t H(\cdot) + \nabla H(\cdot) (f(t,x) + g(t,x) u) + W(t,x) = \alpha_w(-H(\cdot)) W(t,x) \label{eq:cbf_condition}
    \end{multline}
    will cause the system trajectory to asymptotically approach the set $\{x\in\reals^n \mid H(t,x) \in [-\alpha_w^{-1}(2), 0] \}$.
\end{theorem}

That is, by varying $\alpha_w$, we can tune how close to the boundary of $\mathcal{S}_H$ trajectories will approach. This theorem is put to use in~\cite{breeden2022docking} to achieve satisfaction of a so-called ``tight-tolerance'' objective with RCBFs. See also~\cite{alan2021safe,alan2023pbf}. The above result is also closely related to the definition of ``physical margin'' in Section~\ref{sec:sampled}, as robustness to unknown disturbances is closely related to robustness to inter-sampling uncertainty.

\subsubsection{Open Problems}

The above results and those in~\cite{breeden2021robust}, and the references therein, demonstrate that CBFs can be applied to relative-degree 2 systems with input constraints. However, there are still many systems that do not satisfy the conditions in the current literature and thus developing CBFs for these systems remains an open problem. Additionally, the above work only applies to one constraint function and CBF at a time. Working with multiple CBFs simultaneously and in the presence of input constraints is also an open problem. We refer the reader to~\cite{tan2022compatibility,breeden_acc23} for some constructive, albeit preliminary, methods towards this problem.

\subsection{Input-Constrained Control Barrier Functions}

As discussed, designing Control Barrier Functions for general nonlinear control-affine systems is challenging. When the system is also input constrained, this becomes further challenging, since there can be regions of the state-space where the CBF condition~\eqref{eq:cbf_online_condition} is instantaneously satisfied, but the system will evenutally reach the boundary of the safe set and exit then exit it. In~\cite{agrawal2021safe} we proposed a technique to isolate such states, and identify an inner safe set that can be rendered forward invariant under the input constraints. 

\begin{figure}
    \centering
    \includegraphics[width=0.4\linewidth]{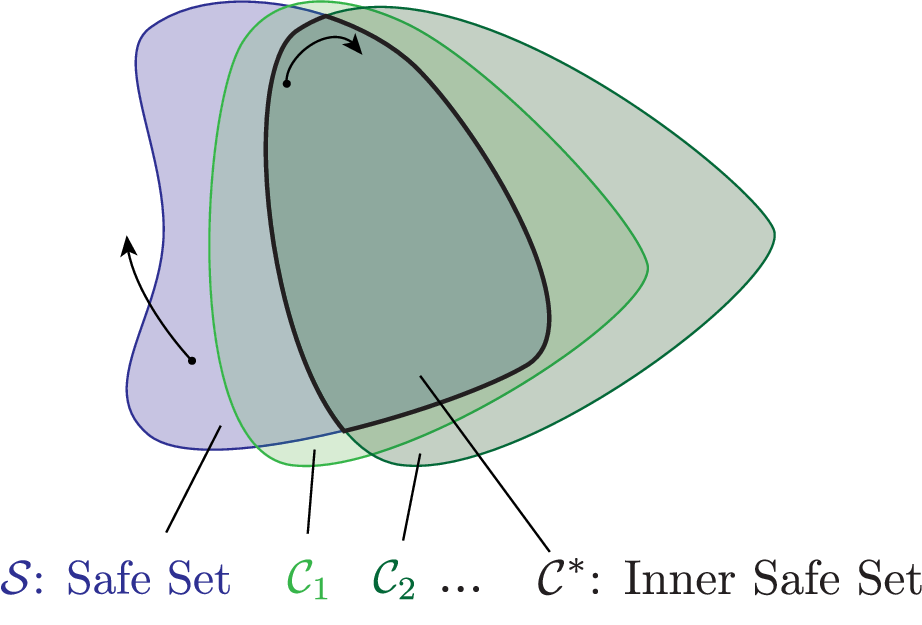}
    \caption{Visual representation of ICCBF method. The safe set $\mathcal S$ and two intermediate sets $\mathcal C_1$ and $\mathcal C_2$ are drawn. The final inner safe set $\mathcal C^*$ is the intersection of each of these sets, and can be rendered forward invariant. }
    \label{fig:SandC}
\end{figure}

Consider the dynamical system~\eqref{eq::dynamics_general} with bounded control inputs $u\in \mathcal U$ and a safe set $\mathcal S$ defined by a function $h: \mathcal X \rightarrow \reals$, as per~\eqref{eq::safeset}. Assume $h$ is not a CBF on $\mathcal S$. We define the following sequence of functions:
\begin{subequations}
\begin{align}
b_0(x) &= h(x) \label{eqn:b0def}\\
b_1(x) &= \textrm{inf}_{u\in \mathcal U} [L_f b_0(x) + L_g b_0(x) u + \alpha_0(b_0(x))] \\
b_2(x) &= \textrm{inf}_{u\in \mathcal U} [L_f b_1(x) + L_g b_1(x) u + \alpha_1(b_1(x))] \label{eqn:b1def}\\
\vdots \nonumber\\
b_N(x) &= \textrm{inf}_{u\in \mathcal U} [L_f b_{N-1}(x) + L_g b_{N-1}(x) u \nonumber\\&{}\quad{}\quad{}\quad{}\quad + \alpha_{N-1}(b_{N-1}(x))] \label{eqn:bNdef}
\end{align}
\end{subequations}
where each $\alpha_i$ is some user-specified class-$\mathcal K$ function, and $N$ is a positive integer. We assume the functions $f, g, h$ are sufficiently smooth such that $b_N$ and its derivative are defined. 
We also define the sets $\mathcal C_0 = \{ x \in \mathcal X : b_0(x) \geq 0\} = \mathcal S$, $\mathcal C_1 = \{ x \in \mathcal X : b_1(x) \geq 0\}$, $\dots$, 
$\mathcal C_N = \{ x \in \mathcal X : b_N(x) \geq 0\}$ and their intersection $\mathcal C^* = \mathcal C_0 \cap \mathcal C_1 \cap ... \cap \mathcal C_N$. We assume the set $\mathcal C^*$ is closed, non-empty and has no isolated points. The sets are visualized in Figure~\ref{fig:SandC}.

\begin{definition}
For the above construction, if there exists a class-$\mathcal K$ function $\alpha_N$ such that
\begin{align}
\sup_{u\in \mathcal U}[L_fb_N(x) + L_gb_N(x) u + \alpha_N(b_N(x))] \geq 0 \quad \forall x\in \mathcal C^*,\label{eqn:ICCBFdef}
\end{align}
then $b_N$ is an \emph{Input Constrained Control Barrier Function (ICCBF)}. 
\end{definition}

Note, this does not require $b_N$ to be a CBF on $\mathcal C_N$. The definition only requires condition~(\ref{eqn:ICCBFdef}) to hold for all $x \in \mathcal C^*$ which is a subset of $\mathcal C_N$.

The main result of~\cite{agrawal2021safe} is stated as:
\begin{theorem}
\label{theorem:main}
Given the input constrained dynamical system (\ref{eq::dynamics_general}), if $b_N$, defined by (\ref{eqn:b0def} - \ref{eqn:bNdef}), is an ICCBF, then any Lipschitz continuous controller $u : \mathcal C^* \rightarrow \mathcal U$ such that $u(x) \in K_{ICCBF}(x)$, where
$K_{ICCBF}(x) =\{ u \in \mathcal U: L_fb_N(x) + L_gb_N(x) u \geq -\alpha_N(b_N(x)) \}$
renders the set $\mathcal C^* \subseteq \mathcal S$ forward invariant.
\end{theorem}
\begin{remark}
\label{remark:HOCBFcomp}
Time-invariant Higher Order CBFs, as in~\cite{hocbf}, are a special case of ICCBFs. For instance, in systems of relative degree 2, $L_gh(x) = 0$ for all $x \in \mathcal S$. In this case, in the construction of ICCBFs we have $b_1(x) = \textrm{inf}_{u\in \mathcal U} [ L_fh(x) + L_gh(x) u +\alpha_0(h(x))]= \textrm{inf}_{u\in \mathcal U} [ L_fh(x)  +\alpha_0(h(x))]= L_fh(x)  + \alpha_0(h(x))$
which is exactly the function defined in~\cite{hocbf}. This repeats for any relative degree greater than 2, and thus for a system with relative degree $\rho$, the first $\rho$ expressions of ICCBFs are identical to those of HOCBFs. Moreover, ICCBFs can handle systems with non-uniform relative degree, by choosing $N$ greater or equal to the largest relative degree of the system in $\mathcal S$. 
\end{remark}

\subsubsection{Example: Adaptive Cruise Control}

\begin{figure*}
    \centering
    \includegraphics[width=\linewidth]{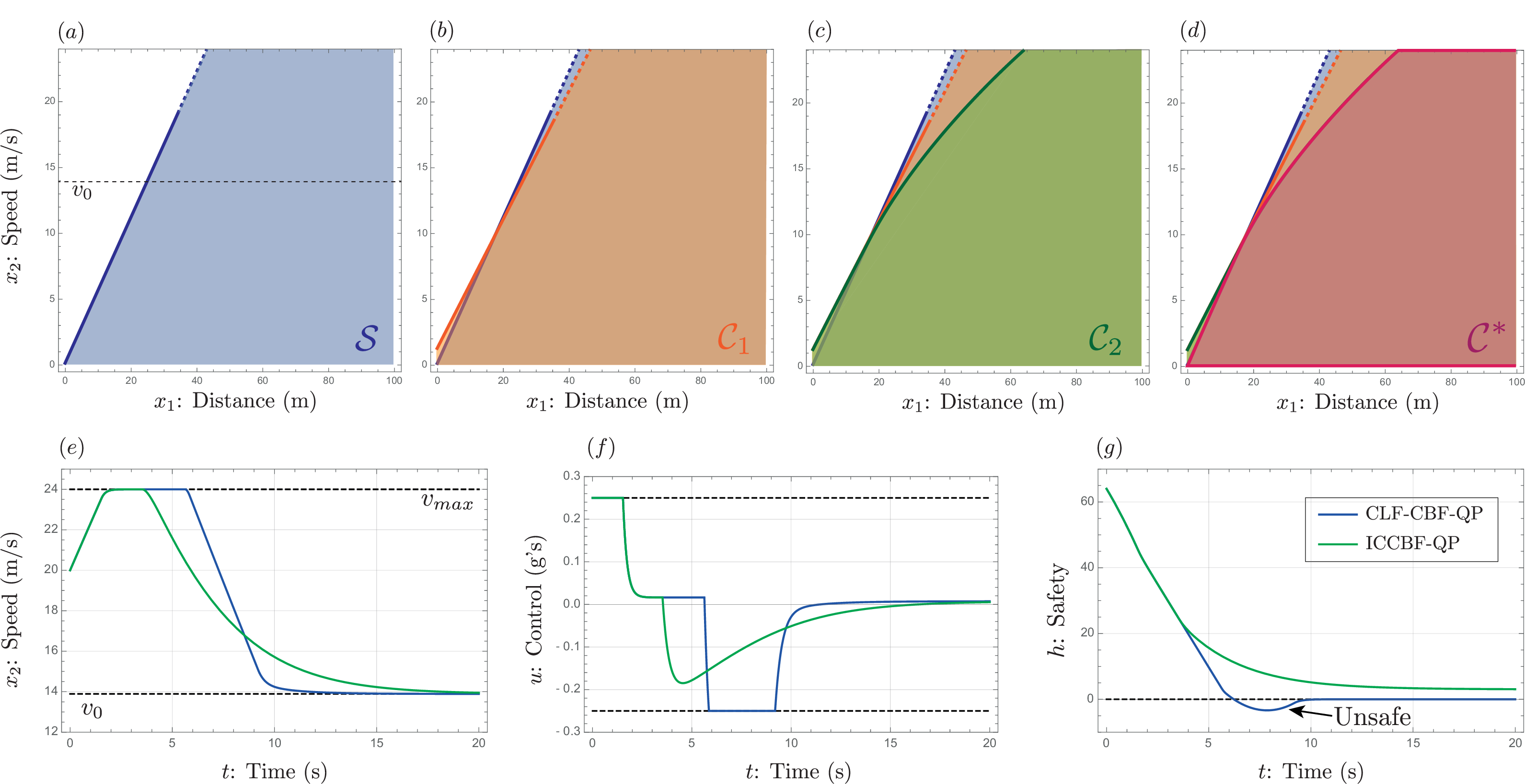}
    \caption{Figures (a-d): State-space diagrams indicating the sets (a)~$\mathcal S$, (b)~$\mathcal C_1$, (c)~$\mathcal C_2$ and (d)~$\mathcal C^*$. The horizontal dashed line in (a) indicates $v_0$, the speed of the car in-front. Figure~(d) represents the inner safe set $\mathcal C^*$ that is rendered forward invariant. Figures (e-g): Simulation results for speed, control input and safety under the CLF-CBF-QP controller~\cite{ames2014} and the ICCBF-QP.}
    \label{fig:sets_acc}
\end{figure*}

As a demonstration, we apply ICCBFs to the Adaptive Cruise Control (ACC) problem of~\cite{ames2014}. Consider a point-mass model of a vehicle moving in a straight line. The vehicle is following a vehicle $d$ distance in-front, moving at a known constant speed $v_0$. The objective is to design a controller to accelerate to the speed limit but prevent the vehicles from colliding. The dynamics model and safety constraints are as in~\cite{ames2014}, $\mathcal S = \{x \in \mathcal X : h(x) = x_1 - 1.8 x_2 \geq 0\}$. In addition, we impose the input constraints $\mathcal U = \{ u :  |u| \leq 0.25 \}$, representing a maximum acceleration or deceleration of 0.25~g. One can verify that $\mathcal S$ cannot be rendered forward invariant under the input constraints, and therefore we use the ICCBF construction technique to design an inner safe set. 

We (arbitrarily) choose $N = 2$ and the class~$\mathcal K$ functions $\alpha_0(h) = 4 h$,  $\alpha_1(h) = 7 \sqrt{h}$, $\alpha_2(h) = 2h$ to define the functions $b_1, b_2$ and sets $\mathcal C_1, \mathcal C_2$. To (approximately) verify that $b_2$ is an ICCBF, we used a nonlinear optimization to determine that~\eqref{eqn:ICCBFdef} was satisfied (see~\cite{agrawal2021safe} for details). The sets $\mathcal C_0,\mathcal  C_1, \mathcal C_2, \mathcal C^*$ are visualized in Figure~\ref{fig:sets_acc}.

Figure~\ref{fig:sets_acc}(e-g) compares the CLF-CBF-QP controller of~\cite{ames2014} (blue) to our proposed controller (green),
\begin{subequations}
\begin{align}
\pi(x) = \underset{u \in \mathbb \R}{\text{argmin }} & \frac{1}{2} (u-\pi_d(x))^2\\
\text{s.t. } & L_fb_2(x) + L_gb_2(x) u \geq - 2 b_2(x)\\
& u \in \mathcal U
\end{align}
\end{subequations}
where $\pi_d(x)$ is the desired acceleration, computed using a Control Lyapunov Function $V(x) = (x_2 - v_{max})^2$, where $v_{max}=24$ is the speed limit. 

The (standard) CLF-CBF-QP reaches the input-constraint at $t=5.9$~seconds, and thus the input constraints force the system to leave the safe set. In constrast, the ICCBF-QP remains feasible and safe for the entire duration, by applying brakes early, at $t=2.9$~seconds. Thus, by explicitly accounting for input constraints ICCBF-QP controller keeps the input-constrained system safe.

In the interest of space, the reader is referred to~\cite{agrawal2021safe} for additional details and examples on how to construct ICCBFs for relatively simple systems. How to systematically construct ICCBFs is part of our ongoing work.

\section{Adaptation: How to prevent loss of controllability, and how to reduce conservatism of the system response?}
\label{section::adaptation}
In this section, we present some results that involve online adaptation of CBFs towards two main challenges: The first for assuring that a candidate CBF will remain a valid CBF throughout the system trajectories, and the second for reducing the conservatism of the system response by allowing trajectories to approach closer to the boundary of the safety set. 

\subsection{Online Verification via Consolidated CBFs}

Verifying a \textit{candidate} CBFs as \textit{valid}, i.e., proving that the CBF condition is satisfiable via available control authority in perpetuity, is a challenging and rather underdeveloped problem. For isolated or single-CBF constraints, verifying or finding valid CBFs under either unlimited~\cite{ames2017control}, or bounded control authority~\cite{Breeden2021InputConstraints, Xiao2022Adaptive}, or by considering only one constraint at a time, either by assumption~\cite{Cortez2022RobustMultiple} or construction in a non-smooth manner~\cite{Glotfelter2017Nonsmooth,Huang2020SwitchedCBF} is a fairly studied task. However, these methods do not extend to multiple constraints. Some recent works synthesize and/or verify a CBF using sum-of-squares optimization~\cite{clark2021verification}, supervised machine learning~\cite{robey2020learning,srinivasan2020synthesis}, and Hamilton-Jacobi-Bellman reachability analysis~\cite{tonkens2022refining,yang2023safe}, but are limited to offline tools. 

In our recent work~\cite{black2022adaptation}, we consider a multi-agent system, each of whose $A$ constituent agents is modeled by the following class of nonlinear, control-affine dynamical systems:
\begin{equation}\label{eq.agent_dynamics}
    \dot{\bb{x}}_i = f_i(\bb{x}_i) + g_i(\bb{x}_i)\bb{u}_i,
\end{equation}
where $\bb{x}_i \in \mathbb R^n$ and $\bb{u}_i \in \mathcal{U}_i \subseteq \mathbb R^m$ are the state and control input vectors for the i$^{th}$ agent, with $\mathcal{U}_i$ the input constraint set, and where $f_i: \mathbb R^n \rightarrow \mathbb R^n$ and $g_i: \mathbb R^{n} \rightarrow \mathbb R^{n \times m}$ are known, locally Lipschitz, and not necessarily homogeneous $\forall i \in \mathcal{A} = [A]$. The concatenated state vector is $\bb{x} = [\bb{x}_1\^T,\hdots,\bb{x}_A\^T]\^T \in \mathbb R^N$, the concatenated control input vector is $\bb{u} = [\bb{u}_1\^T,\hdots,\bb{u}_A\^T]\^T \in \mathcal{U} \subseteq \mathbb R^M$, and as such the full system dynamics are
\begin{equation}\label{eq.multiagent_system}
    \dot{\bb{x}} = F(\bb{x}(t)) + G(\bb{x}(t))\bb{u}(\bb{x}(t)), \quad \bb{x}(0) = \bb{x}_0,
\end{equation}
where $F = [f_1\^T,\hdots,f_A\^T]\^T: \mathbb R^N \rightarrow \mathbb R^N$ and $G = \textrm{diag}([g_1,\hdots,g_A]): \mathbb R^{N} \rightarrow \mathbb R^{N \times M}$. 
Consider also a collection of $c>1$ state constraints, each described by a function $h_s \in \mathcal{C}^1: \mathbb R^N \rightarrow \mathbb R$ for $s \in [c]$. Each $h_s$ is a candidate CBF (hereafter referred to as a constituent constraint function) and defines a set
\begin{equation}\label{eq.multiple_constraint_sets}    \mathcal S_s = \{\bb{x} \in \mathbb R^N \mid h_s(\bb{x}) \geq 0\},
\end{equation}
that obeys the same structure as~\eqref{eq::safeset}. The following assumption is required, otherwise it is impossible to satisfy all constraints jointly.
\begin{assumption}\label{ass.nonempty_safe_set_intersection}
    The intersection of constraint sets is non-empty, i.e., $\mathcal S = \bigcap_{s=1}^c\mathcal S_s \neq \emptyset$.
\end{assumption}


\subsubsection{Definition of Consolidated CBFs}
Define a positive gain vector $\bb{k} = [k_1 \; \hdots \; k_c]\^T \in \mathbb R^c_+$. A consolidated CBF (C-CBF) candidate $H: \mathbb R^N \times \mathbb R_+^c \rightarrow \mathbb R$ takes the following form:
\begin{equation}\label{eq.consolidated_cbf}
    H(\bb{x}, \bb{k}) = 1 - \sum_{s=1}^c\phi\Big(h_s(\bb{x}), k_s\Big),
\end{equation}
where $\phi \in \mathcal{C}^1: \mathbb R_{+} \times \mathbb R_{+} \rightarrow \mathbb R_+$ belongs to class $\mathcal{L}\mathcal{L}$ and satisfies\footnote{For example, the decaying exponential function, i.e., $\phi(h_s,k_s)=e^{-h_sk_s}$, satisfies the requirements over the domain $\mathbb R_{+} \times \mathbb R_{+}$.} $\phi(h_s,0)=\phi(0,k_s)=\phi(0,0)=1$. It follows that the set $\mathcal{C}(\bb{k}) = \{\bb{x} \in \mathbb R^N \; | \; H(\bb{x}, \bb{k}) \geq 0\}$ is a subset of $\mathcal{S}$ (i.e., $\mathcal{C}(\bb{k}) \subset \mathcal{S}$), where the level of closeness of $\mathcal{C}(\bb{k})$ to $\mathcal{S}$ depends on the choices of gains $\bb{k}$. This may be confirmed by observing that if any $h_s(\bb{x}) = 0$ then $H(\bb{x}) \leq 1 - 1 - \sum_{j=1, j\neq s}^c \phi(h_j(\bb{x}), k_j) < 0$, and thus for $H(\bb{x}) \geq 0$ it must hold that $h_s(\bb{x}) > 0$, for all $s \in [c]$.

Now, if $H$ is a \textit{valid} C-CBF over the set $\mathcal{C}(\bb{k})$, then $\mathcal{C}(\bb{k})$ is forward invariant and thus the trajectories of~\eqref{eq.multiagent_system} remain safe with respect to each constituent safe set $\mathcal{S}_s$, $\forall s \in [c]$. For a static gain vector (i.e., $\dot{\bb{k}} = \mathbf{0}_{c \times 1}$) the function $H$ is a CBF on the set $\mathcal{S}$ if there exists $\alpha_H \in \mathcal{K}_\infty$ such that the following condition holds for all $\bb{x} \in \mathcal{S} \supset \mathcal{C}(\bb{k})$:
\begin{equation}\label{eq.ccbf_condition_static_k}
    L_FH(\bb{x}, \bb{k}) + L_GH(\bb{x}, \bb{k})\bb{u}(\bb{x}) \geq -\alpha_H(H(\bb{x}, \bb{k})),
\end{equation}
where from~\eqref{eq.consolidated_cbf} it follows that
\begin{align}
    L_FH(\bb{x}) &= -\sum_{s=1}^c\frac{\partial \phi}{\partial h_s}L_Fh_s(\bb{x}), \label{eq.LfH_static_k} \\
    L_GH(\bb{x}) &= -\sum_{s=1}^c\frac{\partial \phi}{\partial h_s}L_Gh_s(\bb{x}) \label{eq.LgH_static_k}.
\end{align}
Again taking $\phi(h_s, k_s) = e^{-h_sk_s}$ as an example, we obtain that $\frac{\partial \phi}{\partial h_s} = -k_se^{-h_sk_s}$, in which case it is evident that the role of the gain vector $\bb{k}$ is to weight the constituent constraint functions $h_s$ and their derivative terms $L_Fh_s$ and $L_Gh_s$ in the CBF condition~\eqref{eq.ccbf_condition_static_k}. In this case, a higher value $k_s$ indicates a weaker weight in the CBF dynamics, as the exponential decay overpowers the linear growth. Due to the combinatorial nature of these gains, for an arbitrary $\bb{k}$ there may exist some $\bb{x} \in \mathcal{C}(\bb{k})$ such that $L_GH(\bb{x}) = \mathbf{0}_{1\times M}$, which lead to the state exiting $\mathcal{C}(\bb{k})$ (and potentially $\mathcal{S}$ as a result). 
Using online adaptation of $\bb{k}$, however, it may be possible to achieve $L_GH(\bb{x}(t)) \neq \mathbf{0}_{1\times M}$ for all $t \geq 0$, which motivates the following problem.
\begin{problem}\label{prob.adaptation}
    Given a C-CBF candidate $H: \mathbb R^N \times \mathbb R_+^c \rightarrow \mathbb R$ defined by~\eqref{eq.consolidated_cbf}, design an adaptation law $\dot{\bb{k}} = \kappa(\bb{x}, \bb{k})$ such that $L_GH(\bb{x}(t)) \neq \mathbf{0}_{1\times M}$, $\forall t \geq 0$.
\end{problem}


\subsubsection{C-CBF Weight Adaptation Law}
\begin{assumption}\label{ass.LgH_nonzero}
    Let the intersection of constraint sets be denoted $\mathcal{S}$; the matrix of controlled constituent function dynamics $\bb{L}_G \in \mathbb R^{c \times M}$ is not all zero, i.e.,
    \begin{equation}\label{eq.G_matrix}
        \bb{L}_G(\bb{x}) \triangleq \begin{bmatrix} L_Gh_1(\bb{x}) \\ \vdots \\ L_Gh_c(\bb{x}) \end{bmatrix} \neq \mathbf{0}_{c \times M}, \; \forall \bb{x} \in \mathcal{S}. 
    \end{equation}
\end{assumption}
The above requires non-zero sensitivity of at least one constraint function $h_s$ to the control input $\bb{u}$. It is a mild condition, and is easily satisfiable when at least one $h_s$ is of relative-degree one with respect to the system~\eqref{eq.multiagent_system}. In what follows, it is shown that the ensuing QP-based adaptation law renders a C-CBF as valid.
\\
\textbf{(C-CBF-QP)}
\begin{subequations}\label{subeq.k_dot}
    \begin{align}
        \kappa(\bb{x}, \bb{k}) &= \argmin_{\bb{\mu} \in \mathbb R^c} \; \frac{1}{2}(\bb{\mu} - \bb{\mu}_0(\bb{x}))\^T\bb{P}(\bb{\mu} - \bb{\mu}_0(\bb{x})) \\
        \mathrm{s.t.}& \quad
        \mu_s + \alpha_k(k_s-k_{s,min}) \geq 0, \; \forall s \in [c], \label{subeq.kdot.k_gr0} \\
        &\quad \bb{p}\^T(\bb{x})\bb{Q}(\bb{x})\dot{\bb{p}} + \frac{1}{2}\bb{p}\^T(\bb{x})\dot{\bb{Q}}\bb{p}(\bb{x}) + \alpha_{p}(h_{p}(\bb{x})) \geq 0, \label{subeq.kdot.vperp}
    \end{align}
\end{subequations}
where $\bb{P} \in \mathbb R^{c \times c}$ is a positive-definite gain matrix, $\alpha_k, \alpha_p \in \mathcal{K}_\infty$, $\bb{\mu}_0 \in \mathbb R^c$ is the desired solution, $\bb{k}_{min}= [k_{1,min},\hdots,k_{c,min}]\^T$ is the vector of minimum allowable values $k_{s,min} > 0$, and $\bb{p}(\bb{x}) \triangleq \left[\frac{\partial \phi(\bb{x})}{\partial h_1} \; \hdots \; \frac{\partial \phi(\bb{x})}{\partial h_c}\right]^\top, 
    \bb{Q}(\bb{x}) \triangleq \bb{I} - 2\bb{N}\bb{N}\^T - \bb{N}\bb{N}\^T\bb{N}\bb{N}\^T \label{eq.Q_matrix}$ with $h_p(\bb{x}) = \frac{1}{2}\bb{p}\^T(\bb{x})\bb{Q}(\bb{x})\bb{p}(\bb{x}) - \varepsilon$, $\varepsilon > 0$, and
    $\bb{N} = \bb{N}(\bb{x}) \triangleq [\bb{n}_1(\bb{x}) \; \hdots \; \bb{n}_r(\bb{x})]$,
such that 
$\{\bb{n}_1(\bb{x}),\hdots,\bb{n}_r(\bb{x})\}$ constitutes a basis for the null space of $\bb{L}_G\^T(\bb{x})$, i.e., $\mathcal{N}(\bb{L}_G\^T(\bb{x})) = \mathrm{span}\{\bb{n}_1(\bb{x}),\hdots,\bb{n}_r(\bb{x})\}$, where $\bb{L}_G(\bb{x})$ is given by~\eqref{eq.G_matrix}.
\begin{theorem}\label{thm.kdot_HCBF}
    Suppose that there exist $c>1$ constraint functions $h_s: \mathbb R^N \rightarrow \mathbb R$ defining sets $\mathcal S_s = \{\bb{x} \in \mathbb R^N \; | \; h_s(\bb{x}) \geq 0\}$, $\forall s \in [c]$, that Assumptions \ref{ass.nonempty_safe_set_intersection} and \ref{ass.LgH_nonzero} hold, and that $\mathcal{U} = \mathbb R^M$. If $\bb{k}(0)$ is such that $L_GH(\bb{x}(0)) \neq \mathbf{0}_{1 \times M}$, then under $\dot{\bb{k}} = \kappa(\bb{x},\bb{k})$ the controlled C-CBF dynamics are non-vanishing provided that~\eqref{subeq.k_dot} is feasible, i.e., $L_GH(\bb{x}(t)) \neq \mathbf{0}_{1 \times M}$, $\forall t \geq 0$.
\end{theorem}
It follows that if the premises of Theorem \ref{thm.kdot_HCBF} hold, then the function $H$ defined by~\eqref{eq.consolidated_cbf} is a CBF for $\mathcal{C}(\bb{k}(t)) = \{\bb{x} \in \mathbb R^N \; | \; H(\bb{x}, \bb{k}) \geq 0\}$, $\forall t \geq 0$.

\subsubsection{Simulation Results: Multi-Robot Coordination}

Consider 3 non-communicative, but responsive robots ($i \in \mathcal{A}_{nc} \setminus \mathcal{A}_{ncnr}$) in a warehouse environment seeking to traverse a narrow corridor intersected by a passageway occupied with 6 non-responsive agents ($j \in \mathcal{A}_{ncnr}$). The non-responsive agents may be, e.g., humans walking or some other dynamic obstacles. 
Each robot is modeled according to a kinematic bicycle model described by~\cite[Ch. 2]{Rajamani2012VDC}:
which we omit here in the interest of space; the reader can refer to the complete case study in~\cite{black2022adaptation}.
The safety constraints for each robot are to 1) obey the velocity restriction, 
2) remain inside the corridor area 
and 3) avoid collisions with all other robots. 
As such, each robot has three constituent constraint sets (defined explicitly in~\cite{black2022adaptation}),
the intersection of which constitutes the safe set for each robot. 
The robots $i \in \mathcal{A}_{nc} \setminus \mathcal{A}_{ncnr}$ are controlled using a C-CBF based decentralized controller with constituent functions $h_c$, $h_s$, $h_r$, an LQR based nominal control input (see~\cite[Appendix 1]{black2023future}), and initial gains $\bb{k}(0) = \mathbf{1}_{10 \times 1}$. The non-responsive agents used a similar LQR  controller to move through the passageway in pairs of two, with the first two pairs passing through the intersection without stopping and the last pair stopping at the intersection before proceeding. 

\begin{figure}[!ht]
    \centering
        \includegraphics[clip,width=0.9\linewidth]{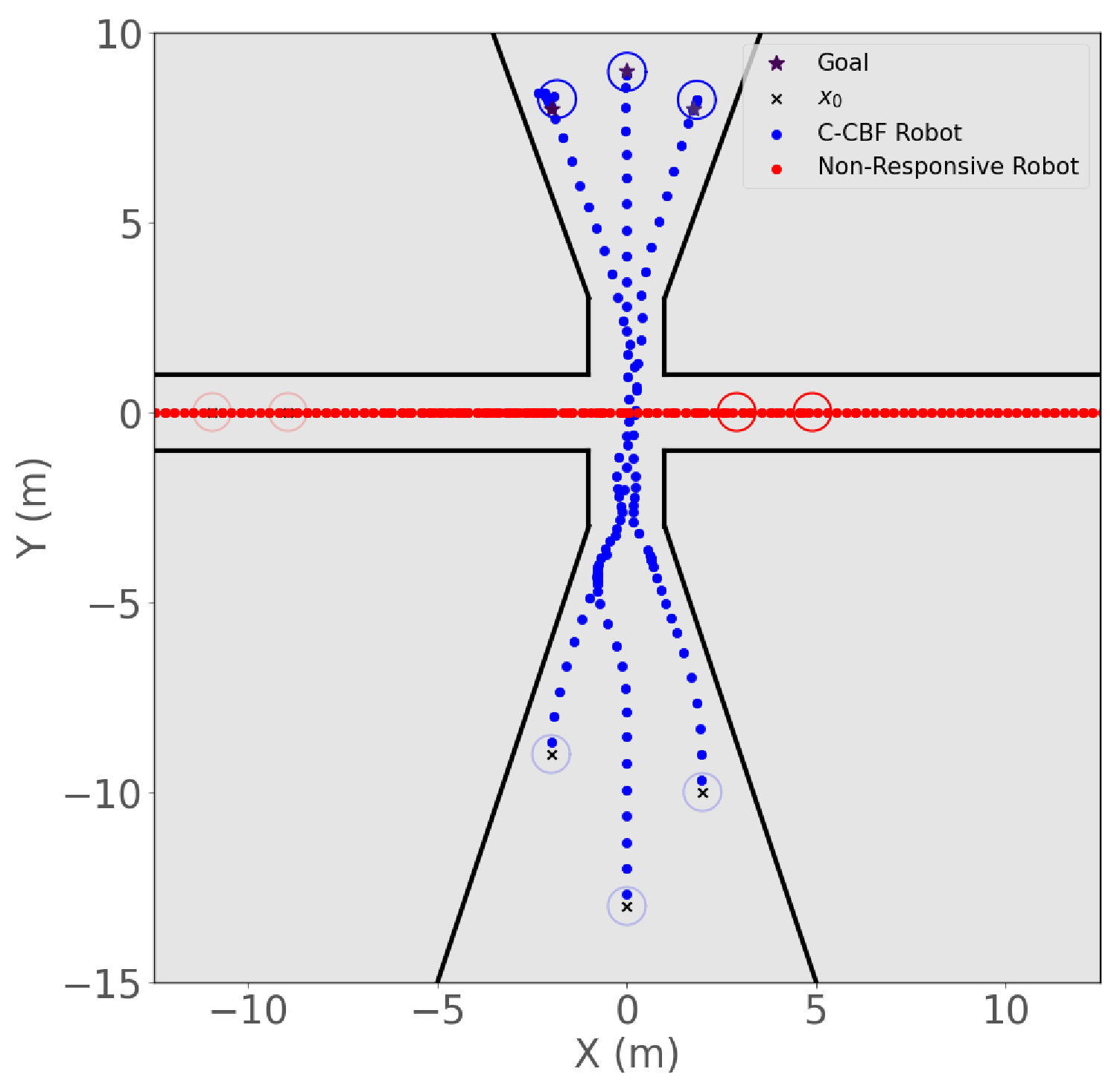}
    \caption{\small{XY paths for the warehouse robots (blue) and non-responsive agents (red) in the warehouse control problem.}}\label{fig.warehouse_robot_trajectories}
    \vspace{-3mm}
\end{figure}

\begin{figure}[!ht]
    \centering
    \includegraphics[clip,width=0.95\columnwidth]{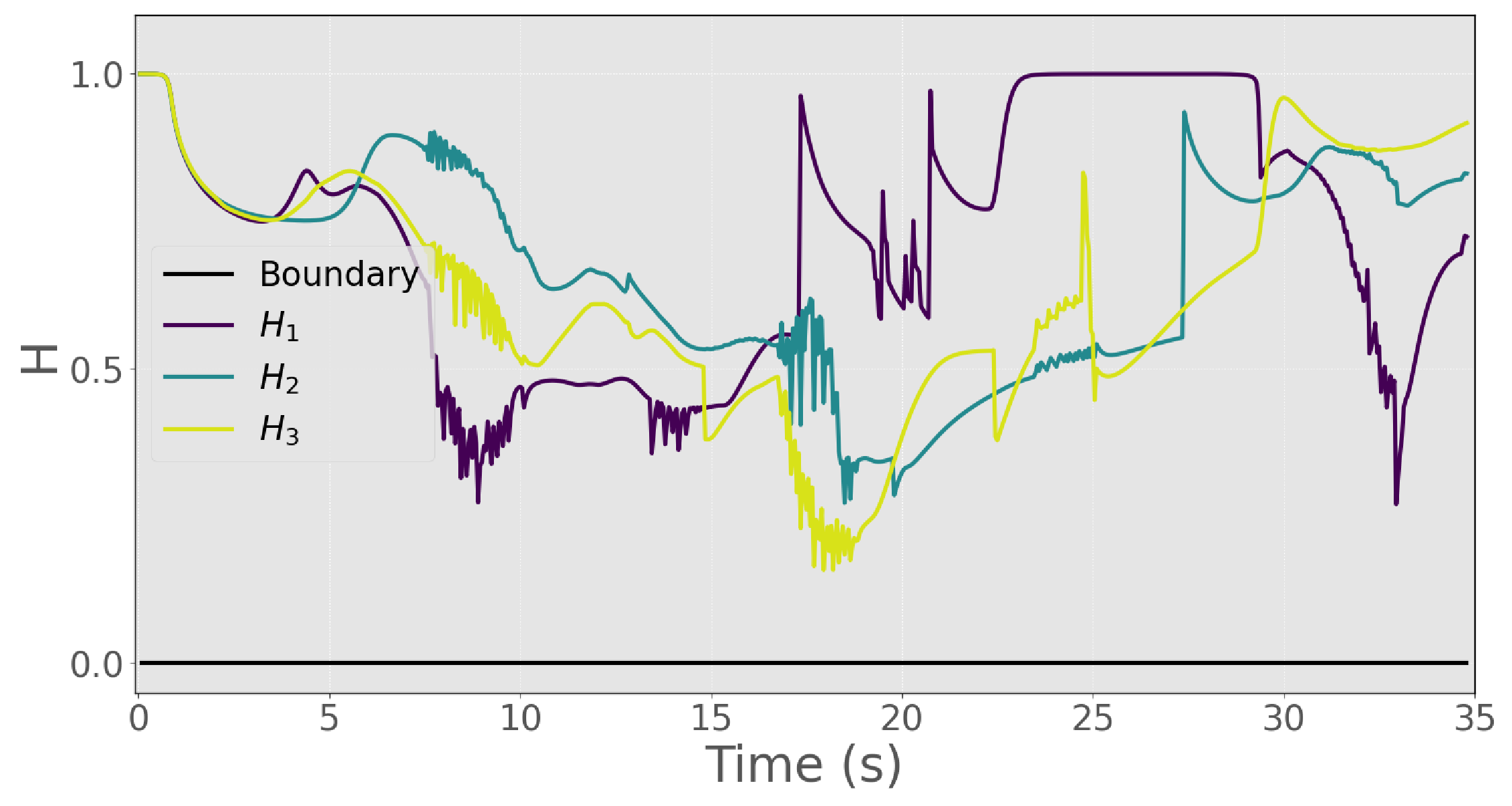}
    \caption{\small{Evolution of C-CBF $H$ for warehouse robots 1, 2, and 3.}}\label{fig.warehouse_ccbfs}
    \vspace{-1mm}
\end{figure}

As shown in Figures \ref{fig.warehouse_robot_trajectories} and \ref{fig.warehouse_ccbfs}, the non-communicative robots traverse both the narrow corridor and the busy intersection to reach their goal locations safely. These results demonstrate that the C-CBF-based adaptive controllers maintained safety and control viability at all times amongst 10 state constraints.

\subsection{Parameter Adaptation with Rate Tunable CBFs}
In a recent parallel thread of work, instead of designing adaptive laws for the control coefficient $L_gh_i(t,x)$, we consider the adaptation of the parameters introduced through the class-$\mathcal K$ function in the CBF condition. In~\cite{parwana2022trust} we introduce a new notion of a Rate-Tunable Control Barrier Function (RT-CBF), which allows consideration of parametric class $\mathcal K$ functions, and adaptation of their parameters online so that the response of the controller becomes less or more conservative, without jeopardizing safety. It is also noteworthy that this adaptation facilitates the consideration and satisfaction of multiple time-varying barrier constraints, by making them easier to tune for performance, especially when they do not represent similar physical quantities (e.g., when imposing constraints on the rotational dynamics, and constraints on the translational dynamics for a quadrotor). 

Designing the parameter dynamics is a non-trivial task, especially in the presence of multiple constraints. We have studied pointwise sufficient conditions on the rate of change of parameters for enforcing feasibility. Although the pointwise design is suboptimal, it is shown empirically to improve upon the standard CBF-QP controllers. It also allows the incorporation of user-designed rules (e.g., heuristics) for updating the class $\mathcal K$ function and project it to a set of feasible update rules. As a case study, we design RT-CBFs for decentralized control for multi-robot systems in~\cite{parwana2022trust}. Specifically, we design the parameter dynamics based on a trust factor, which in turn is defined on the instantaneous ease of satisfaction of the CBF constraints, and illustrate how this can be applied to robots of heterogeneous dynamics.

\subsubsection{Definition of Rate-Tunable CBFs}

For ease of understanding, we illustrate our theory with examples that only consider linear class $\mathcal K$ functions of the following form in the ensuing. 
\begin{align}
    \alpha_i^k(z) = \nu_i^k z, \quad \nu_i^k \in \reals^+.
\end{align}

Since we allow parameters $\nu_i^k$ to vary with time, the derived barrier functions in~\eqref{eq::hocbf_derived_barriers} for, for example, a second-order barrier function are given as follows
\begin{subequations}
    \begin{align}
            \psi_i^0 &= h_i, \\
\psi_i^1 &= \dot \psi^0 + \nu_i^1 \psi_i^0, \\
\psi_i^2 &= \ddot \psi^0 + \dot \nu_i^1 \psi_i^0 + \nu_i^1 \dot \psi_i^0 + \nu_i^2( \psi_i^0 + \nu_i^1 \psi_i^0), \label{eq::derived_barriers_example_2}
    \end{align}
\label{eq::derived_barriers_example}
\end{subequations}
where $\psi_i^2 = \dot \psi_i^1 + \nu_i^2 \psi_i^1\geq 0$ is the CBF condition~\eqref{eq::qp_hocbf_constraint} that is used to design the control input. We denote the parameters and their derivatives contributing to the CBF condition~\eqref{eq::derived_barriers_example_2} as $\Theta_{\alpha_i} \in \reals^{n_\alpha}$ and the objective is to design $\dot \Theta_{\alpha_i}$. 
For example, for~\eqref{eq::derived_barriers_example}, $\Theta_{\alpha_i^k}=[\nu_i^1, \dot \nu_i^1, \nu_i^2]^T$ and 
\begin{align}
\dot \Theta_{\alpha_i} = \begin{bmatrix}
    \dot \nu_i^1 \\
    \ddot \nu_i^1 \\
    \dot \nu_i^2
\end{bmatrix} = \begin{bmatrix}
    0 & 1 & 0 \\ 
    0 & 0 & 0 \\
    0 & 0 & 0 \\
\end{bmatrix}\Theta_{\alpha_i} + \begin{bmatrix}
    0 & 0 \\
    1 & 0 \\
    0 & 1
\end{bmatrix}\begin{bmatrix}
    \ddot \nu_i^1 \\
    \dot \nu_i^2
\end{bmatrix}
\label{eq::Theta_dynamics}
\end{align}
Note that the derivatives to be designed, namely $\ddot \nu_i^1, \dot \nu_i^2$, do not appear in the CBF condition~\eqref{eq::derived_barriers_example_2} that is imposed in QP to design the control input. This allows for decoupling the design of control input and the parameter dynamics.

Consider the system dynamics in (\ref{eq::dynamics_general}) augmented with the state $\Theta_\alpha\in\reals^{n_\alpha}$ that obeys the dynamics
\begin{align}
\begin{bmatrix} \dot x \\ \dot \Theta_{\alpha_i} \end{bmatrix} = \begin{bmatrix} f(x) + g(x) u \\ f_{\alpha_i}(x, \Theta_{\alpha_i}) \end{bmatrix},
\label{eq::augmented}
\end{align}
where $f_{\alpha_i}:\mathcal{X}\times \mathcal A_i \rightarrow \mathcal A_i$ is a locally Lipschitz continuous function w.r.t. $(x,\Theta_{\alpha_i})$, where $\mathcal{A}_i
\subset \reals^{n_{\alpha_i}}$ is a compact set, $i \in \{1,\dots,N\}$.

Assume that $\Theta_{\alpha_i}(t)\in \mathcal{A}_i$, $\forall t\geq 0$, where $\mathcal{A}_i
\subset \reals^{n_{\alpha_i}}$ is a compact set. 
Let the set of allowable states $\s_i(t)$ at time $t$ be defined as the 0-superlevel set of a function $h_i:\mathbb R^+\times\mathcal{X}\rightarrow \reals$ as in~\eqref{eq::safeset}. Suppose $h_i$ has relative degree $r_i$ w.r.t the control input $u$ and define functions $\psi^k_i$ as: 
\begin{subequations}
\begin{align}
   \psi^0_i(t,x) &= h_i(t,x), \\
   \psi^k_i(t,x) &= \dot \psi^{k-1}_i(t,x,\Theta_{\alpha_i^k}) + \alpha^k_i (\psi^{k-1}_i(t,x,\Theta_{\alpha_i^k})),  \label{eq::rt_hobcf_barriers} \\
   & \quad \quad \quad \quad \quad k\in \{1,2,..,r_i-1\}.\nonumber
\end{align}
\end{subequations}

\begin{definition}(Rate-Tunable CBF) 
A (single) constraint function $h_i:\mathbb R^+ \times \mathcal X\rightarrow \reals$ is a Rate-Tunable Control Barrier Function (RT-CBF) for the set $\mathcal{S}_i(t)$ under the augmented system~\eqref{eq::augmented} if for every initial state $x(0)\in \mathcal{S}_i(0)$, there exists  $\Theta_{\alpha_i}(0)\in \mathcal{A}_i$ 
such that $\forall t\geq 0$  
\begin{align}
\sup_{u\in \mathcal{U}} \left[
  \dot \psi_i^{r_i-1}(t,x,u, \Theta_{\alpha_i}) + \alpha_i^{r_i} (\psi_i^{r_i-1}(t,x,\Theta_{\alpha_i})) 
  \right] \geq 0.
\label{eq::rtcbf_derivative_condition}
\end{align}
\label{definition::RT-CBF}
\end{definition}
Note that for $\dot \Theta_{\alpha_i}\equiv0$ and $\Theta_{\alpha_i}(0)\equiv \Theta_{\alpha_i}$ (a constant independent of $x(0)$), we recover the definition of the classical CBF. In that regard, RT-CBF is a weaker notion of a classical CBF, which allows for tuning the response of the system.
    Note \eqref{eq::rtcbf_derivative_condition} is required to be satisfied for all $t\geq 0$ and not for all $x\in \mathcal{X}$ as required in vanilla CBF \eqref{eq::cbf_derivative} and HOCBF \eqref{eq::hobcf_cbf_derivatrive} conditions. This difference is essential as we allow for the initial parameter value $\Theta_{\alpha_i}(0)$ that is dependent on the initial state $x(0)$.

\begin{remark}
While several works employ heuristics to tune the parameters of the CBF condition~\eqref{eq::cbf_derivative} so that a solution to the CBF-QP~\eqref{eq::HOCBF-CLF-QP} exists for all $t>0$~\cite{zeng2021safety,wang2022ensuring,garg2019control}, most of these are equivalent to treating the parameter $\nu$ of a linear $\classK$ function $\alpha(h) = \nu h$ as an optimization variable. However, a formal analysis encompassing all these heuristics and other possible ways to adapt the $\classK$ function has been lacking so far, and RT-CBFs aim to bridge this gap in theory and application.
\end{remark}

In~\cite{parwana2023tac} we show that under mild assumptions (existence and uniqueness of the system trajectories), the existence of a RT-CBF is a necessary and sufficient condition for safety. We also show that several existing parameter adaptation schemes~\cite{ma2022learning, parwana2022recursive, zeng2021enhancing,garg2019control,garg2022control} fall under the framework of the proposed RT-CBFs.
We then state the following theorem that illustrates how the tuning of the parameter $\Theta$ can be used to shape the response of the CBF-QP controller.

\begin{theorem}\cite{parwana2023tac}
\label{theorem::alpha_minimum_norm}
    Consider the system~\eqref{eq::dynamics_general}, a first-order candidate barrier function $h$, a function $\alpha(h,x): \reals\times\reals^n \rightarrow \reals$ and the following CBF-QP controller with unbounded control input
    \begin{subequations}
    \begin{align}
        u_{QP} = \min_{u\in \reals^m} \quad &  ||u-u_{r}(t,x)||^2 \\
        \textrm{s.t.} \quad &  \dot h(t,x,u) \geq -\alpha(x,h(t,x))
    \end{align}
        \label{eq::cbf_qp_simple}
\end{subequations}
where $u_{r}: \reals^+ \times \reals^n\rightarrow $ is the reference (nominal) control input. Let $u_d(t,x):\reals^+ \times \reals^n\rightarrow \mathcal{U}$ be any desired safe response of the system and $u_r\neq u_d$ w.l.o.g. Then the following choice of function $\alpha(\cdot, \cdot)$ minimizes the norm $||u_d - u_{QP}||$ 
\begin{align}
\alpha(x,h) \!= \!\left\{  
      \begin{array}{ll}
         \!\sqrt{2}\frac{L_gh}{||L_gh||}(u_r\! -\! u_d)  &  \textrm{if } ||L_gh||>0\\ 
         \!\quad  - L_fh - L_ghu_r & \textrm{     and } L_gh u_d\!>\!L_gh u_r\\ \\
         \! - L_fh - L_ghu_r & \textrm{if } ||L_gh||=0 \\ \\
          \!\sqrt{2}\frac{L_gh}{||L_gh||}(u_r\! -\! u')  &  \textrm{if } ||L_gh||>0 \\
          \!\quad   - L_fh - L_ghu_r &
          \textrm{     and } L_gh u_d\!<\!L_gh u_r
      \end{array}
\right.   
\label{eq::alpha_minimum_norm}
\end{align}
where 
\begin{align}
u' = u_d + \frac{1}{\sqrt{2}}L_gh^T L_gh (u_r-u_d).
\label{eq::u_special_case}
\end{align}
\end{theorem}

The result of Theorem \ref{theorem::alpha_minimum_norm} albeit simple gives us some important insights. First, for different desired responses $u_d$ (such as conservative or aggressive) at state $x$ and time $t$, the function $\alpha$ can be used to steer $u_{QP}$ close to $u_d$. Second, to achieve the aforementioned steering, the function $\alpha$ cannot be just a $\classK$ of the barrier function $h$ as~\eqref{eq::alpha_minimum_norm} depends not only on $h$ but also on $x$. In our framework of RT-CBF, the parameter $\Theta$ is a function of $t,x,h$ and thus can fulfill this objective at the points where $\alpha$ in~\eqref{eq::alpha_minimum_norm} is differentiable. 
We consider the following RT-CBF-QP controller with parametric class $\mathcal K$ functions
\\
\textbf{(RT-CBF-QP)}
\begin{subequations}
        \begin{align}
     u = \arg \min_{u\in \mathcal{U}} \quad &  ||u-u_r||^2 + M \delta^2 \\
         \textrm{s.t.} \quad & \dot{V}(t,x,u) \leq - k V(t,x) + \delta, \label{eq::rt_qp_hoclf_constraint} \\
           \dot \psi_i^{r_i-1}(t,x,u, \Theta_{\alpha_i}) &+ \alpha_i^{r_i} (\psi_i^{r_i-1}(t,x,\Theta_{\alpha_i})) \geq 0, \\ 
          & \quad \quad \quad \quad i\in \{1,2,..,N\}\nonumber. \label{eq::qp_rt_hocbf_constraint} \;
        \end{align}
        \label{eq::RT-HOCBF-CLF-QP}
\end{subequations}

While we have established the necessity of RT-CBFs, finding a valid update law $\dot \Theta$ for parameters, much like finding a valid CBF in the sense of~\eqref{eq::cbf_derivative}, is non-trivial. In~\cite{parwana2023tac} we present some suboptimal and heuristic methods for ensuring that~\eqref{eq::RT-HOCBF-CLF-QP} admits a solution for all time. In the interest of space we omit the presentation of the algorithm from the current paper, and present directly some illustrative results.

\subsubsection{Simulation Results: Adaptive Cruise Control}
We simulate the Adaptive Cruise Control (ACC) problem with a decelerating leader. Let the ego agent move with velocity $v$ and a leader agent move at velocity $v_L$ with distance $D$ between them. The dynamics are 
\begin{align}
\begin{bmatrix}
    \dot v \\
    \dot v_L \\ 
    \dot D
\end{bmatrix} = \begin{bmatrix}
    \frac{-1}{M}F_r(v) + \frac{1}{M} u \\
    a_L \\
    v_L - v
\end{bmatrix}
\end{align}
where $M$ denotes the mass of ego agent, $F_r(v)$ is the resistance force 
and $a_L$ is the acceleration of the leader. The safety objective of the ego agent is to maintain a minimum distance $\bar D$ with the leader. The desired velocity is given by $v_d$. Additionally, the velocity and control input are constrained by $v_{min} \leq v \leq v_{max}$, $-c M g \leq u \leq c M g$, where $c$ is acceleration coefficient and $g$ is acceleration due to gravity. To ensure safety, we formulate three barrier functions $h_1 = D - \bar D$, $h_2 = v - v_{min}$, $h_3 = v_{max}-v$, where $h_1$ is second-order barriers and $h_2, h_3$ are first-order barrier. The CLF is chosen as $V = (v-v_d)^2$. We further choose linear \classK functions and apply the Algorithm in~\cite{parwana2023tac} to ensure safety.\footnote{The following parameters are used for the simulation: $M=1650$ kg, $f_0=0.1, f_1=5, f_2=0.25$, $g=9.81m/s^2$, $v_{mijn}=0$, $v_{max}=0$, $v_d=24m/s$, $c=0.4$, $\bar D=10$.} 

\begin{figure}[h!]
\centering
         \includegraphics[scale=0.5]{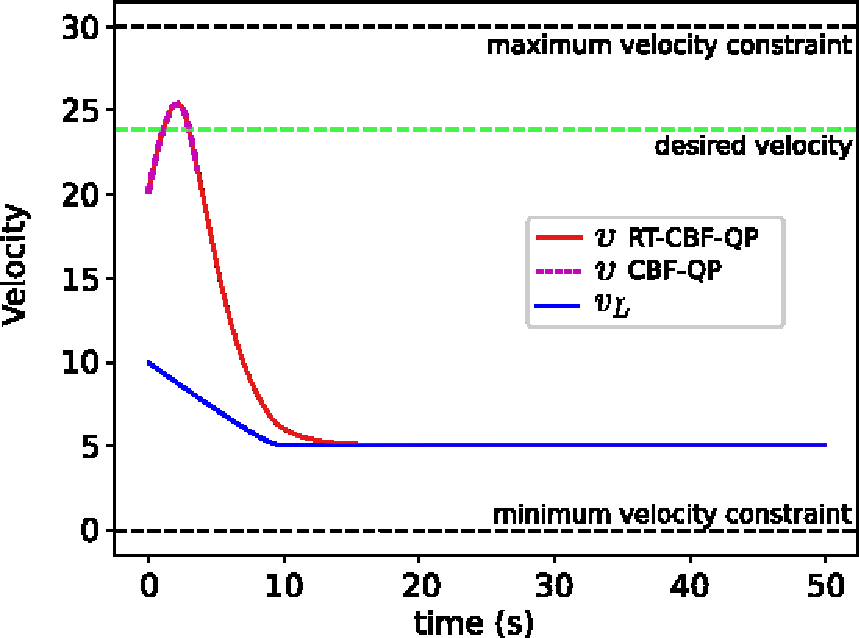}
         \caption{\small{Velocity variation for the ACC problem. The dotted pink line ends when CBF-QP becomes infeasible. The leader decelerates for the first 10 seconds and then moves at constant velocity. The ego-agent approaches the boundary of the safe set and hence does not have enough bandwidth to steer back towards its desired velocity of $24m/s$. }}
         \label{fig::acc_vel}
\end{figure}

\begin{figure}[h!]
\centering
         \includegraphics[scale=0.5]{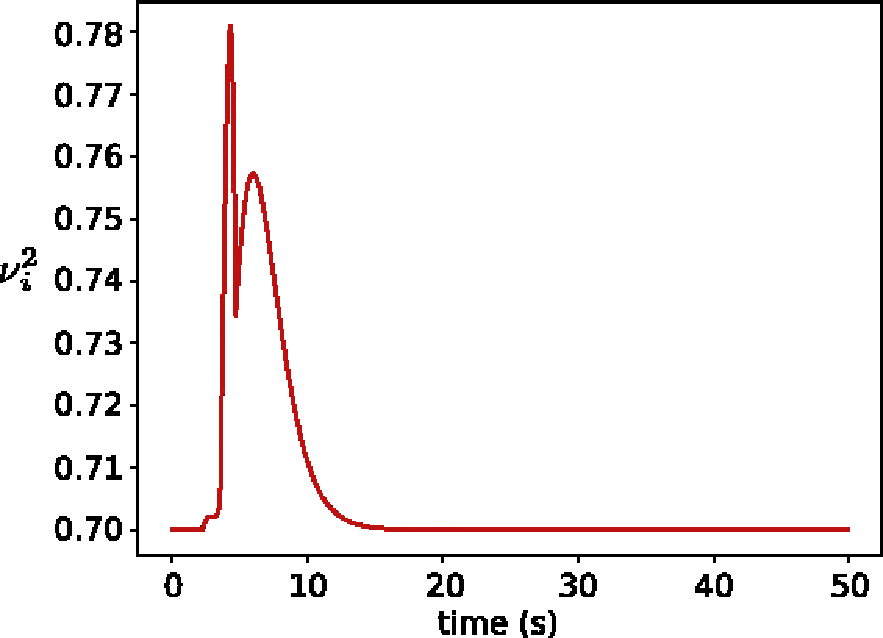}
         \caption{\small{The variation of parameter $\nu_1^2$ with time for the ACC problem. $\nu_i^1$ starts increasing as control bounds are approached. It is brought back to its nominal value of $0.7$ eventually}}.
         \label{fig::acc_param}
\end{figure}

\begin{figure}[h!]
\centering
         \includegraphics[scale=0.5]{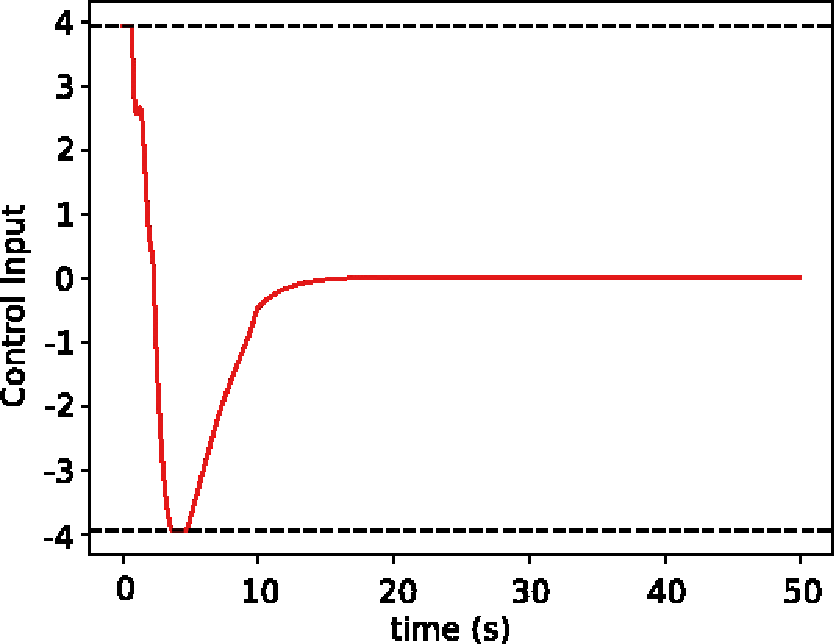}
         \caption{\small{The variation of normalized control input with time for the ACC problem. Normalization is done as $v/M^2$. }}
         \label{fig::acc_control}
\end{figure}

The objective of the CBF-QP to ACC is designed in the same way as~\cite{ames2017control}. We compare CBF-QP with RT-CBF-QP in Fig. \ref{fig::acc_vel} (for more comparisons please refer to~\cite{parwana2023tac}). 
The simulation is run for the same initial  parameters in both cases, for $50$ seconds. The CBF-QP (with fixed parameters) becomes infeasible before the simulation finishes as shown in Fig. \ref{fig::acc_vel}. The RT-CBF-QP on the other hand (with adapted parameters) can maintain feasibility and safety at all times. 
The variations of $\nu_1^2$, and $u$ with time are shown in Figs. \ref{fig::acc_param} and \ref{fig::acc_control} respectively. Note that $\nu_1^2$ starts increasing as the control input bound is approached. 
Other parameters do not change in this example and their variation is thus not shown
This example illustrates that, for a chosen barrier function $h$, the feasibility of CBF-QP controllers is highly dependent on parameters, but online adaptation can help circumvent this issue.

\section{Prediction: How to reduce myopic behavior?}
\label{section::predictive}
We review some results that aim to mitigate the ``myopic" nature of CBF-based controllers, i.e., the fact that the control input is optimal only pointwise and does not consider the system trajectories over a finite horizon ahead. The notion of a ``future-focused CBF" is introduced in~\cite{black2023future} as a solution to the unsignaled intersection-crossing problem for mixed traffic (communicating and non-communicating vehicles), so that the vehicles avoid collisions that are predicted to occur in the future. In the interest of space, we omit to present this work in detail, and refer the interested reader to~\cite{black2023future}. The following section details a predictive approach related to~\cite{black2023future} that is applicable to more general systems.

\subsection{Bird's Eye CBFs}
In this section, we consider the problems of A) designing CBFs for systems for which it is difficult to find a CBF using existing methods, and B) designing CBF-based controllers to act more proactively to maintain safety. To this end, we propose a special form of CBF that we call a ``bird's eye CBF'' (BECBF), introduced in~\cite{breeden2022predictive}. We previously called this a predictive CBF, but we note that in the broader CBF literature, the term ``predictive CBF'' is usually synonymous with ``backup CBF'', whereas the following work is distinct from any backup-type formulation, e.g.~\cite{gurriet2018towards,squires2018backup,wabersich2022predictive,chen2021backup,chen2021backup_compare,Breeden2021InputConstraints,wiltz2023construction}.


This form of CBF was specifically developed with the intent of controlling satellites in Low Earth Orbits. In this environment, a small control input applied early can have a large effect on the system trajectory over time. However, if the control input is not applied until two satellites are near collision, then the satellites must apply a very large control input to alter their trajectories in time to avoid collision. This sort of ``last-second'' behavior is undesirable and wastes fuel. Moreover, in this environment, obstacles may be moving extremely fast, so it is important to incorporate the future positions of the obstacles into the CBF formulation. Thus, we sought a control law that could maintain safety proactively using predictions about the future. The resultant CBF, while inspired by satellite orbits, is extremely general, and has proven especially useful in collision avoidance settings, such as the cars at an intersection also simulated in this section (see also~\cite{black2023future} for much more extensive simulations for this particular application). We call this tool a bird's eye CBF by analogy to having a bird's eye view of the environment and thus being able to 1) see far away obstacles entering the environment, and 2) choose avoidance maneuvers that take into account the complete size/shape of the obstacle instead of just the location of its boundary.

At this point, we emphasize that this subsection is still ongoing work. While the BECBF design herein works well for the following simulations, we note that achieving the regularity conditions specified below is a nontrivial challenge that we are still studying. 

In this section, we consider the model
\begin{equation}
    \dot{x} = f(t,x) + g(t,x) u
\end{equation}
with time $t\in\mathcal{T}\subseteq\reals$, state $x\in\mathcal{X}\subseteq\reals^n$, and control $u\in\mathcal{U} = \reals^m$. Suppose we are given a constraint function $h:\mathcal{T}\times\mathcal{X}\rightarrow\reals$ and safe set
\begin{equation}
    \mathcal{S}(t) = \{x\in\mathcal{X} \mid h(t,x) \leq 0 \label{eq:pred_safe_set}
\end{equation}
Similar to Section~\ref{section::HOCBF}, our goal is to find a CBF $H:\mathcal{T}\times\mathcal{X}\rightarrow\reals$ and a CBF set $\mathcal{S}_H = \{x\in\mathcal{X} \mid H(t,x) \leq 0\}$. The CBF that we design will be based on finite time predictions, and thus might be non-differentiable when the ``time-of-interest'' in this horizon switches from an endpoint of this horizon to the interior of this horizon. Thus, we relax the definition of CBF to absolutely continuous functions.
\begin{definition}[\textbf{CBF} (Absolutely Continuous)] \label{def:cbf_ac}
An absolutely continuous function ${\varphi} : \mathcal{T} \times\mathcal{X} \rightarrow \reals$, denoted ${\varphi}(t,x)$, is a \emph{control barrier function (CBF)} if there exists $\alpha \in \mathcal{K}$ (not necessarily locally Lipschitz continuous) such that
\begin{equation}
    \inf_{u\in\reals^m} \big[ \underbrace{\partial_t {\hspace{1pt}\varphi}( t,x ) + L_{f(t,x) + g(t,x) u} {\hspace{1pt}\varphi}( t,x )}_{=\frac{d}{dt}[{\varphi}(t,x)]} \big] \leq \alpha(-{\varphi}(t,x)) \,, \label{eq:cbf_definition}
\end{equation}
for almost every $x \in \mathcal{S}_\varphi(t), t \in \mathcal{T}$, where ${\mathcal{S}_\varphi(t)\triangleq\{x\in\mathcal{X}\mid{\varphi}(t,x)\leq 0 \} \,.}$
\end{definition}

\subsubsection{CBF Construction}

To construct the BECBF, suppose that we are given a nominal control input $\mu:\mathcal{T}\times\mathcal{X}\rightarrow\mathcal{U}$. This can be any control law, and does not need to be ``safety-encouraging'' as is required for the backup CBF. We assign the hypothetical flow of the system according to this control input to the function $p:\mathcal{T}\times\mathcal{T}\times\mathcal{X}\rightarrow\mathcal{X}$ satisfying
\begin{align}
    \frac{\partial p(\tau, t, x)}{\partial \tau} &= f(\tau, p(\tau, t, x)) + g(\tau, p(\tau, t, x)) \mu(\tau, p(\tau, t, x)) \,, \nonumber \\
    p(t,t,x) &= x \,. \label{eq:def_p}
\end{align}
We call $p$ the \emph{path function}; $p$ encodes the predicted future trajectories of the system from any initial state. Assume that $p$ is continuously differentiable. The idea of the BECBF is to use the path function to analyze whether the nominal control law $\mu$ is safe in the future, and if not, then to use the sensitivities (i.e. derivatives) of $p$ to choose a control input $u\neq \mu$ that is safe, i.e. a control input that leads to forward invariance of some subset $\mathcal{S}_H \subset \mathcal{S}$.

\begin{figure}
    \centering
    \includegraphics[width=3.5in,trim={0in, 0in, 0in, 0.31in},clip]{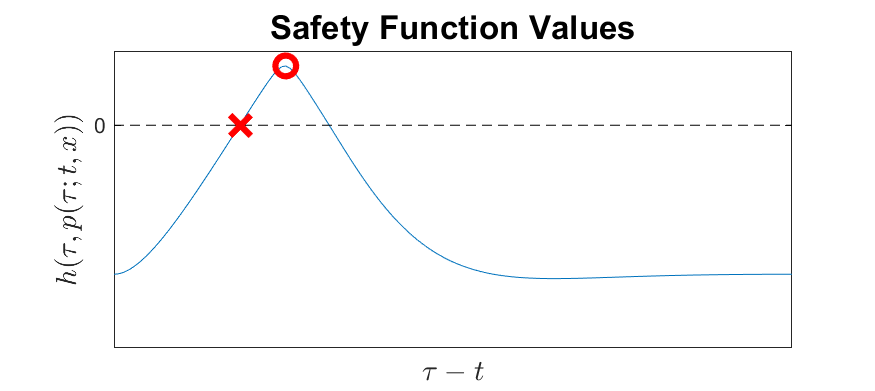}
    \caption{A possible trajectory of $\phi(\tau) = h(\tau,p(\tau,t,x))$ for a specific $(t,x)$. The trajectory has one maximizer $\boldsymbol{M}$ (red circle), and two roots. $\boldsymbol{R}$ is defined as the first root (red ``x'').}
    \label{fig:pointsofinterest}
\end{figure}

To perform this analysis, for some $(t,x)$, consider the curve $\phi(\tau) = h(\tau, p(\tau,t,x))$ containing the future values of $h$ along the path function for $\tau\in[t,t+T]$ for some fixed $T\in\reals_{>0}$. We are interested in two ``points-of-interest'' along $\phi$: 1) the maximum value of $\phi$ over $[t, t+T]$, and 2) the time at which $\phi$ first exceeds zero (i.e. the first root). For simplicity, assume that $\phi$ has a unique maximizer and no more than two roots, as shown in Fig.~\ref{fig:pointsofinterest}. Then define the quantities
\begin{align}
    \boldsymbol{K}(t,x) & = \max_{\tau\in[t,t+T]} h(\tau, p(\tau,t,x)) \,, \label{eq:def_K0} \\
    \boldsymbol{M}(t,x) & = \argmax_{\tau\in[t,t+T]} h(\tau, p(\tau,t,x)) \,, \label{eq:def_M0} \\
    \boldsymbol{R}(t,x) & = \begin{cases} \displaystyle\underset{\{\tau\in[t,t+T] \mid h(\tau, p(\tau,t,x)) = 0\}}{\min\hspace{1.1in}} \tau & \boldsymbol{K}(t,x) \geq 0 \\ \boldsymbol{M}(t,x) & \boldsymbol{K}(t,x) < 0 \end{cases} \,. \label{eq:def_R0}
\end{align} 
That is, $\boldsymbol{K}$ is the maximum of $\phi$ over $[t,t+T]$, and $\boldsymbol{M}$ is the time at which $\boldsymbol{K}$ occurs (which was assumed to be unique). If $\boldsymbol{K} \geq 0$, then $\boldsymbol{R}$ is the first root of the curve $\phi$, and otherwise, $\boldsymbol{R}$ is equal to $\boldsymbol{M}$. We note that the original paper~\cite{breeden2022predictive} considered a wider set of possible curves $\phi$, including curves with several maximizers and/or maximizer intervals. However, we do not cover those cases here 1) to encourage simplicity, and 2) because these more general trajectories are likely to conflict with the assumed regularity conditions of Theorem~\ref{thm:is_a_cbf}.

Using~\eqref{eq:def_K0}-\eqref{eq:def_R0}, we then define the BECBF as the sum of the maximal (i.e. worst point) safety metric along the predicted trajectory and a relaxation function $\beta$ of the time until the safety metric first becomes positive (i.e. unsafe), if ever. Let $\beta:\reals_{\geq 0}\rightarrow\reals_{\geq 0}$ be a class-$\mathcal{K}$ function, and then choose
\begin{equation}
    H(t,x) = \boldsymbol{K}(t,x) - \beta(\boldsymbol{R}(t,x) - t) \,. \label{eq:jcbf0}
\end{equation}
Intuitively,~\eqref{eq:jcbf0} says that a state $(t,x)$ belongs to $\mathcal{S}_H$ if the time at which the hypothetical trajectory $p$ first becomes unsafe is sufficiently far in the future, as measured by $\beta$, that we can adjust the value of $\boldsymbol{K}$ before the system reaches any unsafe states. As such, the selection of $\beta$ will substantially impact the size of $\mathcal{S}_H$ and the amount of control effort utilized to correct the trajectory, as well as how early this control effort is applied.

The main results of~\cite{breeden2022predictive} are then as follows.

\begin{lemma}[\textbf{Zero Sublevel Set of $H$}]
\label{thm:subset}
For some constraint function $h:\mathcal{T}\times\mathcal{X}\rightarrow\reals$, let $\mathcal{S}$ be as in~\eqref{eq:safe_set} and $H$ as in~\eqref{eq:jcbf0} with $\mathcal{S}_H \triangleq \{ x \in \mathcal{X} \mid H(t,x) \leq 0\}$. Then $\mathcal{S}_H(t)\subseteq\mathcal{S}(t), \forall t \in \mathcal{T}$.
\end{lemma}

\begin{theorem}[\textbf{$H$ is a CBF}]
\label{thm:is_a_cbf}
Let the derivative $\beta':\reals_{\geq0}\rightarrow\reals_{\geq0}$ of $\beta$ in~\eqref{eq:jcbf0} be strictly positive on $(0,T)$. Assume that $H$ in~\eqref{eq:jcbf0} is absolutely continuous, and further assume that $\boldsymbol{M}$ in~\eqref{eq:def_M0} is continuously differentiable whenever $\boldsymbol{M}(t,x) \in (t,t+T)$. Let $\mathcal{U} = \reals^m$. Suppose that there exists $\gamma$ such that $\frac{d}{d\tau}[\phi(\tau)] \equiv \frac{d}{d\tau}[h(\tau,p(\tau;t,x))] \leq \gamma$ for all $\tau\in\mathcal{T},t\in\mathcal{T},x\in\mathcal{X}$. Suppose also that $\| \frac{\partial h(\eta, p(\eta; t,x))}{\partial x}\frac{\partial p(\eta; t,x)}{\partial x} g(t,x) \| \neq 0$ for all $\eta \in (t,t+T)\setminus \boldsymbol{M}(t,x)$ and for all $t\in\mathcal{T},x\in\mathcal{X}$, and that $\frac{\partial h(\boldsymbol{M}, p(\boldsymbol{M};t,x))}{\partial x}\transpose \frac{\partial h(\boldsymbol{R}, p(\boldsymbol{R};t,x))}{\partial x} \geq 0$ whenever $\boldsymbol{R}(t,x)\neq\boldsymbol{M}(t,x)$ and for all $t\in\mathcal{T},x\in\mathcal{X}$. Then $H$ in~\eqref{eq:jcbf0} is a CBF as in Definition~\ref{def:cbf_ac}.
\end{theorem}

That is, under mild assumptions, the function $H$ in~\eqref{eq:jcbf0} is indeed a CBF, and its zero sublevel set $\mathcal{S}_H$ is a subset of the constraint set $\mathcal{S}$. In brief, the assumptions of Theorem~\ref{thm:is_a_cbf} are 1) the slope of $\phi$ is bounded, 2) the trajectories encoded in $p$ have non-zero sensitivity to the initial state (i.e. the system is controllable), and 3) the sensitivities of $p$ satisfy a consistency condition so that decreasing $\boldsymbol{K}$ does not cause $\boldsymbol{R}$ to occur sooner. See~\cite{breeden2022predictive} for further discussion. We also assume that $\mathcal{U}$ is equivalent to $\reals^m$, though in practice, the function $\beta$ can be tuned to achieve input constraint satisfaction. This is a powerful theorem because it implies the form~\eqref{eq:jcbf0} can be applied to any system and any constraint function of any relative degree. It is also advantageous in practice, as will be clear in the simulations.


Next, we provide two lemmas on how to compute the derivatives of the BECBF~\eqref{eq:jcbf0}. In these lemmas, the notation $\frac{\partial p(a,b,c)}{\partial \lambda_x}$ means $\frac{\partial p(\tau,t,x)}{\partial x} \big|_{(\tau,t,x) = (a,b,c)}$.

\begin{lemma}[\textbf{Derivative of $\boldsymbol{K}$}] \label{lemma:K_der}
    For some $(t,x)\in\mathcal{T}\times\mathcal{X}$, suppose that $h$ in~\eqref{eq:pred_safe_set} is continuously differentiable in a neighborhood of $\boldsymbol{M}(t,x)$ and $p^* = p(\boldsymbol{M}(t,x),t,x)$. Assume that $\boldsymbol{M}$ in~\eqref{eq:def_M0} is Lipschitz continuous. Let $\phi'$ be the derivative $\frac{d}{d\tau}[\phi(\tau)]$ of $\phi(\tau) = h(\tau,p(\tau,t,x))$. Then
    \begin{multline}
        \frac{d \boldsymbol{K}(t,x)}{dt} = \frac{\partial h(\boldsymbol{M}, p^*)}{\partial\lambda_x} \frac{\partial p(\boldsymbol{M}, t, x)}{\partial \lambda_x} g(t,x) (u - \mu(t,x))  \\ + 
        \begin{cases} \frac{\partial h(\boldsymbol{M},p^*)}{\partial \lambda_t} + \frac{\partial h(\boldsymbol{M},p^*)}{\partial \lambda_x} \phi'(\boldsymbol{M}) & \boldsymbol{M} \in \{t, t+T\} \\ 0 & \boldsymbol{M} \in (t,t+T)\end{cases} 
        \label{eq:K_der}
    \end{multline}
\end{lemma}

\begin{lemma}[\textbf{Derivative of $\boldsymbol{R}$}] \label{lemma:R_der}
    For some $(t,x)\in\mathcal{T}\times\mathcal{X}$, suppose that $\boldsymbol{R}(t,x) \neq \boldsymbol{M}(t,x)$ and that $h$ is continuously differentiable in a neighborhood of $\boldsymbol{R}(t,x)$ and $p^\circ = p(\boldsymbol{R}(t,x),t,x)$. Let $\phi'$ be the derivative $\frac{d}{d\tau}[\phi(\tau)]$ of $\phi(\tau) = h(\tau,p(\tau,t,x))$. Then
    \begin{equation}
        \frac{d\boldsymbol{R}(t,x)}{dt} = \frac{1}{-\phi'(\boldsymbol{R})}\frac{\partial h(\boldsymbol{R},p^\circ)}{\partial \lambda_x} \frac{\partial p(\boldsymbol{R},t,x)}{\partial \lambda_x} g(t,x)(u - \mu(t,x))\label{eq:R_der}
    \end{equation}
\end{lemma}

That is, while we could compute the derivatives of $H$ in~\eqref{eq:jcbf0} numerically (which would require recomputing $\boldsymbol{K}$, $\boldsymbol{M}$, and $\boldsymbol{R}$ at least $n$ times), we also have explicit expressions for the derivatives of $H$ in terms of the derivatives of $h$ and $p$. Moreover, these derivatives are zero if $u = \mu$ and $\boldsymbol{M}$ is not an endpoint. Unfortunately, the expression for $\frac{d\boldsymbol{R}}{dt}$ in the case of $\boldsymbol{R} = \boldsymbol{M}$ is prohibitively complex, so this should be computed numerically. Note that if $\boldsymbol{M}(t,x)$ is continuously differentiable at $x$, then
\begin{equation}
    \frac{d\boldsymbol{M}(t,x)}{dt} = \frac{\partial \boldsymbol{M}(t,x)}{\partial \lambda_x} g(t,x) (u - \mu(t,x)) \label{eq:M_der}
\end{equation}
has the same structure as~\eqref{eq:K_der}-\eqref{eq:R_der}, where $\frac{\partial \boldsymbol{M}}{\partial \lambda_x}$ is instead computed numerically. One can also replace $\boldsymbol{R}-t$ in~\eqref{eq:jcbf0} with the current safety metric $h$, as was done in~\cite[Eq.~22]{black2023future}. This results in simpler derivative expressions than~\eqref{eq:K_der}-\eqref{eq:R_der}, but is less applicable to the satellites scenario for which this was intended.

\subsubsection{Simulations}

To demonstrate the advantages of the BECBF over conventional CBFs, we consider two case studies. Our first case study involves two cars passing through a four-way intersection. We assume that the cars are fixed in their lanes $l_1:\reals\rightarrow\reals^2$ and $l_2:\reals\rightarrow\reals^2$, respectively, with locations $z_1\in\reals$ and $z_2\in\reals$ along their lanes. Thus, the position of car~1 on the road is $l_1(z_1)\in\reals^2$ and the position of car~2 is $l_2(z_2)\in\reals^2$. For simplicity, we model the cars as double-integrators: $\ddot{z}_1 = u_1$ and $\ddot{z}_2 = u_2$, resulting in state vector $x = [z_1;\, \dot{z}_1;\, z_2;\, \dot{z}_2]$ and control vector $u = [u_1; u_2]$. Suppose the cars nominally want to travel in their lanes at velocities $v_1$ and $v_2$, respectively, so the nominal control input is $\mu = [\mu_1; \mu_2]$ where $\mu_i(t,x) = k(v_i - \dot{z}_i)$ for some gain $k > 0$. The path function is then $p = [p_1;\,p_2]$ where $p_i$ is the solution to~\eqref{eq:def_p} under $\mu$. The function $p$ is computed explicitly here, though we note that numerical solutions to~\eqref{eq:def_p} are also fine.
\begin{equation}
    p_i(\tau;t,x) = \begin{bmatrix} z_i + v_i(\tau - t) + \frac{\dot{z}_i - v_i}{k}(1 - e^{-k(\tau - t)}) \\ v_i + (\dot{z}_i - v_i)e^{-k(\tau-t)} \end{bmatrix}.
\end{equation}
Let the safety constraint be $h = \rho - \|l_1(z_1) - l_2(z_2)\|$.

This system meets all the assumptions of Theorem~\ref{thm:is_a_cbf}. The system also meets the assumptions of the formulas~\eqref{eq:K_der}-\eqref{eq:M_der}, except on the critical manifolds $l_1(p_1^*) = l_1(p_2^*)$ (where $h$ is not continuously differentiable) and $\boldsymbol{K}(t,x) = 0$ (where~\eqref{eq:def_R0} switches cases). These manifolds are of Lebesgue measure zero and can be ignored, so it is straightforward to apply the CBF $H$ in~\eqref{eq:jcbf0} in a QP
\begin{equation}
    u = \argmin_{u\in\reals^2} \| u - \mu(t,x) \|^2 \textrm{  s.t. } \frac{dH(t,x)}{dt} \leq \alpha(-H(x)) \,. \label{eq:qp_pred}
\end{equation}
Note that this is a centralized control law, since $u_1$ and $u_2$ are computed simultaneously.

\begin{figure}
    \centering
    \includegraphics[width=3.5in]{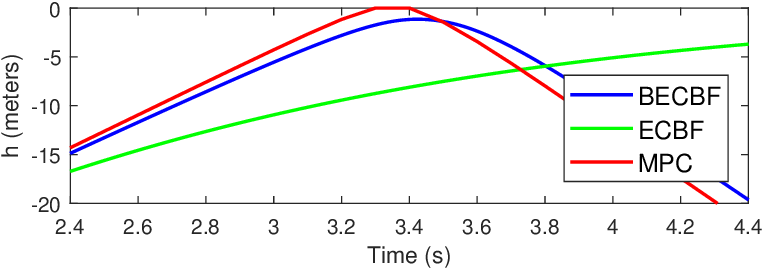}
    \caption{The values of $h$ during the three simulations. The BECBF and MPC trajectories are similar, whereas the ECBF trajectory slowly converges to zero as the vehicles come to a complete stop.}
    \label{fig:safety_pred}
\end{figure}

We then simulated two cars approaching an intersection with the control law~\eqref{eq:qp_pred}, with the same control law with an ECBF~\cite{nguyen2016exponential} in place of the BECBF, and with a nonlinear MPC control law. The results are best demonstrated by the video \url{https://youtu.be/0tVUAX6MCno}, and the safety function values are shown in Fig.~\ref{fig:safety_pred}. The BECBF and MPC cases performed generally similarly, with both cars approaching close and then continuing on opposite sides of the intersection. On the other hand, the ECBF caused both agents to come to a complete stop, so neither agent made it through the intersection. On average, the control computation time was 0.0011~s for the ECBF, 0.0061~s for the BECBF, and 0.40~s for the MPC approach, all in MATLAB on a 3.5 GHz CPU, though the code for all three cases could likely be further optimized. For more information and simulation code, see~\cite{breeden2022predictive}.

\begin{figure}
    \centering
    \includegraphics[width=3.5in]{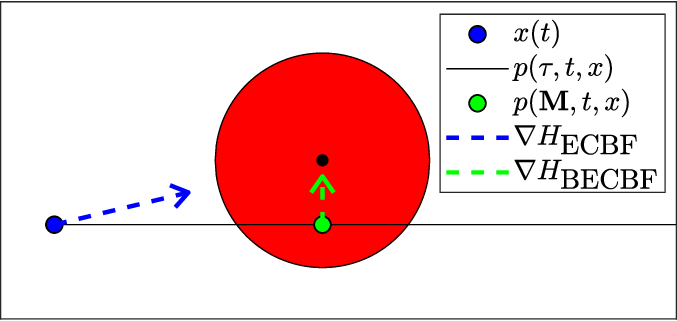}
    \caption{This shows a possible predicted trajectory through a static obstacle (red zone). A conventional CBF like the ECBF has a gradient vector pointing from the current agent position to the center of the obstacle as in the blue vector. By contrast, the BECBF gradient is the green vector, which takes into account the entire obstacle instead of just the closest point. This is why we call~\eqref{eq:jcbf0} a ``bird's eye CBF''.}
    \label{fig:2d_static}
\end{figure}

Thus, the BECBF achieved similar performance to the nonlinear MPC with substantially reduced computation time. This was possible because the BECBF only performs one prediction of $\phi$ per control cycle and only varies the current control input, whereas MPC varies all the states and control inputs within a horizon. Compared to the ECBF, the QP with BECBF chooses the control input that most encourages safety at the moment where the two cars are closest together---in this case, that means one car decelerating and one car accelerating so that the cars are never too close together. By contrast, the QP with ECBF always chooses the control that most encourages safety at the present---in this case, that means applying a control input that is opposite the present direction of the other car (and constrained to be along the lanes $l_1$ and $l_2$), which is a deceleration for both cars. The differences between these directions for a static obstacle is illustrated in Fig.~\ref{fig:2d_static}.

Next, we consider the case study of satellites in low earth orbit with state $x = [r\transpose\;\dot{r}\transpose]\transpose\in\reals^6$ and dynamics $\ddot{r} = -\frac{\mu r}{\|r\|^3} + u$. In this example, $\mu = 3.986(10)^{14}$ is very large (much larger than the example in Section~\ref{section::HOCBF}), so the uncontrolled term of $\ddot{r}$ will always be much larger than $u$ and the system state will evolve rapidly. This means that 1) conventional approaches like the ECBF or Theorem~\ref{thm:variable_cbf} will not work very well, and 2) a small control input will have a large effect on system trajectories over time.

\begin{figure}
    \centering
    \includegraphics[width=3.5in]{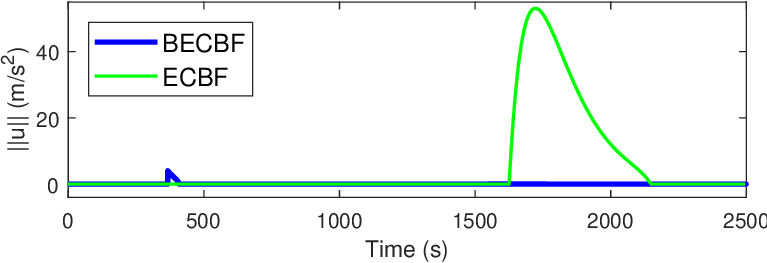}
    \caption{Control thrusts of a satellite with safety determined by an ECBF or BECBF.}
    \label{fig:control_space}
\vspace{10pt}
    \centering
    \begin{tikzpicture}
        \node[anchor=south west,inner sep=0] (image) at (0,0) {\includegraphics[width=3.5in]{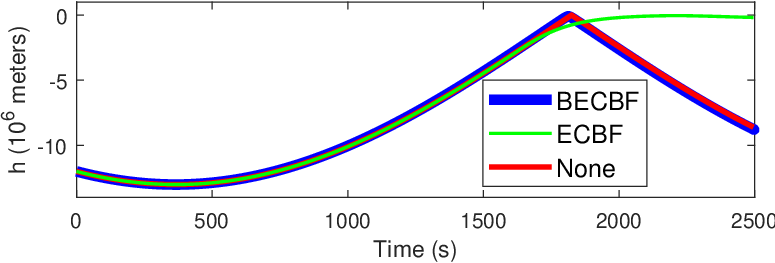}};
        \node[anchor=south west,inner sep=0] (image) at (1.25,1.15) 
        {\includegraphics[width=0.75in,trim={0in, .46in, 0in, .14in},frame,clip]{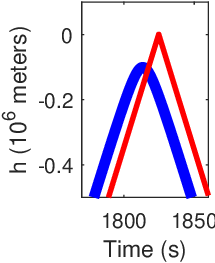}}; 
    \end{tikzpicture}
    \caption{The values of $h$ during the two simulations in Fig.~\ref{fig:control_space} and a simulation with zero control input. The BECBF minimally modifies the original trajectory so the red and blue lines are very similar, but the red line exceeds the safe set by 1000 m (see zoomed in inset or the animation), while the blue and green lines remain safe.}
    \label{fig:safety_space}
\end{figure}

Let $h(t,x) = \rho - \| r - r_o(t) \|$ where $\rho = 1000$ and $r_o$ is the location of an uncontrolled piece of debris. Let $\mu = 0$, so the satellite nominally applies no control input, unless $u\neq 0$ is necessary to ensure safety. Our simulation scenario places the controlled satellite and the debris initially very far apart, but in orbits that eventually intersect if no control action is taken. Simulations under the BECBF, under an ECBF, and under no control action are shown in Figs.~\ref{fig:control_space}-\ref{fig:safety_space} and in the video  \url{https://youtu.be/HhtWUG63BWY}. Note how the BECBF trajectory (blue) takes a small control action as soon as the unsafe prediction enters the horizon at $t=367\textrm{ s}$, and then is very similar to the nominal trajectory (red). On the other hand, the ECBF trajectory (green) takes control action much later, when over 10 times as much thrust is required. This avoidance problem could in theory also be solved with MPC, but would require a very fine discretization, because $h > 0$ occurs for only 0.14 seconds since the satellites are moving so fast.  Thus, utilizing the same length of prediction horizon would require more than $10^4$ samples, making the MPC problem intractable.

Thus, in addition to choosing a better control direction than the ECBF, the BECBF also provides a mechanism for tuning how early we want the system to detect and respond to predicted collisions. The BECBF is also able to take into account the future locations of time-varying obstacles with known paths rather than just their current positions. Thus, the BECBF can make use of paths like satellite orbits, and in the cars example, the BECBF can take into account whether each car intends to continue straight or turn. Note that both of these examples included fairly simple path functions, where the controlled agents were always in motion so that $\phi$ had a strictly negative hessian, and thus $\boldsymbol{K}$, $\boldsymbol{M}$, and $\boldsymbol{R}$ were always well-defined and differentiable. In the future, we seek to consider more general paths that may challenge the regularity assumptions presently made.

\section{Practical Challenges: Output Feedback Control and Sampled-Data Control with Control Barrier Functions}
\label{section::implementation}

\subsection{Output Feedback Control}
Synthesizing safe controllers for nonlinear systems using output feedback can be a challenging task, since observers and controllers designed independently of each other may not render the system safe. In our recent work~\cite{agrawal2022safe} we present two observer-controller interconnections that ensure that the trajectories of a nonlinear system remain safe despite bounded disturbances on the system dynamics and partial state information. The first approach utilizes Input-to-State Stable observers, and the second uses Bounded Error observers. Using the stability and boundedness properties of the observation error, we construct novel Control Barrier Functions that impose inequality constraints on the control inputs which, when satisfied, certify safety. We propose quadratic program-based controllers to satisfy these constraints, and prove Lipschitz continuity of the derived controllers.

\subsubsection{Tunable-Robust CBFs}
We consider nonlinear control-affine systems of the form:
\begin{subequations}
\label{eqn:sys}
\begin{align}
\dot x &= f(x) + g(x) u + g_d(x) d(t), \label{eqn:main_dyn}\\
y &= c(x) + c_d(x) v(t), \label{eqn:measurements}
\end{align}
\end{subequations}
where $x \in \mathcal X \subset \mathbb R^n$ is the system state, $u \in \mathcal U \subset \mathbb R^m$ is the control input, $y \in \mathbb R^{n_y}$ is the measured output, $d: \mathbb R^+ \to \mathbb R^{n_d}$ is a disturbance on the system dynamics, and $v: \mathbb R^+ \to \mathbb R^{n_v}$ is the measurement disturbance. We assume $d$ and $v$ are piecewise continuous, bounded disturbances, $\sup_{t}\|d(t)\|_\infty = \bar d$ and $\|v(t)\|_\infty \leq \bar v$ for some known $\bar d, \bar v < \infty$. The functions $f : \mathcal X \to \mathbb R^n$, $g: \mathcal X \to \mathbb R^{n \times m}$, $c : \mathcal X \to \mathbb R^{n_y}$, $g_d : \mathcal X \to \mathbb R^{n \times n_d}$, and $c_d : \mathcal X \to \mathbb R^{{n_y} \times n_v}$ are all assumed to be locally Lipschitz continuous. Notice that $g_d(x) d(t)$ accounts for either matched or unmatched disturbances.

We seek to establish observer-controller interconnections of the form:
\begin{subequations}
\label{eqn:observer-controller}
\begin{align}
\dot {\hat x} &= p(\hat x, y) + q(\hat x, y) u, \label{eqn:observer}\\
u &= \pi(t, \hat x, y), \label{eqn:controller}
\end{align}
\end{subequations}
where $p : \mathcal X \times \mathbb R^{n_y} \to  \mathbb R^n$, $q: \mathcal X \times \mathbb R^{n_y} \to \mathbb R^{n \times m}$ are locally Lipschitz in both arguments. The feedback controller $\pi : \mathbb R^+ \times \mathcal X \times \mathbb R^p \to \mathcal U$ is assumed piecewise-continuous in $t$ and Lipschitz continuous in the other two arguments. Then, the closed-loop system formed by~(\ref{eqn:sys}, \ref{eqn:observer-controller}) is
\begin{subequations}
\label{eqn:closed_loop_sys}
\begin{align}
\dot x &= f(x) + g(x)u + g_d(x) d(t),\\
\dot {\hat x} &= p(\hat x, y) + q(\hat x, y) u,\\
x(0) &= x_0, \ \hat x(0) = \hat x_0,
\end{align}
\end{subequations}
where $y$ and $u$ are defined in~\eqref{eqn:measurements} and~\eqref{eqn:controller} respectively. Under the stated assumptions, there exists an interval $\mathcal I = \mathcal I(x_0, \hat x_0) = [0, t_{max}(x_0, \hat x_0))$ over which solutions to the closed-loop system exist and are unique~\cite[Thm 3.1]{khalil2002nonlinear}.

Safety is defined as the true state of the system remaining within a \emph{safe set}, $\mathcal S \subset \mathcal X$, for all times $t \in \mathcal I$, where the safe set $\mathcal S$ is defined as the super-level set of a continuously-differentiable function $h: \mathcal X \to \mathbb R$ as in~\eqref{eq::safeset}.

A state-feedback controller\footnote{In \emph{state-feedback} the control input is determined from the true state, $u = \pi(t, x)$. In \emph{estimate-feedback} the input is determined from the state estimate and measurements, $u = \pi(t, \hat x, y)$.} $\pi : \mathcal R^+ \times \mathcal X \to \mathcal U$ \emph{renders system~\eqref{eqn:sys} safe} with respect to the set $\mathcal S$, if for the closed-loop dynamics $\dot x = f(x) + g(x) \pi(t, x) + g_d(x) d(t)$, the set $\mathcal S$ is \emph{forward invariant}, i.e., $x(0) \in \mathcal S \implies x(t) \in \mathcal S \ \forall t \in \mathcal I$. In output-feedback we define safety as follows:
\begin{definition}
An observer-controller pair~\eqref{eqn:observer-controller} \emph{renders system~\eqref{eqn:sys} safe} with respect to a set $\mathcal S \subset \mathcal X$ from the initial condition sets $\mathcal X_0, \hat{\mathcal X}_0 \subset \mathcal S$ if for the closed-loop system~\eqref{eqn:closed_loop_sys} it holds that $x(0) \in \mathcal X_0 \text{ and } \hat x(0) \in \hat{\mathcal X}_0 \implies x(t) \in \mathcal S \quad \forall t \in \mathcal I.$
\end{definition}

Now, inspired by~\cite{jankovic2018robust} and~\cite{alan2021safe}, we define the following CBF to account for disturbances and measurement noise.
\begin{definition}
A continuously differentiable function $h : \mathcal X \to \mathbb R$ is a \emph{Tunable Robust CBF} (TR-CBF) for  system~\eqref{eqn:sys} if there exists a class $\mathcal K$ function $\alpha$, and a continuous, non-increasing function $\kappa: \mathbb R^+ \to \mathbb R$ with $\kappa(0) = 1$, s.t.
\begin{align}
\sup_{u \in \mathcal U} & \ L_fh(x) +  L_gh(x) u + \alpha(h(x)) \notag \\ 
&  \ \geq \kappa(h(x)) \|{L_{g_d}h(x)}\| \bar d,  \ \forall x \in \mathcal S.
\end{align} 
\end{definition}
Examples include $\kappa(r) = 1$ and $\kappa(r) = 2/(1 + \exp(r))$.
Given a TR-CBF $h$ for~\eqref{eqn:sys}, the set of safe control inputs is
\begin{align}
\label{eqn:Ktrcbf}
K_{trcbf}(x) = \{ &u \in \mathcal U :  \ L_fh(x) + L_gh(x) u -\notag\\ &\quad \kappa(h(x)) \|{L_{g_d}h(x)}\| \bar d \geq -\alpha(h(x))\},
\end{align}
and a safe state-feedback controller is obtained by solving a QP, as in~\cite[Eq. 30]{jankovic2018robust}. The main question is:
\begin{problem}
Given a system~\eqref{eqn:sys} with disturbances of known bounds $\|{d(t)}\|_\infty \leq \bar d$, $\|{v(t)}\|_\infty \leq \bar v$, and a safe set $\mathcal S$ defined by~\eqref{eq::safeset}, synthesize an interconnected observer-controller~\eqref{eqn:observer-controller} and the initial condition sets $\mathcal X_0, \hat{\mathcal X}_0$ to render the system safe.
\end{problem}


\subsubsection{Observer-Controller Interconnection}

We review Approach 2 of our work in~\cite{agrawal2022lcss}, where we consider the class of Bounded-Error Observers: 
\begin{definition}
An observer~\eqref{eqn:observer} is a \emph{Bounded-Error (BE) Observer}, if there exists a bounded set $\mathcal D(\hat x_0) \subset \mathcal X$ and a (potentially) time-varying bounded set $\mathcal P(t, \hat x) \subset \mathcal X$ s.t. $x_0 \in \mathcal D(\hat x_0) \implies x(t) \in  \mathcal P(t, \hat x) \ \forall t \in \mathcal I$.
\end{definition}

The idea is to find a common, safe input for all $x \in \mathcal P(t, \hat x)$:
\begin{theorem}
\label{thm:intersections}
For system~\eqref{eqn:sys}, suppose the observer~\eqref{eqn:observer} is a Bounded-Error observer. Suppose the safe set $\mathcal S$ is defined by a continuously differentiable function $h : \mathcal X \to \mathbb R$, where $h$ is a Tunable Robust-CBF for the system. Suppose $\pi: \mathbb R^+ \times \mathcal X \to \mathcal U$ is an estimate-feedback controller, piecewise-continuous in the first argument and Lipschitz continuous in the second, s.t. 
\begin{align}
\label{eqn:intersections_cond}
\pi(t, \hat x) \in  \bigcap_{x \in \mathcal P(t, \hat x)} K_{trcbf}(x),
\end{align} 
where $K_{trcbf}$ is defined in~\eqref{eqn:Ktrcbf}. Then the observer-controller renders the system safe from the initial-condition sets $x(0) \in \mathcal{X}_0 = \mathcal D(\hat x_0)$ and $\hat x_0 \in \hat{\mathcal{X}}_0 = \{ \hat x: \mathcal P(0, \hat x_0) \subset \mathcal S \}$.
\end{theorem}

\begin{figure}[tbp]
   \centering
   \includegraphics[width=0.8\columnwidth]{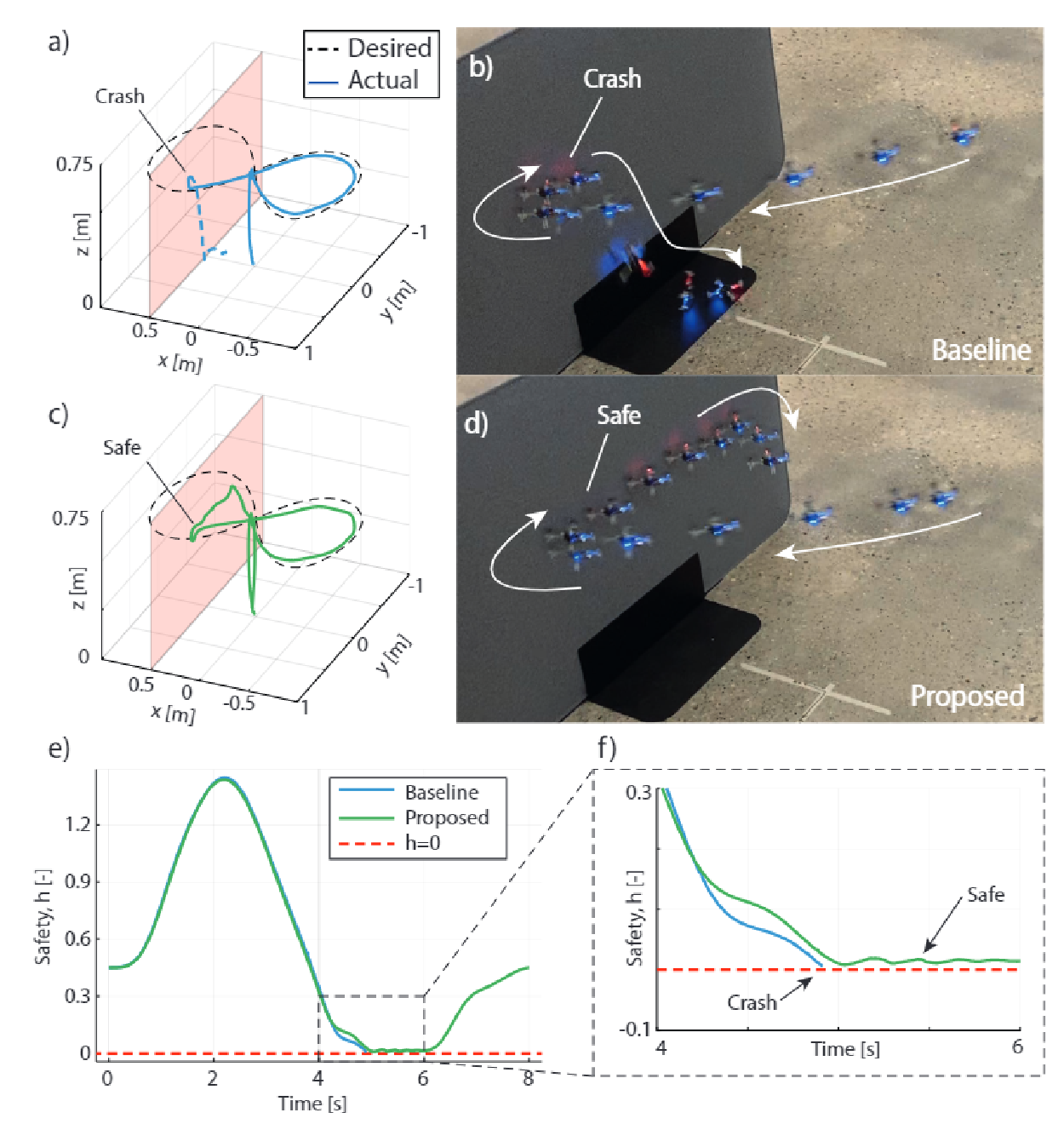} 
   \caption{Experimental results. The quadrotor is commanded to track a figure-of-eight trajectory, while avoiding the physical barrier at $x=0.5$~m. Ground truth trajectories are plotted in (a, c) for the baseline CBF and proposed controllers respectively. Snapshots from the experiment are show in~(b,~d). (e, f) Plots of the safety value, $h$ over time for both trajectories. Video: \url{https://www.youtube.com/watch?v=83o7CEgTEwo}}
   \label{fig:fig1}
\end{figure}

In general, designing a controller satisfying~\eqref{eqn:intersections_cond} can be difficult. Under certain assumptions, one can define certain forms of a Quadratic Program that defines a controller that meets the desired properties. In the interest of space, we refer the reader to~\cite{agrawal2022lcss} for the details, and we directly present some experimental results of the safe observer-controller interconnection. The objective for a quadrotor is to fly in a ``figure of eight" trajectory, but to not crash into a physical barrier placed at $x=0.5$~meters. An Extended Kalman Filter is used as the bounded error observer. To design the controller, first $\pi_{des}(t, \hat x)$ is computed using an LQR controller, which computes desired accelerations wrt to an inertial frame to track the desired trajectory. This command is filtered using a safety-critical QP, either the baseline CBF-QP (Figure~\ref{fig:fig1}a) or the proposed QP using Approach 2 (Figure~\ref{fig:fig1}c). 
The trajectories from the two flight controllers are compared in~Figure~\ref{fig:fig1}. In the baseline controller, the quadrotor slows down as it approaches the barrier, but still crashes into barrier. In the proposed controller, the quadrotor remains safe,~Figure~\ref{fig:fig1}e.


\subsection{Zero-Order Hold (ZOH) Control} \label{sec:sampled}
In this section, we consider one of the major challenges that continuous-time CBF-based controllers (such as those derived in the previous sections) face in practice, namely that physical systems evolve in continuous time, under control inputs that are implemented at discrete time instances, such as zero-order-hold (ZOH) controllers with fixed time-step. One can easily construct counter-examples showing that the control laws developed from the CBF condition in~\citep{ames2017control,wang2017safety} are no longer safe when the controller is executed in discrete steps. On the other hand, a controller implemented under discrete-time CBFs~\cite{Discrete_CBFs,Blanchini} may not satisfy the continuous safety condition between time steps~\citep{CBF_Planning}. 

In our paper~\cite{breeden_acc21} we study conditions for forward invariance of safe sets under ZOH controllers. We define two types of margins, the \textit{controller margin} and the \textit{physical margin}, to compare the conservatism of the conditions developed. In~\cite{breeden_acc21}, we present extensions to the approaches in~\cite{Cortez,usevitch2021adversarial,Singletary} that reduce conservatism as measured by these margins, while similarly relying on proving that the continuous-time CBF condition is always satisfied. We also present a novel condition inspired instead by discrete-time set invariance conditions, and compare the conservatism of all the approaches studied using the above margins. For brevity, we only present the prior state-of-the-art and this last approach here, and we refer the reader to~\cite{breeden_acc21} for details about the other approaches. We also build upon the following approach further in~\cite{breeden2023aiaa}.

\subsubsection{Problem Formulation}

We consider the system
\begin{equation}
	\dot{x} = f(x) + g(x) u \,, \label{eq:model}
\end{equation}
with state $x \in \mathbb R^n$, control input $u\in U \subset \mathbb R^m$ where $U$ is compact, and locally Lipschitz continuous functions $f: \mathbb R^n \rightarrow \mathbb R^n$ and $g: \mathbb R^n \rightarrow \mathbb R^{n \times m}$. Define $u_\textrm{max} \triangleq \max_{u\in U} ||u||$. Let $h: \mathbb R^n \rightarrow \mathbb R$ where $h \in C^1_{loc}$, and define a safe set $S$ as 
\begin{equation}
	S \triangleq \{ x \in \mathbb R^n \mid h(x) \leq 0\} \,. \label{eq:safe_set}
\end{equation}

For a continuous control law $u(x)$, Theorem \ref{thm: CBF safety} (with the adjusted sign of the CBF) can be used for guaranteeing safety of dynamical systems. 
%
To apply Theorem \ref{thm: CBF safety}, we must ensure~\eqref{eq:cbf_online_condition} (with the inequality reversed) is satisfied along $x(t)$ for all $t \geq 0$. However, suppose instead that the state $x$ is only measured discretely (and thus the control policy $u(x)$ is updated in a discrete fashion too) at times $t_k = kT, k=0,1,2,\cdots$ for a fixed time-step $T\in\mathbb R_{>0}$. 
Consider a ZOH control law\footnote{Under $u$ as in~\eqref{eq:piecewise_control} for a compact set $U$, uniqueness of the maximal closed-loop solution $x(t)$ (and hence $x_k$) is guaranteed by~\cite[Thm. 54]{uniqueness}.}
\begin{equation}
    u(t) = u_k,\; \forall t \in [t_k, t_{k+1}) \,, \label{eq:piecewise_control}
\end{equation}
where $u_k = u_k(x_k)\in U$ and $x_k = x(t_k)$, $\forall k\in \mathbb N$. 
Satisfaction of~\eqref{eq:cbf_online_condition} only discretely is not sufficient for safety per Theorem \ref{thm: CBF safety}.
Thus, we seek a condition similar to~\eqref{eq:cbf_online_condition} under which safety can be guaranteed when the control input is updated only at discrete times. We consider the following problem.

\begin{problem} \label{problem:design}
    Design a function $\phi:\mathbb R_{>0}\times\mathbb R^n \rightarrow \mathbb R$ such that any bounded, piecewise-constant control input $u\in U$ of the form~\eqref{eq:piecewise_control} satisfying
    \begin{equation}
        L_f h(x_k) + L_g h(x_k) u_k \leq \phi(T, x_k), \label{eq:zoh_cbf_condition}
    \end{equation}
    at the sampled states $x_k=x(kT), k\in\mathbb N$ renders $S$ forward invariant along the closed-loop trajectories of~\eqref{eq:model}.
\end{problem}
\noindent We call~\eqref{eq:zoh_cbf_condition} the ZOH-CBF condition. 
The following result, adapted from~\cite{Cortez}, provides one form of the function $\phi$ that solves Problem \ref{problem:design} (see also~\cite{usevitch2021adversarial}).

\begin{lemma}[{\hspace{-0.2pt}\cite[Thm. 2]{Cortez}}]
	\label{prior:zoh}
	Let the set $S$ in~\eqref{eq:safe_set} be compact and $\alpha\in\mathcal{K}$ be locally Lipschitz continuous. Let $l_{L_fh}, l_{L_gh}, l_{\alpha(h)}$ be the Lipschitz constants of $L_fh, L_gh, \alpha(-h)$, respectively. Then the function $\phi_0^g:\mathbb R_{>0}\times\mathbb R^n$, defined as
	\begin{equation}
		\phi_0^g(T,x) \triangleq \alpha(-h(x)) - \frac{l_1 \Delta}{l_2}\left( e^{l_2 T} - 1\right) \,, \label{eq:phi0}
	\end{equation}
	solves Problem \ref{problem:design}, 	where $l_1 = l_{L_f h} + l_{L_g h}u_\textrm{max} + l_{\alpha(h)}, l_2 = l_{L_f h} + l_{L_g h}u_\textrm{max}$, and $\Delta = \sup_{x\in S,u\in U}||f(x) + g(x)u||$.
\end{lemma}

In practice, the form of the function $\phi_0^g$ in~\eqref{eq:phi0} is conservative in the sense that many safe trajectories may fail to satisfy~\eqref{eq:zoh_cbf_condition} for $\phi=\phi_0^g$. To overcome this limitation, we define two metrics to quantify the conservatism of the solutions to \ref{problem:design} and then develop novel solutions to Problem \ref{problem:design} that are less conservative compared to~\eqref{eq:phi0}.

\subsubsection{Comparison Metrics} 

In this work, we consider functions $\phi$ of the form:
\begin{equation}
    \phi
    (T,x)= \alpha(-h(x)) - \nu(T,x) \,, \label{eq:phi_form}
\end{equation}
where $\alpha$ is a class-$\mathcal{K}$ function that vanishes as $h(x)\rightarrow 0$, and $\nu:\mathbb R_{>0}\times\mathbb R^n\rightarrow\mathbb R$ is a function of the discretization time-step $T$ and the state $x$ that does not explicitly depend on $h$. This motivates our first metric of comparison, defined as follows.

\begin{definition}[\textbf{Controller margin}] \label{def:control_margin}
    The function $\nu$ in~\eqref{eq:phi_form} is called the \textnormal{controller margin}.
\end{definition}

Note that $\nu$ is the difference between the right-hand sides of conditions~\eqref{eq:cbf_online_condition} and~\eqref{eq:zoh_cbf_condition}, and is a bound on the discretization error that could occur between time steps. At a given state $x \in S$, a larger controller margin will necessitate a larger control input to satisfy~\eqref{eq:zoh_cbf_condition} (hence, the name ``controller margin"). A sufficiently large controller margin might also necessitate inadmissible control inputs, and thus make a CBF no longer applicable to a system. Thus, it is desired to design functions $\phi$ whose controller margins are small. For a given $T$, we call a solution $\phi_a$ less conservative than $\phi_b$ if the controller margins of $\phi_a$ and $\phi_b$ satisfy $\nu_a(T,x) \leq \nu_b(T,x), \forall x \in S$.

The controller margin is called \textit{local} (denoted as $\nu^l(T,x)$) if $\nu$ varies with $x$, and \textit{global} (denoted as $\nu^g(T)$) if $\nu$ is independent of $x$. The superscripts $l$ and $g$, respectively, denote the corresponding cases, and $\nu$ is denoted with the same sub/superscripts as the corresponding $\phi$ function. For instance, \begin{equation}
\nu_0^g(T)=\frac{l_1\Delta}{l_2}(e^{l_2T}-1) \label{eq:nu0}
\end{equation} is the controller margin of $\phi_0^g$ defined in~\eqref{eq:phi0}, and is a global margin because it is independent of $x$.

Note that condition~\eqref{eq:cbf_online_condition} imposes that the time derivative of $h$ vanishes as $h$ approaches the boundary of the safe set. In contrast, the ZOH-CBF condition~\eqref{eq:zoh_cbf_condition} causes the time derivative of $h$ to vanish at a manifold in the interior of the safe set. Inspired from this, we define a second metric of comparison, which captures the maximum distance between this manifold and the boundary of the safe set.

\begin{definition}[\textbf{Physical margin}] \label{def:physical_margin}
    For a solution $\phi$ of Problem \ref{problem:design} with the form~\eqref{eq:phi_form}, the \textnormal{physical margin} is the function $\delta:\mathbb R_{>0}\rightarrow\mathbb R$ defined as 
    \begin{equation}
        \delta(T) \triangleq \sup_{\{x\in S \;\mid\; \phi(T,x) = 0\}} -h(x) \,. \label{eq:physical_margin}
    \end{equation}
\end{definition}
\noindent Intuitively, $\delta$ quantifies the effective shrinkage of the safe set due to the error introduced by discrete sampling.
The condition~\eqref{eq:zoh_cbf_condition} may exclude closed-loop trajectories from entering the set $S_\delta = \{x\; |\; -\delta \leq h(x)\leq 0\}$, while the condition~\eqref{eq:cbf_online_condition} does not. 
Thus, a smaller physical margin $\delta$ implies a smaller subset $S_\delta$ of the safe set where system trajectories may not be allowed to enter.
In our paper~\cite{breeden_acc21} we develop three solutions to Problem \ref{problem:design} that have lower controller and/or physical margins than $\phi_0^g$, in both local and global forms, which follow from either continuous-time CBF conditions such as~\eqref{eq:cbf_online_condition}, or discrete-time CBF conditions~\cite{Discrete_CBFs,Blanchini}. In the interest of space, in this tutorial paper, we will refer to only one of the results and interested readers are referred to~\cite{breeden_acc21} for a thorough analysis and comparison among all three methods and their relation to the state-of-the-art.

\subsubsection{A Less Conservative Methodology}
Rather than choosing $\phi$ so as to enforce~\eqref{eq:cbf_online_condition} between sample times, as is done in~\cite{Cortez,usevitch2021adversarial,Singletary}, here we start from a discrete-time CBF condition and apply it to an approximation of the continuous-time dynamics. One sufficient discrete-time CBF condition, as shown in~\cite{Discrete_CBFs}, is
\begin{equation}
	h(x_{k+1}) - h(x_k) \leq -\gamma h(x_k), {\color{black}\; \forall k \in \mathbb{N}} \label{eq:discrete_cbf}
\end{equation}
for some $\gamma \in (0,1]$.
In general, this condition is not control-affine. However, its linear approximation is control-affine and thus amenable to inclusion in a QP. The error of a linear approximation of a twice differentiable function is bounded by the function's second derivative. For brevity, define 
    {\color{black}$\psi(x,u) \triangleq \nabla[\dot{h}(x)]\,(f(x) + g(x)u)$} 
which represents the second derivative of $h$ between time steps.
Since $f,g,\nabla[h]$ are assumed locally Lipschitz, $\psi$ is defined almost everywhere. Let $\mathcal{R}(x,T)$ denote the set of states reachable from some $x(0)\in S$ in times $t\in[0, T)$. {\color{black}Define the bound}
\begin{equation}
	\eta(T,x) \triangleq \max\left\lbrace\left(\sup_{z\in\mathcal{R}(x,T)\setminus\mathcal{Z},u\in U} \psi(z,u)\right), 0\right\rbrace \label{eq:eta} \,,
\end{equation}
where $\mathcal{Z}$ is any set of Lebesgue measure zero {\color{black}(to account for CBFs that are not twice differentiable everywhere)}. A solution to Problem~\ref{problem:design} is then as follows.

\begin{theorem} \label{thm:nu3}
    The function $\phi_3^l:\reals_{>0}\times\reals^n\rightarrow\reals$, defined as
	\begin{equation}
		\phi_3^l(T,x) \triangleq -\frac{\gamma}{T}h(x) - \underbrace{\frac{1}{2}T\eta(T,x)}_{\nu_3^l(T,x)} \, \label{eq:phi3}
	\end{equation}
	solves Problem \ref{problem:design}, for any $\gamma \in (0,1]$.
\end{theorem}

In~\cite{breeden_acc21}, we provide a detailed discussion and proofs on how this method is less conservative as compared to~\eqref{eq:phi0} and to the other methods derived. Here, we only demonstrate this by simulation. 
Finally, if one wishes to instead use a global margin to avoid needing to compute $\mathcal{R}$, the function $\phi_3^g:\reals_{>0}\rightarrow\reals$ as follows also solves Problem~\ref{problem:design}:
\begin{equation}
    \phi_3^g (T,x) \triangleq -\frac{\gamma}{T}h(x) - \underbrace{\frac{1}{2}T\sup_{z\in S}\eta(T,z)}_{\nu_3^g(T)} \,. \label{eq:phi3g}
\end{equation}
See also~\cite{breeden2023aiaa} for extensions of~\eqref{eq:phi3g}.

\subsubsection{Simulation Results}\label{sec:ZOH results}

We present a case study involving a robotic agent modeled as the unicycle system
\begin{equation*}
\dot x_1 = u_1 \cos(x_3), \; \dot x_2 = u_1 \sin (x_3), \; \dot x_3 = u_2,
\end{equation*}
where $[x_1,\,x_2]^T$ is the position, $x_3$ is the orientation, and $u_1$,$u_2$ are the linear and angular velocity of the agent; its task is to move around an obstacle at the origin using the CBF~\cite{polar_cbf}
\begin{equation*}
	h = \rho - \sqrt{x_1^2 + x_2^2 - ({\textstyle\wrap_\pi}(x_3 - \sigma \arctan2(x_2, x_1)))^2} \,,
\end{equation*}
where $\rho$ is the radius to be avoided, and $\sigma$ is a shape parameter. We choose $T = 0.1$, $\alpha(\lambda) = \lambda$ for $\phi_0^g$, and $\gamma = 1$ for $\phi_3^l,\phi_3^g$. Other notable parameters are listed in~\cite[Table~1]{breeden_acc21}.

The agent moves under the following controller
\begin{equation}
	u = \textrm{argmin}_{u\in U} \|u - u_\textrm{nom}\| \textrm{ s.t.~\eqref{eq:zoh_cbf_condition} is satisfied} \label{eq:the_qp}
\end{equation}
where $u_\textrm{nom}$ is a nominal control law that ignores the obstacle.  Note that the local margin~\eqref{eq:phi3} requires $\nu_3^l$ to be computed online (or computed offline and stored in a function); this computation took 0.018 seconds per control cycle on a 3.5 GHz computer using MATLAB R2019b. The global margin~\eqref{eq:phi3g} only requires $\nu_3^g$ to be computed once offline, which took under a minute. In total, 7 solutions ($\phi_0^g$, $\phi_1^g$, $\phi_1^l$, $\phi_2^g$, $\phi_2^l$, $\phi_3^g$, $\phi_3^l$) to Problem \ref{problem:design} were tested; simulation code and details on the other methods $\phi_1^l$, $\phi_1^g$, $\phi_2^l$, $\phi_2^g$ can be found in~\cite{breeden_acc21}.


\begin{table}[t]
    \centering
\begin{tabular}{c| c| c | c }  
$T$ & 0.1 & 0.01 & 0.001 \\
\hline 
$\delta_0^{g}$ & $1.2(10)^{42}$ & 420 & 0.010 \\
$\delta_1^{g}$  & 0.54 & 0.054 & 0.0054 \\
$\delta_2^{g}$ & 0.53 & 0.053 & 0.0053  \\
$\delta_3^{g}$ & 0.013 &  $1.3(10)^{-4}$  &  $1.3(10)^{-6}$ \\
\end{tabular}
    \caption{Global physical margins for selected time-steps $T$}
    \label{tab:physical}
\end{table}

\begin{figure}
	\centering
	\begin{tikzpicture}
        \node[anchor=south west,inner sep=0] (image) at (0,0) {\includegraphics[width=0.6\columnwidth,trim={0.25in, 0in, 0.5in, 0.3in},clip]{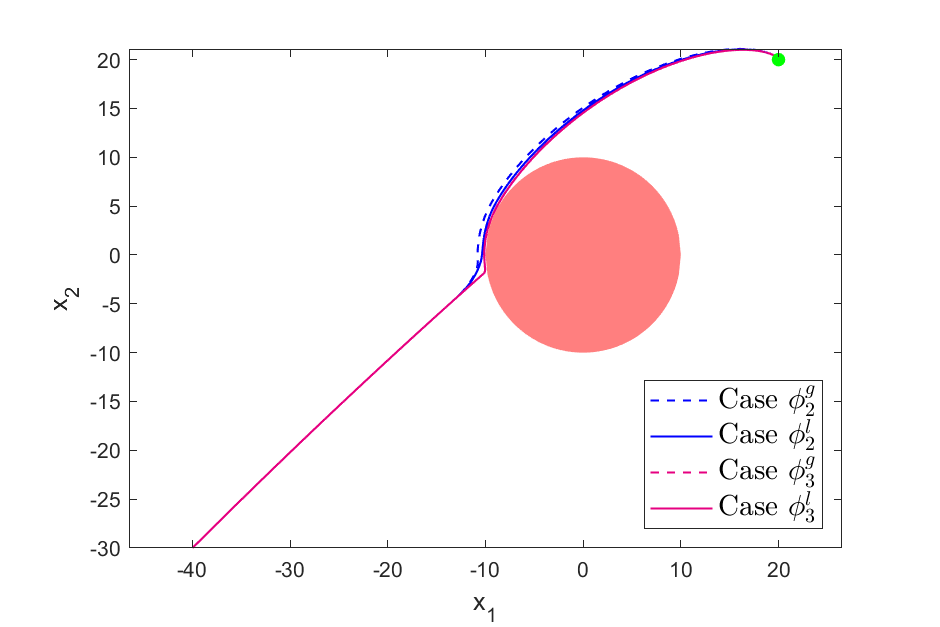}}; %
        \node[anchor=south west,inner sep=0] (image) at (1.1,2.92) 
        {\includegraphics[width=0.94in,trim={3in, 2.1in, 2.5in, 1.3in},clip,frame]{figures/Uni_Trajectories.eps}
        };
    \end{tikzpicture}
	\caption{{The trajectories of the unicycle for 4 of the 7 margin functions}}
	\label{fig:unicycle_trajectories}
\end{figure}

\begin{figure}
\centering
\includegraphics[width=0.7\columnwidth]{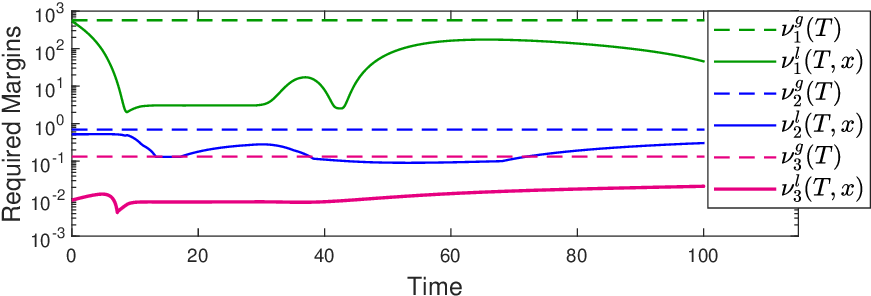}	

\includegraphics[width=0.7\columnwidth]{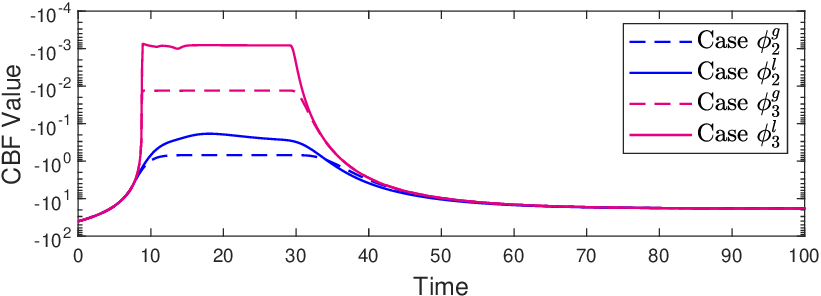}
	\caption{{Top: Controller margins for the unicycle system. Bottom: CBF values along the 4 unicycle trajectories in Fig.~\ref{fig:unicycle_trajectories}}}
	\label{fig:unicycle_barriers}
\end{figure}

The trajectories of the unicycle around an obstacle are plotted in Fig.~\ref{fig:unicycle_trajectories}, where the green marker is the target location. Only four methods are shown because using $\phi_0^g$,$\phi_1^g$,$\phi_1^l$ resulted in the agent turning away from the target due to excessive conservatism. 
The instantaneously required controller margins $\nu$ for all 7 methods, computed for $x(t)$ along the $\phi_3^l$ trajectories from Fig.~\ref{fig:unicycle_trajectories}, are plotted in the top plot in Fig.~\ref{fig:unicycle_barriers}. Note the logarithmic scale; the prior work $\nu_0^g$ is omitted because $\nu_0^g = 1.3(10)^{50}$, indicating that this method cannot be used at the chosen time step of $T = 0.1$.


\begin{figure}
    \centering
    \begin{tikzpicture}
        \node[anchor=south west,inner sep=0] (image) at (0,0) {\includegraphics[width=0.7\columnwidth,trim={0.25in, 0in, 0.5in, 0.34in},clip]{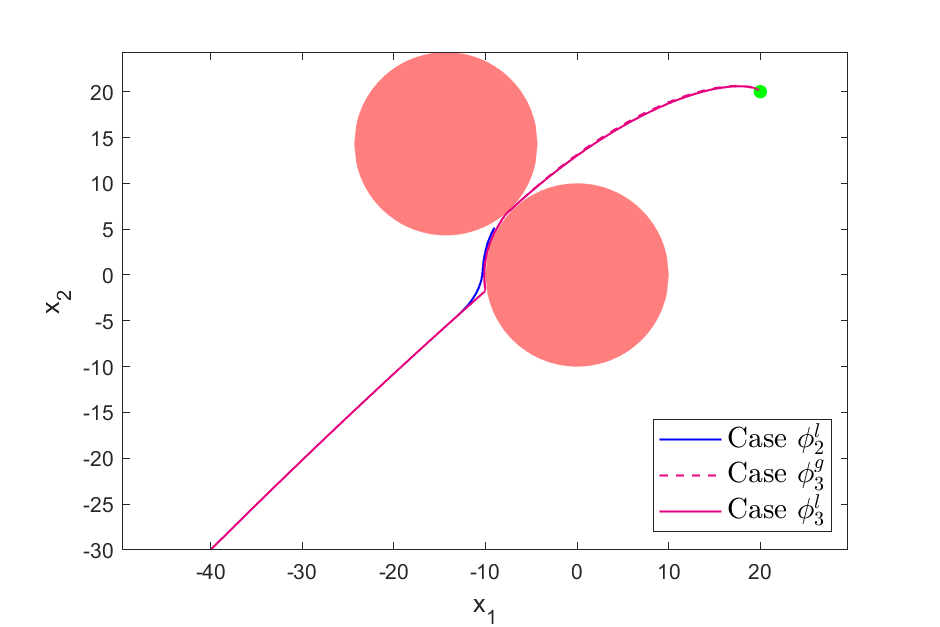}};
        \node[anchor=south west,inner sep=0] (image) at (1.1,2.4) {\includegraphics[width=0.77in,trim={3in, 2in, 2.5in, 1in},clip,frame]
        {figures/Uni_TwoObstacles.eps}
        };
    \end{tikzpicture}
    \caption{{A simulation with two tightly-spaced obstacles, in which controllers using margins $\phi_3^l$ and $\phi_3^g$ permit passage through the obstacles, while the other functions force the agent to stop.}} 
    \label{fig:two_obstacles}
\end{figure}

Next, consider the physical margins. As~\eqref{eq:physical_margin} includes a supremum over $S$, the physical margin is inherently a global quantity. Physical margins for the prior work $\phi_0^g$ and the three new global margins $\phi_1^g$, $\phi_2^g$, and $\phi_3^g$ are computed in Table~\ref{tab:physical}. Note how the physical margin for $\phi_3^g$ decreases quadratically with $T$, whereas the other methods decrease linearly with $T$; the reasoning for this is elaborated upon in~\cite{breeden_acc21}. Finally, the CBF values during the simulations from Fig.~\ref{fig:unicycle_trajectories} are also shown in Fig.~\ref{fig:unicycle_barriers}. The peaks of the dashed lines agree well with the theoretical values in the first column of Table~\ref{tab:physical}. 
%
Noting these physical margins, we added a second constraint that forced the unicycle to navigate through a narrow corridor only 0.3 units wide, shown in Fig.~\ref{fig:two_obstacles}. The unicycle operating under $\phi_3^g$ or $\phi_3^l$ makes it through the obstacles, while the best of the other methods ($\phi_2^l$) could not.

\section{Adversarially-Robust Control Barrier Functions}
\label{section::adversarial}
We present applications of CBFs to the control of multi-agent/multi-robot systems in the presence of adversaries. We first review our results on Adversarially-Robust CBFs, which provide safety conditions for sampled-data distributed control when agents are behaving adversarially. Then we present an application of the method to resilient control against adversarial (chasing) agents.
\subsection{Sampled-Data Distributed Control}
The idea of assuring forward invariance of a set  under sampled-data implementation can be utilized and extended to control of sampled-data multi-agent systems. In our work~\cite{usevitch2021adversarial, usevitch2022tac} we consider a class of functions describing safe sets that have high relative degree with respect to (w.r.t.) the system dynamics, where the control inputs of the agents do not appear for one or more time derivatives of the safe-set function. We also consider asynchronous sampling times with clock disturbances, the presence of adversarially-behaving agents, and functions describing safe sets that have high relative degree w.r.t. the system dynamics. Our goal is to establish a framework under which a set of normally-behaving agents in a system with sampled-data dynamics can collaboratively render a safe set forward invariant despite the actions of adversarial agents. Our analysis considers asychronous sampling times and distributed calculation of agents' control inputs.

More specifically, we consider a group of $N \in \mathbb Z_{>0}$ agents, with the set of agents denoted by $\mathcal V$ and each agent indexed $\{1,\ldots,N\}$.
Each agent $i \in \mathcal V$ has the state $x_i \in \mathbb R^{n_i}$, $n_i \in \mathbb Z_{>0}$ and input $u_i \in \mathbb R^{m_i}$, $m_i \in \mathbb Z_{>0}$.

In addition, at each $t_i^k \in \mathcal{T}_i$ the agent $i$ applies a zero-order hold (ZOH) control input $u(t_i^k)$ that is constant on the time interval $t \in [t_i^k, t_i^{k+1})$.
For brevity, we denote $x_i^{k_i} = x_i(t_i^k)$ and $u_i^{k_i} = u_i(t_i^k)$.
The sampled-data dynamics of each agent $i \in \mathcal V$ under its ZOH controller on each interval $t \in [t_i^{k}, t_i^{k+1})$ is as follows:
\begin{align}
\begin{aligned}
\label{eq:generalsystem}
	\dot{x}_i(t) &= f_i(x_i(t)) + g_i(x_i(t)) u_i(t_i^k) + \phi_i(t). \\
\end{aligned}
\end{align}
The functions $f_i$, $g_i$ may differ among agents, but are all locally Lipschitz on their respective domains $\mathbb R^{n_i}$.
The functions $\phi_i: \mathbb R \rightarrow \mathbb R^{n_i}$, $i \in \mathcal V$,
are locally Lipschitz in $t$ and
model disturbances to the system~\eqref{eq:generalsystem}. 
Each $\phi_i$ is bounded as per the following assumption:
\begin{assumption}
For all $i \in \mathcal V$, the disturbances $\phi_i(t)$ satisfy $\|{\phi_i(t)}\| \leq \phi_{i}^{\max} \in \mathbb R_{\geq 0},\ \forall t \geq 0$.
\end{assumption}

Safety of the multi-agent system can be collaboratively preserved by defining a multi-agent CBF $h$ as
\begin{align}
    \begin{aligned}
\label{eq:setdefinition}
    S = \{\vec{x} \in \mathbb{R}^{\bar{n}} : h(\vec{x}) \leq 0\}, \\
    \partial S = \{\vec{x} \in \mathbb{R}^{\bar{n}} : h(\vec{x}) = 0\}, \\
    \text{int}(S) = \{\vec{x} \in \mathbb{R}^{\bar{n}} : h(\vec{x}) < 0\}.
\end{aligned}
\end{align}
Forward invariance of $h$ can be guaranteed by satisfying the condition 
\begin{align}
\label{eq:Nagumosum}
    \sum_{i \in \mathcal{V}} \pth{L_{f_i} h^{x_i}(\vec{x}) + L_{g_i} h^{x_i}(\vec{x}) u_i + L_{\phi_i} h^{x_i}(\vec{x})} \leq -\alpha(h(\vec{x})),
\end{align}
for a class-$\mathcal{K}_\infty$ function $\alpha$, which follows from a comparison result~\cite{konda2020characterizing}.
Prior literature assumed the cooperation of all agents to collaboratively ensure the satisfaction of the safety condition for forward invariance of $h$. Our work dropped this assumption and considered the presence of a subset of \emph{adversarial agents} $\mathcal{A} \subset \mathcal{V}$ that apply the following control input for all sampling times $t_j^k$, $k \in \mathbb{Z}_{\geq 0}$, $j \in \mathcal{A}$:
\begin{align}
\label{eq:adversarialu}
    \hspace{-4pt}u_j^{\max}(\vec{x}^{k_j}) = \arg\max_{u \in \mathcal{U}_j} \bkt{L_{f_j}h^{x_j}(\vec{x}^{k_j}) + L_{g_j} h^{x_j}(\vec{x}^{k_j}) u}.
\end{align}
Agents that are not adversarial are called \emph{normal}. The set of normal agents is denoted $\mathcal{N} = \mathcal{V} \backslash \mathcal{A}$.
Dividing the left-hand side (LHS) of~\eqref{eq:Nagumosum} into normal and adversarial parts yields the following sufficient condition for set invariance in the presence of adversaries:
\begin{align}
    &\sum_{j \in \mathcal{A}} \pth{L_{f_j} h^{x_j}(\vec{x}) + L_{g_j} h^{x_j}(\vec{x}) u_j^{\max} + L_{\phi_j} h^{x_j}(\vec{x})} + \label{eq:Nagumo_adversarial}\\
    &\sum_{i \in \mathcal{N}} \pth{L_{f_i} h^{x_i}(\vec{x}) + L_{g_i} h^{x_i}(\vec{x}) u_i + L_{\phi_i} h^{x_i}(\vec{x})} \leq -\alpha(h(\vec{x})). \nonumber
\end{align}
The form of~\eqref{eq:Nagumo_adversarial} reflects sampled-data adversarial agents seeking to violate the set invariance condition in~\eqref{eq:Nagumosum} by maximizing their individual contributions to the LHS sum.
The normal agents must compute control inputs that render the set $S$ forward invariant using the sufficient condition in~\eqref{eq:Nagumo_adversarial} despite the worst-case behavior of the adversarial agents in $\mathcal{A}$.

In addition to the ZOH sampled-data dynamics and adversarial actions, the normal agents must take into account asynchronous sampling times.
The assumption of identical, synchronous sampling times typically does not hold in practice.
Also, a distributed system may not have access to a centralized entity to solve a QP computing control inputs for all normal agents.
Finally, safety may be defined in terms of CBFs which may have high relative degree with respect to the agents' dynamics.
It is therefore necessary to consider heterogeneous sampling times, distributed methods for computing local control inputs, and CBFs with high relative degree.


In consideration of these challenges, our work~\cite{usevitch2022tac} defined the following safety-preserving control set for each normal agent $i \in \mathcal{N}$:
\begin{align}
    \label{eq:Kiki}
    K_i^{\psi}(\vec{x}^{k_i}) &= \brc{u_i \in \mathcal{U}_i : \psi_q^d(\vec{x}^{k_i}) \leq 0},\\
    &\hspace{1em} = \Big\{u_i \in \mathcal{U}_i : L_{f_i} (\psi_{q-1}^d)^{x_i}(\vec{x}^{k_i}) + L_{g_i} (\psi_{q-1}^d)^{x_i}(\vec{x}^{k_i}) u_i \\
    &\hspace{1em} +\sum_{l \in \N \backslash \{i\}} \bkt{L_{f_l} (\psi_{q-1}^d)^{x_l}(\vec{x}^{k_i}) + L_{g_l} (\psi_{q-1}^d)^{x_l}(\vec{x}^{k_i}) \hat{u}_l^{k_i}} \\
     &\hspace{1em} + \sum_{j \in \A} \gamma_j^{\max}(\vec{x}^{k_i}) + \alpha_q(\psi_{q-1}^d(\vec{x}^{k_i}))  + \eta'(\Gamma_i + \delta^{\max}) \leq 0 \Big\}.
\end{align}
In this equation, the functions $(\psi_{q-1}^d)^{x_i}$ for each agent $i \in \N$ each represent the last of a series of functions $\psi_{i+1} \triangleq \dot{\psi}_i(\vec{x}) + \alpha_{i+1}(\psi_i(\vec{x})$ typically defined to account for CBFs with high relative degree~\cite{xiao2019control}.
The term $\gamma_j^{\max}$ is a function describing the worst-case adversarial behavior with respect to the CBF safety condition:
\begin{align}
    \gamma_i^{\max}(\vec{x}) &= \max_{u_i \in \mathcal{U}_i} \bkt{L_{f_i}h^{x_i}(\vec{x}) + L_{g_i}h^{x_i}(\vec{x}) u_i}. \label{eq:vLPs}
\end{align}
\begin{figure}
    \centering
    \includegraphics{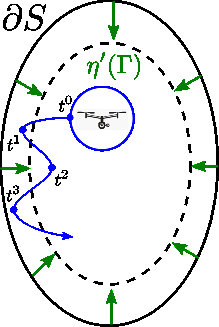}
    \caption{An illustration of the effect of the $\eta$ term in the definition of the safety-preserving safe set. The $\eta$ function ensures that safety is maintained for all times despite the sampled-data dynamics and heterogeneous sampling / communication times. This is accomplished by contracting the original safe set. The amount of contraction is proportional to several factors including the sampling period, the agents' $f$ and $g$ dynamics functions, and an upper bound on temporal disturbances to sampling times.}
    \label{fig:enter-label}
\end{figure}
The function $\eta'$ in~\eqref{eq:Kiki} accounts for the evolution of system dynamics between sampling times and heterogeneous sampling instances between agents, and is defined as
\begin{align}
\label{eq:etaprimedef}
    \eta'(\Gamma) = \pth{c_f' + c_g' u_{\max} + c_\alpha' + c_{\hat{\gamma}}'} \eps(\Gamma).
\end{align}
The variables $\Gamma_i$ and $\delta^{\max}$ in~\eqref{eq:Kiki} account for heterogeneous nominal sampling periods and an upper bound on disturbances to the nominal sampling times, respectively.

For each normal agent $i \in \mathcal{N}$, applying a control input $u_i^{k_i} \in K_i)\vec{x}^{k_i}$ guarantees that the trajectory of the combined normal agents' states remains within the safe set for all times within the sampling interval $t \in [t_i^k, t_i^{k+1})$: 

\begin{theorem}
\label{thm:async}
Consider the system~\eqref{eq:generalsystem} with sampling times described by 
$\mathcal{T}_i = \{t_i^0, t_i^1, \ldots\} \text{ s.t. } t_i^{k+1} - t_i^{k} = \Gamma_i + \delta_i(k),\ \forall k \in \mathcal{Z}_{\geq 0},$.
If at sampling time $t_i^k$ for $k \geq 0$, $i \in \mathcal{N}$ it holds that $\vec{x}^{k_i} \in S$, then for any $u_i^{k_i} \in K_i(\vec{x}^{k_i})$ the trajectory $\vec{x}(t)$ satisfies $\vec{x}(t) \in S$ for all $t \in [t_i^k, t_i^{k+1})$.
\end{theorem}

The safe control inputs $u_i^{k_i} \in K_i)\vec{x}^{k_i}$ can be computed for each normal agent $i$ in a distributed manner using a convex quadratic programming formulation:
\begin{alignat}{3}
		&u_i(\vec{x}^{k_i}) = 
		&&\underset{u_i \in \mathcal{U}_i}{\arg\min} \|u_i - u_{i,\textup{nom}}^{k_i}\|_2^2 \label{eq:QPdistributed} \\
		& \hspace{1em}\text{s.t.} & &   \pth{L_{f_i}h^{x_i}(\vec{x}^{k_i}) + L_{g_i}h^{x_i}(\vec{x}^{k_i})u_i} + \nonumber\\
		& & & \sum_{l \in \mathcal{N} \backslash \{i\}} \pth{L_{f_l}h^{x_l}(\vec{x}^{k_i}) + L_{g_l}h^{x_l}(\vec{x}^{k_i})\hat{u}_l^{k_i}} + \nonumber\\
		& & & \hspace{2em} \sum_{j \in \mathcal{A}} \gamma_j^{\max}(\vec{x}^{k_i}) + \alpha(h(\vec{x}^{k_i})) + \nonumber\\
		& & & \hspace{2em}\eta(\Gamma_i + \delta^{\max}) \leq 0. \nonumber
\end{alignat}

Further details about these results can be found in~\cite{usevitch2022tac}.

\begin{figure*}
    \centering
    \includegraphics[width=0.32\textwidth]{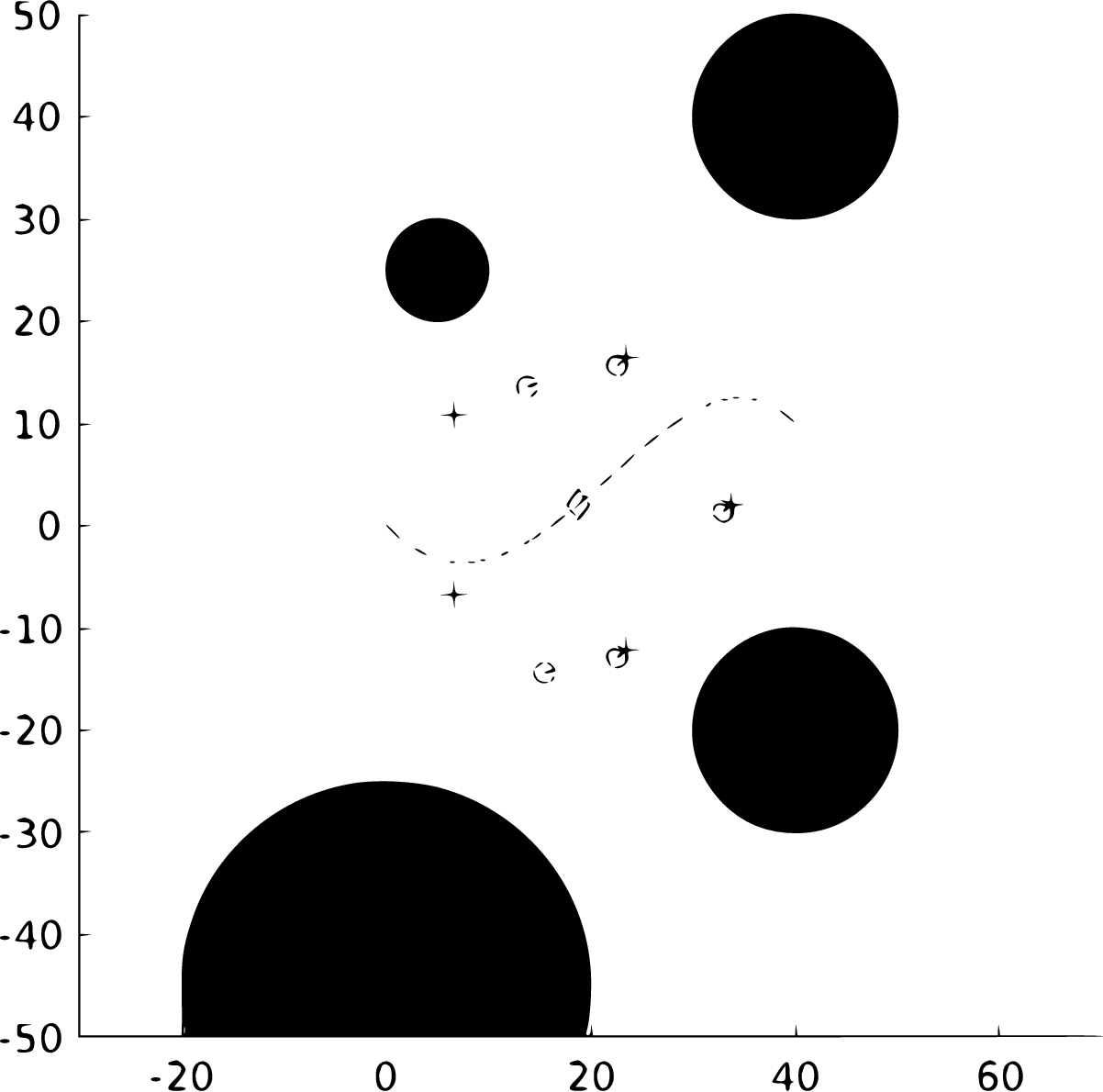} 
    \includegraphics[width=0.32\textwidth]{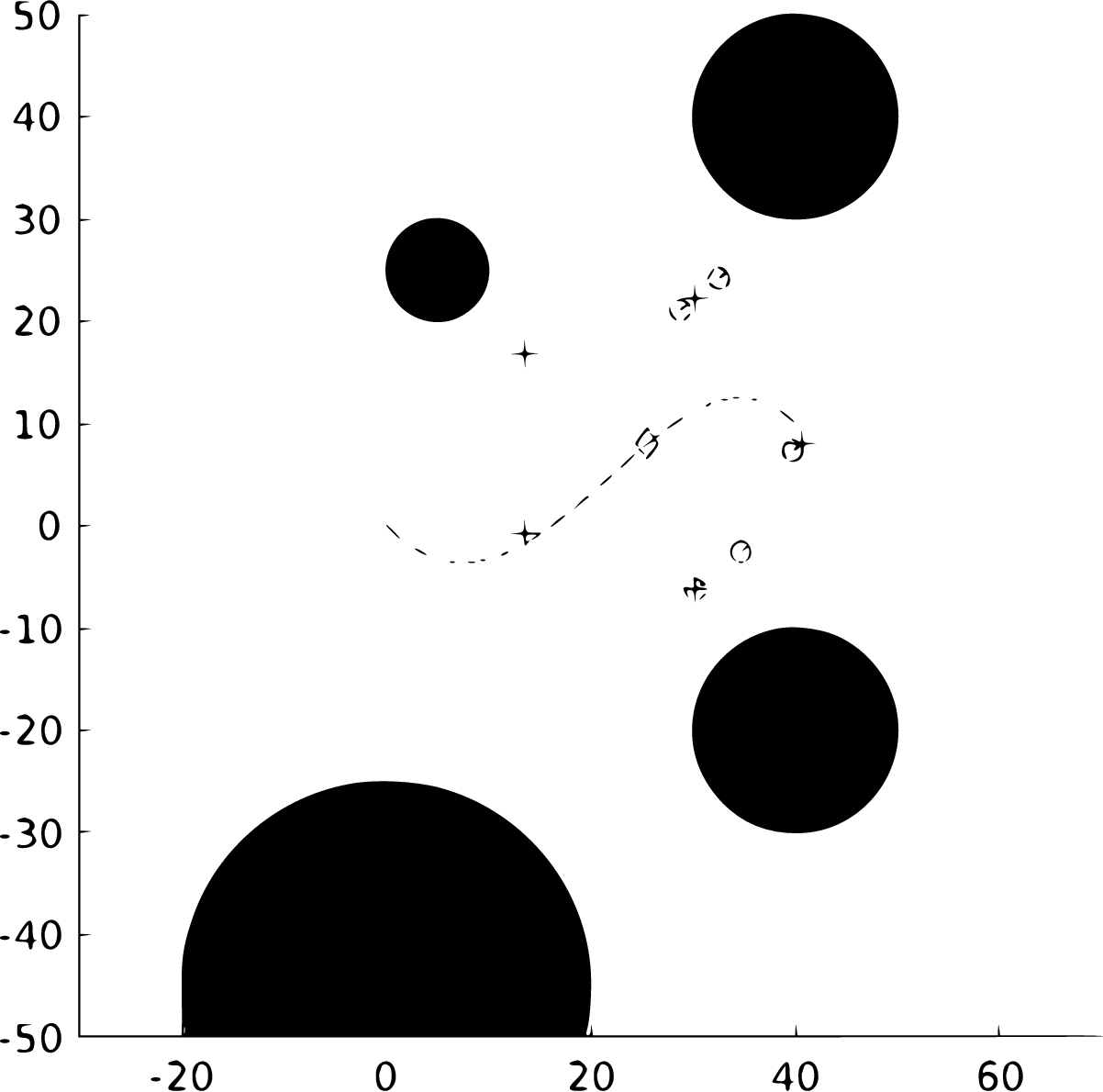} 
    \includegraphics[width=0.32\textwidth]{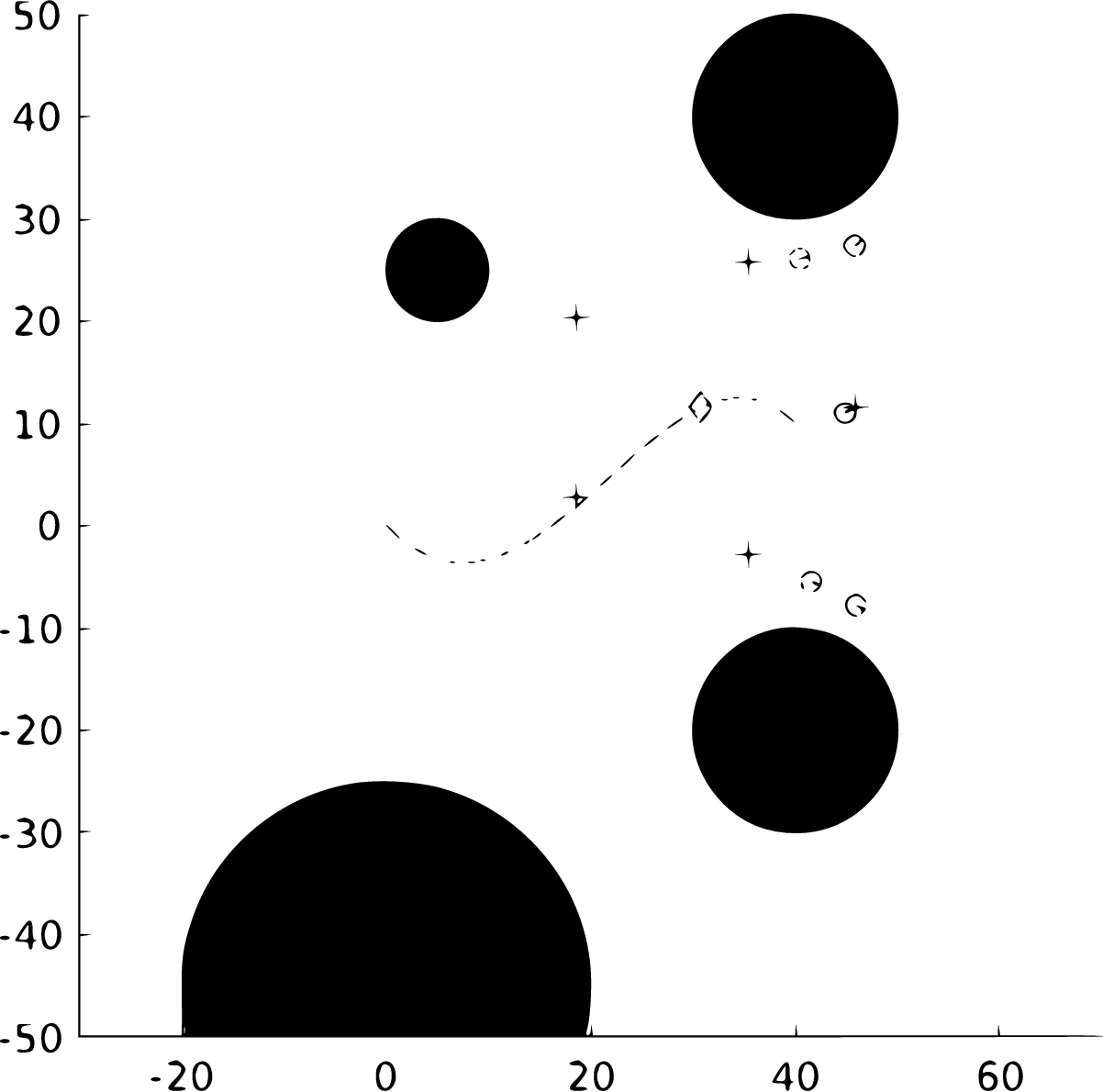} 
    
    
    \caption{Still frames from a simulation involving adversarially resilient CBFs. Normal agents are represented by blue circles and adversarial agents are represented by red circles. The dotted red lines around the blue circles represent normal agents' safety radii. The time-varying formation trajectory is represented by the dotted magenta line; the magenta diamond represents the center of formation. Black crosses represent agents' nominal local time-varying formational points.}
    \label{fig:Sim1}
\end{figure*}


\section{Non-Smooth Control}
\label{section::stronginvariance}

Many forward invariance results in prior literature on CBFs assume that the agent's control input is Lipschitz continuous. At the same time, the majority of CBF results rely upon computing safe control inputs via a parametric convex quadratic program (QP). Demonstrating that the Lipschitz continuity of a parametric QP is nontrivial in general.

Our work in~\cite{usevitch2020strong} addresses a simple question: Can forward invariance be assured for discontinuous control inputs?
Our results in~\cite{usevitch2020strong} answer this in the affirmative.
Unlike prior literature, our work approaches the problem using the notions of Clarke tangent cones and transversality. We demonstrate that a constrained control input simultaneously rendering these subsets invariant can be generated by simply solving a feasibility problem with compact linear constraints. The control input is only required to be Lebesgue measurable and is not required to be continuous. 

This work considers control affine systems in the form
\begin{align}
\label{eq:controlaffine}
\begin{aligned}
    \dot{x}(t) &= f(x(t)) + g(x(t))u(t), \\
    u(t) &\in \mathcal{U} \subset \mathbb{R}^m\ \forall t \geq t_0.
\end{aligned}
\end{align}
where the functions $f: \mathbb{R}^n \rightarrow \mathbb{R}^n$ and $g : \mathbb{R}^n \rightarrow \mathbb{R}^{n \times m}$ are assumed to be locally Lipschitz on $\mathbb{R}^n$ and the set $\mathcal{U}$ is a compact, convex polytope with $\textnormal{int}(\mathcal{U}) \neq \emptyset$ which has the form
\begin{align}
\begin{aligned}
\label{eq:upolytope}
\mathcal{U} &= \{u \in \mathbb{R}^m : A_u u \leq b_u \}, \\ A_u &\in \mathbb{R}^{p \times m},\ b_u \in \mathbb{R}^{p \times 1},
\end{aligned}
\end{align}
where $A_u$, $b_u$ are constant. Our analysis requires the notion of \emph{strict} CBFs, which are defined as follows:
\begin{definition}
The continuously differentiable function $h : \R^n \rarr \R$ is called a \textbf{strict} CBF for the set $S \subset \R^n$ defined as $S = \{x\; |\; h(x)\leq 0\}$ if the following holds:
\begin{align}
\label{eq:necessary}
    \inf_{u \in \mathcal{U}} \bkt{L_f h(x) + L_g h(x) u} < 0\ \forall x \in \partial S,
\end{align}
where $f,g$ are defined as in~\eqref{eq:controlaffine}.
\end{definition}

Our work approaches the problem by designing a differential inclusion of the form
\begin{align}
    G(x) = \{f(x) + g(x)u : u \in K(x) \},
\end{align}
where the set-valued map $K : \R^n \rarr \mathcal{P}(\R^m)$ satisfies $K(x) \subseteq \mathcal{U}$ for all $x \in \R^n$. 
When a single safe set is being considered, $K$ is defined as
\begin{align}
\begin{aligned}
\label{eq:Kforreal}
        &K(x) = \brc{u \in \R^m : \bmx{A_S(x) \\ A_u} u \leq \bmx{b_S(x) \\ b_u}},
\end{aligned}
\end{align}
To ensure that the set-valued maps are locally Lipschitz, we consider a \emph{$\gamma$-contraction} of $K$ and $G$ defined as
\begin{align}
    \label{eq:Kgamma}
        K_{\gamma}(x) &= \textnormal{int}(K(x)) - \gamma B(0,1), \\
        &= \{u \in K(x) : d_{K^c}(u) \geq \gamma \},\ K^c = \R^m \backslash K(x), \\
        G_{\gamma}(x) &= \brc{v \in \R^n :  v = f(x) + g(x) u,\  u \in K_{\gamma}(x)}.
    \end{align}
The following result demonstrates conditions under which forward invariance of a safe set under a discontinuous control input can be guaranteed.

\begin{theorem}
\label{thm:main1}
    Consider the system 
    \begin{align} 
    \label{eq:theoremsystem}
    \dot{x}(t) \in G_{\gamma}(x(t)).
    \end{align}
    Let $S$ be a safe set for some strict control barrier function $h$.
    Let $x(\cdot)$ be any trajectory of~\eqref{eq:theoremsystem} under a Lebesgue measurable control input $u(\cdot)$ with $x_0 = x(0) \in \textnormal{int}(S \cap \Omega)$. Let $[0,T(x_0))$ be the (possibly empty) maximal interval such that $x(t) \in \textnormal{int}(\Omega)$ for all $t \in [0,T(x_0))$. Then $x(t) \in S$ for all $t \in [0,T(x_0))$.
\end{theorem}

It is common to consider multiple safe sets $S_{i} = \{x\; |\; h_i(x)\leq 0\}$ simultaneously. In other words, we seek to render the composed set $S_I = \bigcap_{i=1}^{N_h} S_i$ strongly invariant.
Towards this end we define the set-valued map
\begin{align}
\begin{aligned}
\label{eq:multiKforreal}
        &\widehat{K}(x) = \brc{u \in \R^m : \bmx{\hat{A}_S(x) \\ A_u} u \leq \bmx{\hat{b}_S(x) \\ b_u}}, \\
        &\hat{A}_S : \R^n \rarr \R^{q \times m},\ \hat{b}_S : \R^n \rarr \R^{q}.
\end{aligned}
\end{align}
where $\hat{A}_S, \hat{b}_S$ are defined as
\begin{align}
    \hat{A}_S(x) &= \bmx{L_g h_1(x) \\ \vdots \\ L_g h_{N_h}(x)}, \nonumber \\
    \hat{b}_S(x) &= \bmx{-\alpha_1(h_1(x)) - L_f h_1(x) \\ \vdots \\ -\alpha_{N_h}(h_{N_h}(x)) - L_f h_{N_h}(x)}. \label{eq:ABmultiple}
\end{align}
We similarly define the $\hat{\gamma}$-contractions $\widehat{K}_{\hat{\gamma}}(x) = \textnormal{int}(\widehat{K}(x)) - \hat{\gamma} B(0,1)$ and $\widehat{G}_{\hat{\gamma}}(x) = \brc{v \in \R^n :  v = f(x) + g(x) u,\  u \in \widehat{K}_{\hat{\gamma}}(x)}$. The following result establishes conditions under which forward invariance of the composed set $S_I$ holds under possibly discontinuous control inputs.

\begin{theorem}
    \label{thm:main2}
    Consider the system 
    \begin{align} 
    \label{eq:theoremsystem2}
    \dot{x}(t) \in \widehat{G}_{\hat{\gamma}}(x(t)).
    \end{align}
    Consider the set $S_I = \bigcap_{i = 1}^{N_h} S_i$ and suppose that the transversality condition holds for the pair $(S_i, S_j)$ for all $i, j\in I(x)$.
    Let $x(\cdot)$ be any trajectory of~\eqref{eq:theoremsystem2} under a Lebesgue measurable control input $u(\cdot)$ with $x_0 = x(0) \in \textnormal{int}(S_I \cap \widehat{\Omega})$. Let $[0,T(x_0))$ be the (possibly empty) maximal interval such that $x(t) \in \textnormal{int}(\widehat{\Omega})$ for all $t \in [0,T(x_0))$. Then $x(t) \in S_I$ for all $t \in [0,T(x_0))$.
\end{theorem}
A full exposition of the details can be found in~\cite{usevitch2020strong}.

\section{Recent Extensions and Ongoing/Future Work}
\label{section::conclusions}
This tutorial paper presented only a small fraction of the rich literature on the theory and applications of Control Barrier Functions, focusing on challenges such as safety under spatiotemporal and input constraints, safety constraints of higher relative degree, robustness to disturbances/noise, performance via prediction and online parameter adaptation, safety under sampled-data implementation and output feedback control, safety against adversarial inputs, and safety under non-smooth inputs. 

Many open problems remain, including how to scale the results to multi-agent systems without adding conservatism, how to online adapt to unknown nonlinear disturbances, and how to systematically define CBFs that guide certifiable safe learning and exploration. Notable recent extensions of the authors' work include: how adversarially-robust CBFs can be used for the detection and mitigation of adversarial agent effects~\cite{mustafa2022adversary}; how CBFs can be used in case of actuator failures and/or cyber attacks on actuators~\cite{garg2022control,garg2023model}; how CBFs ensure safety when an autonomous agent is exploring an unknown environment to maximize clarity (or minimize uncertainty) about the environment~\cite{agrawal2023cdc}; how one can define CBFs and CLFs for systems with impulsive actuators and dwell time constraints~\cite{breedencdc2023}; how to synthesize safe controllers using Koopman-based identification of nonlinear models~\cite{black2023ifac}; how to incorporate risk awareness for the safety verification of stochastic systems~\cite{black2023icra}; how to define multi-rate architectures for safety verification across the planning and control layers of differentially-flat systems~\cite{agrawal2022lcss}; and how learning-based methods can be used for safety of large-scale MAS using distributed CBFs~\cite{zhang2023neural}. Ongoing work includes the investigation of methods for computing safety certificates/CBFs online with provable guarantees, of methods that enable adaptable resilience against adversaries, and of methods that guide the exploration of multiple robots in unknown constrained environments.    

\section{Acknowledgements}
The authors would like to acknowledge the support of the (1) Air Force Office of Scientific Research (AFOSR) under award number FA9550-17-1-0284; (2) National Science Foundation (NSF) under the award numbers 1942907, 1931982 and the Graduate Research Fellowship Program; (3) Office of Naval Research (ONR) under grant number N00014-20-1-2395. The views and conclusions contained herein are those of the authors only and should not be interpreted as representing those of ONR, the U.S. Navy or the U.S. Government; (4) Automotive Research Center (ARC) in accordance with Cooperative Agreement W56HZV-19-2-0001 U.S. Army CCDC Ground Vehicle Systems Center (GVSC) Warren, MI; (6) Fran\c{c}ois Xavier Bagnoud Foundation.

\bibliographystyle{elsarticle-num} 
\bibliography{journal1,sources}





\end{document}